\documentclass[11pt,a4paper]{amsart}
\usepackage[dvips]{epsfig}
\usepackage{graphics}
\usepackage{latexsym}
\usepackage{verbatim}
\usepackage{amsmath}
\usepackage{amsthm}
\usepackage{mathtools}
\usepackage{mathrsfs}
\usepackage{amssymb}
\usepackage{MnSymbol}
\usepackage{stmaryrd}
\usepackage{hyperref}
\usepackage[]{hyperref}
\hypersetup{
    colorlinks=true,%
    filecolor=black,%
    linkcolor=black,%
    urlcolor=black,  %
    citecolor = red
}
\usepackage [table]{xcolor}
\usepackage{multirow}
\usepackage{float}
\usepackage{tikz}
\usepackage{subfig}
\usepackage{algorithm2e}

\graphicspath{{pix/}}	
\newtheorem{theorem}{Theorem}[section]

\newtheorem{proposition}[theorem]{Proposition}

\newtheorem{remark}[theorem]{Remark}

\newcommand\CC{\hbox{C\kern -.58em {\raise .54ex \hbox{$\scriptscriptstyle |$}}
  \kern-.55em {\raise .53ex \hbox{$\scriptscriptstyle |$}} }}
\newcommand\NN{\hbox{I\kern-.2em\hbox{N}}}

\newcommand\ZZ{{{\rm Z}\kern-.28em{\rm Z}}}

\newcommand{\cfl}{\nu^n}
\newcommand{\ds}{\displaystyle}

\usepackage{algorithmic}

\begin{document}

\title[A hybrid finite volume method for the advection equations]{A hybrid finite volume method for advection equations and its applications in population dynamics}

\author{Chang Yang}

\address{Department of Mathematics, Harbin Institute of Technology, 92 West Dazhi Street, Nan Gang District,
Harbin 150001, China}

\email[C.~Yang]{yangchang@hit.edu.cn}

\author{ L\'eon Matar Tine}

\address{
Universit\'e de Lyon, Universit\'e Lyon 1, CNRS UMR 5208, Institut Camille Jordan 43 blvd du 11 novembre 1918, F-69622 Villeurbanne-Cedex, France.
\&
Inria Team Dracula, Inria Center Grenoble Rhone-Alpes, France.
}

\email[L.M.~Tine]{ leon-matar.tine@univ-lyon1.fr }

\hyphenation{bounda-ry rea-so-na-ble be-ha-vior pro-per-ties
cha-rac-te-ris-tic}

\maketitle

\begin{abstract}
We present in this paper a very adapted finite volume numerical scheme for transport type-equation. The scheme is an hybrid one combining an 
anti-dissipative method with down-winding approach for the flux \cite{bibGLT, bibdespres} and an high accurate method as the WENO5 one \cite{bibJS}. The main 
goal is to construct a scheme able to capture in exact way the numerical solution of transport type-equation without artifact like numerical diffusion or 
without ``stairs" like oscillations and this for any regular or discontinuous initial distribution.
This kind of numerical hybrid scheme is very suitable when properties on the long term asymptotic behavior of the solution are of central importance in the modeling what 
is often the case in context of population dynamics where the final distribution of the considered population and its mass preservation relation are required for prediction.

\end{abstract}

\vspace{0.1cm}

\noindent 
{\small\sc Keywords.}  {\small Discontinuity detector; WENO scheme; anti-diffusive method; population dynamics.}


\tableofcontents


\section{Introduction} 
\setcounter{equation}{0}
\label{sec:Intro}
In this paper we are interested in finite volume numerical simulations of PDEs of hyperbolic type with transport term. More precisely we are looking to correct the numerical diffusion which appears in the simulations of the asymptotic profile of such PDEs. Indeed this numerical diffusion is an artifact that is inherent to most of existing numerical schemes even for high-order ones. So, as reported first in \cite{bibCG} and confirmed in \cite{bibGLT} for the Lifshitz-Slyozov equation which is of transport type, capturing numerically the exact asymptotic profile for transport equation is a real challenge because numerical diffusion smooths out the fronts and leads to an artificial profile. In the context of biology modelling, specially in population dynamics, most of existing models \cite{bibP, bibGreer, bibCush} contain at least a transport part which takes into account the growth of the considered population, it is crucial to recover the exact asymptotic profile in order to predict the behavior of the population or to estimate some parameters for instance its growth, division or death rates by using measures that are based on the profile.            
Some authors address the question consisting to correct this inherent numerical diffusion by establishing adequate schemes such as the WENO (Weighted Essentially Non-Oscillatory) scheme \cite{bibJS}, anti-dissipative scheme the ADM (Anti Dissipative Method) \cite{bibGLT}, etc. All these schemes define their numerical fluxes in order to minimize at best the artifacts.

We propose here a hybrid finite volume  scheme where the construction of the numerical flux appears as a combination of the WENO and ADM fluxes. This scheme is very suitable for firstly removing de numerical artifacts but also correct the stair treads appearing in ADM scheme. More precisely we use a discontinuity detector at each grid point and use the WENO order 5 scheme when the solution is regular near the point otherwise we apply the ADM. As a validation of this hybrid scheme we use two kind of test cases. The first one is a classical (academic) test case for transport equation where the initial distribution is considered in one hand as a very oscillatory one as given in \cite{bibJS}; in the other hand, we use in 2D the famous Zalesak' disk test that is given in \cite{bibLeveque}. The second kind of test is based on population dynamics and polymerization process where we consider a population of cells or polymers growing either by nutrients uptake (for cells) or by gain and lost of monomers (for polymers). This application on population dynamics is of great importance because in many cases some predictions on the numerical behavior of the solution allow to investigate inverse problem of estimating relevant parameters of the considered model. So, having a bad numerical reconstruction induces bad parameters estimation.   

The paper is organized as follows. In section \ref{sec:adm0}, we recall the biology context on which we focus our study. In the third section, we detail the derivation of our hybrid scheme which is based on a general conservation laws. The section \ref{sec:detector} is devoted to the numerical results and comparison of our hybrid method against the WENO and ADM schemes. \\    
      
\section{Biological models}
\label{sec:adm0}
\setcounter{equation}{0}
In cell biology as in physics of particles, the evolution dynamics of a group of cells or macro-particles in cell culture or in a bath of micro-particles plays a central role in the understanding and the explanation of some physical and biological behaviors. Often the observed quantities evolve by growth process either by nutrients uptake (the example of the micro-organism {\it Daphnia} which uses the nutrients to grow \cite[sect. 4.3.1]{bibP}, \cite{bibMD, bibRoos}) or by earnings micro-particles by polymerization (for polymers modeling \cite{bibLS}).

Lot of conjectures are based on the observation of these quantities, especially on their evolution dynamics in long term. In modeling point of view this long term evolution dynamic is obtained by the analysis of the asymptotic behavior of the considered quantity.
 
Following the processes taken into account in the model, the asymptotic behavior can either be dependent or be independent of the initial distribution of the considered group of cells or parasites. Indeed, in the case where the considered population evolves only by growth it is proven in Lifshitz-Slyozov equations that the asymptotic behavior depends on the initial distribution \cite{bibCG, bibGLT}. However, when aggregation process or division process is taken into account, the asymptotic behavior is regularized towards a quasi-universal profile as shown in \cite{bibGLT} and then it is independent to the initial distribution.\\ 
A crucial point linked to the modelling of these phenomena is the numerical simulation which can lead, following the used scheme, to bad conjectures on the behavior of the model. These bad conjectures result to some numerical artifacts caused by numerical diffusion inherent to some standard schemes \cite{bibGLT}.

In order to apply the hybrid method proposed later in this paper, we consider the following test case where a size-structured cells population model is taken into account and the cells evolve by gain or loss of micro cells. Let denote by $f(t,x,\xi)$ the size density repartition of cells of size $\xi\geq0$, located at position $x\in\Omega \subset\mathbb{R}^{i}$ ($i=1, 2, 3$) at time $t\geq0$ where $\Omega$ is a smooth bounded domain. Then the model can be written for all $(t,x,\xi)\in\mathbb{R}_{+}\times\Omega\times\mathbb{R}_{+}$ as follows

\begin{equation}\label{eq1}
\left\{ \begin{array}{lll}
\displaystyle\frac{\partial}{\partial t} f(t,x,\xi)+\frac{\partial}{\partial \xi} ((a(\xi)c(t,x)-b)f(t,x,\xi)) =0,\\
\\
f(t=0, x, \xi)=f^{0}(x,\xi).\\
\end{array} \right.
\end{equation}        
Where $c(t,x)$ is the concentration of micro-organisms (nutrients) and it follows a diffusion equation of this form:
\begin{equation}\label{eq2}
\left\{ \begin{array}{lll}
\displaystyle\frac{\partial}{\partial t} \bigl(c(t,x)+\int_{0}^{\infty}\xi\,f(t,x,\xi) d\xi\bigr) = \Delta_x c(t,x), \quad t\geq0,\; x\in\Omega \; ,\\
\\
\displaystyle\frac{\partial}{\partial \nu}c=\nabla c\,\cdot\nu =0, \quad \text{on}\;  \partial \Omega.\\
\end{array}  \right.
\end{equation}  
For this kind of coupling model \eqref{eq1}-\eqref{eq2}, the kinetic coefficients $a(\xi)$, $b$ are interpreted as the rates at which cells gain or loss nutriments (monomers, micro-organisms).

We assume the micro-cells (or monomers) to follow a diffusion equation as depicted in equation \eqref{eq2}. We endowed this diffusion equation with a homogeneous Neumann boundary condition where $\nu$ is the outward unit vector at point $x\in\partial\Omega$.\\

The problem \eqref{eq1}-\eqref{eq2} is a variant of the very known standard Lifshitz-Slyozov system which models the evolution of a population of macro-particles immersed in a bath of monomers \cite{bibLS}. 

The analytical study of \eqref{eq1}-\eqref{eq2} concerning the existence, uniqueness and properties of the solution is rigorously done in \cite{bibTLG} and the main result is based on the following hypothesis:\\

{\bf Hypothesis.} \; The kinetic coefficients $a, b$ are required to satisfy
$b=1$; $a$ is an increasing function with $a(0)=0$ and $a(+\infty)=+\infty$; $a\in C^{0}([0,\infty))\cap C^{1}((0,\infty))$ and for any $\xi_0>0$ there exists $L_{a,0}>0$ such that $0\leq a'(\xi)\leq L_{a,0}$ for $\xi\geq \xi_0 >0$.\\ 
The initial condition satisfy $c(t=0, x) \in L^{\infty}(\Omega)$; $f(t=0,x,\xi)\in L^{\infty}(\Omega; L^1((0,\infty), (1+\xi) d\xi)))$. \\

With this previous hypothesis, the authors in \cite{bibTLG} prove the following statement on the well-posedness of  \eqref{eq1}-\eqref{eq2}:
\begin{theorem}
There exists a weak solution $(c,f)$ of \eqref{eq1}-\eqref{eq2} with, for any $0<T<\infty$,\\
 $c\in L^{\infty}((0,T)\times\Omega) \, \cap \, L^{2}(0, T; H^{1}(\Omega))$, $f\in L^{\infty}((0,T)\times\Omega;L^{1}((0, \infty),(1+\xi) d\xi))$, \,$c\in C^{0}([0,T]; L^2 (\Omega)-weak)$,\, $f\in C^{0}([0,T]; L^{1}(\Omega\times(0,\infty))-weak)$.   
\end{theorem}
In addition they prove thanks to the Neumann boundary condition, the following mass preservation relation:
\begin{equation}\label{eq:mass_conservation}
\frac{d}{dt} \biggl[\int_{\Omega}\int_{0}^{\infty} \xi f(t,x,\xi) \,d\xi \,dx +\int_{\Omega} c(t,x) \, dx \biggr]=0. 
\end{equation}            

The fact that the space variable $x$ acts as a parameter in the size density repartition $f$ implies that \eqref{eq1} is a transport equation and the study of its asymptotic behavior is numerically very challenging.

Indeed, following the chosen model as in \eqref{eq1}--\eqref{eq2},  one needs in the modelling and simulations to recover the evolution dynamics  of the considered population such as the time evolution of the total number of individuals (even in asymptotic time), the total mass of the population or the conservation law fulfilled by the model. For the numerical simulations of the evolution dynamics, a very adapted scheme is required in order to capture in exact way the solution of the system without artifact in order to get the right and essential properties. That's the aim to introduce the following hybrid method.  


\section{A hybrid finite volume method for advection equations}
\label{sec:adm}
\setcounter{equation}{0}

\subsection{Anti-dissipative method}
\label{sec:Method}
In this section, we consider the following  advection equation
\begin{equation}
   \frac{\partial f}{\partial t}\,\,+\,\,\,  \frac{\partial(V f)}{\partial x}\,\,=\,\,0,\quad t\geq0.
 \label{eq:1D}
\end{equation}
where $V(t,x)$ is a given smooth velocity field. Let consider a regular mesh, with constant step $\Delta x>0$: the cells are the intervals 
$[x_{i-1/2},x_{i+1/2}], \; i\in\mathbb{N}$ with $x_{-1/2}=0, \; x_{i+1/2}=(i+1)\Delta x$, and $x_{i}$ denotes the midpoint of the cell: 
$x_{i}=(i+1/2)\Delta x$. We denote by $f^{n}_{i}$ the numerical unknown, which is the approximation of $f^{n}_{i}=\displaystyle\frac{1}{\Delta x} \int^{x_{i+1/2}}_{x_{i-1/2}} f(t^{(n)}, x) dx$, where $t^{(0)}=0<t^{(1)} < \cdot\cdot\cdot <t^{(n)} < t^{(n+1)}$ are times discretization with a possible variable time step $\Delta t=t^{(n+1)}-t^{(n)}$. We denote by $V^{n}_{i-1/2}$ the approximations of the velocity at the cell interfaces: namely, we set
\begin{equation*}
V^{n}_{i-1/2} = V(t^{(n)},x_{i-1/2}), \quad n\in\mathbb{N}, \; i\in\mathbb{N}.
\end{equation*}                                            
The finite volume scheme applied to \eqref{eq:1D} gives the following approximation 
\begin{equation} 
f^{n+1}_i= f^{n}_i -\frac{\Delta t^{(n)}}{\Delta x}(V^{n}_{i+1/2}f^{n}_{i+1/2} - V^{n}_{i-1/2}f^{n}_{i-1/2}). 
\label{discret_new}
\end{equation}
Then the main task is how to define the range of the interface fluxes $f^{n}_{i+1/2}$  with respect to the stability, consistency and positivity constraints of the scheme.

 In order to describe the different constraints, let introduce the following useful notations:
 \begin{itemize}
\item $\cfl = \ds\frac{\Delta t^{(n)}}{\Delta x}$,

\item $
 m^{n}_{i+1/2}=\min(f^{n}_{i},f^{n}_{i+1})$, 
and  $
M^{n}_{i+1/2}=\max(f^{n}_{i},f^{n}_{i+1})$,\\

\item if $V^n_{i+1/2}, V^n_{i-1/2}> 0$ :
  \[
  \begin{array}{lll}
    b^{n}_{i+1/2} & = & \displaystyle\frac{1}{\cfl  V^{n}_{i+1/2}}\big(f^{n}_{i} - \max(f^{n}_{i},f^{n}_{i-1})\big) + \max(f^{n}_{i},f^{n}_{i-1}) \\ 
    & = &  \frac{1}{\cfl V^{n}_{i+1/2}}\big(f^{n}_{i} - M^{n}_{i-1/2}\big) + M^{n}_{i-1/2}, \\
    \\
    B^{n}_{i+1/2} & = & \displaystyle \frac{1}{\cfl  V^{n}_{i+1/2}}\big(f^{n}_{i} - \min(f^{n}_{i},f^{n}_{i-1})\big) + \min(f^{n}_{i},f^{n}_{i-1}) \\
    & = & \frac{1}{\cfl V^{n}_{i+1/2}}\big(f^{n}_{i} - m^{n}_{i-1/2}\big) + m^{n}_{i-1/2},
    \\
    
    \mathscr B^{n}_{i+1/2} & = & 
    \left\{
    \begin{array}{ll}
      \min \left(B^n_{i+1/2}, m^{n}_{i-1/2}\ds\frac{V^{n}_{i-1/2}}{V^{n}_{i+1/2}} + \ds\frac{f^{n}_{i}}{\cfl  V^{n}_{i+1/2}}\right), & \text{ if } m^{n}_{i-1/2}\geq 0, \\ [3mm]
      B^n_{i+1/2}, & \text{otherwise},
    \end{array}
    \right.
  \end{array}
  \]

\item if $V^n_{i+1/2}, V^n_{i-1/2}<0$ :
  \[
  \begin{array}{lll}
    b^{n}_{i-1/2} & = & \displaystyle \frac{1}{\cfl |V^{n}_{i-1/2}|}\big(f^{n}_{i} - \max(f^{n}_{i},f^{n}_{i+1})\big) + \max(f^{n}_{i},f^{n}_{i+1}) \\ [3mm]
    & = & \frac{1}{\cfl |V^{n}_{i-1/2}|}\big(f^{n}_{i} - M^{n}_{i+1/2}\big) + M^{n}_{i+1/2}, \\
\\
B^{n}_{i-1/2} & = & \displaystyle \frac{1}{\cfl |V^{n}_{i-1/2}|}\big(f^{n}_{i} - \min(f^{n}_{i},f^{n}_{i+1})\big) + \min(f^{n}_{i},f^{n}_{i+1}) \\[3mm]
& = & \frac{1}{\cfl |V^{n}_{i-1/2}|}\big(f^{n}_{i} - m^{n}_{i+1/2}\big) + m^{n}_{i+1/2},
\\
\mathscr B^{n}_{i-1/2} & = & \left\{\begin{array}{ll}
  \min
  \left(B^n_{i-1/2}, m^{n}_{i+1/2}\ds  \frac{|V^{n}_{i+1/2}|}{|V^{n}_{i-1/2}|} + \frac{f^{n}_{i}}{\cfl  |V^{n}_{i-1/2}|} \right), 
  & \text{ if } m^{n}_{i+1/2}\geq 0, 
  \\[3mm]
   B^n_{i-1/2}, & \text{ otherwise},
\end{array}
\right.
\end{array}
\]

\item if $V^n_{i+1/2}, V^n_{i-1/2}$ do not have the same sign,
  we set $b^n_{i+1/2}=\mathscr B^n_{i+1/2}=f^n_i$ if $V^n_{i+1/2}>0$
  and $b^n_{i+1/2}=\mathscr B^n_{i+1/2}=f^n_{i+1}$ if
  $V^n_{i+1/2}<0$. 
\item $\mu^n_{i+1/2}=\max(m^n_{i+1/2},b^n_{i+1/2})$, and $\mathscr
  M^n_{i+1/2}=\min(M^n_{i+1/2}, \mathscr B^n_{i+1/2})$. 
\end{itemize}

 \subsubsection{Stability constraints}
 The stability constraints is based on the standard {\it{Courant-Friedrichs-Levy (CFL)}} condition:
 \begin{equation}
 0\leq \ds\frac{\Delta t^{(n)}}{\Delta x}\ds\max_{i}(\vert V^{n}_{i+1/2}\vert) \leq 1.\label{cfl}
 \end{equation}
 From this {\it{CFL}} condition, we consider the case where the velocity is positive in the cell $i$: $ V^{n}_{i-1/2}>0$  and  $V^{n}_{i+1/2}>0$,
 so we have 
$$\frac{1}{\cfl V^{n}_{i+1/2}} -1\geq 0$$
 and in first hand we write 
 $$\big(\frac{1}{\cfl
  V^{n}_{i+1/2}}-1\big)\big(f^{n}_{i} - \min(f^{n}_{i},f^{n}_{i-1})\big) \geq 0. 
$$
Then we deduce  
\begin{equation}\label{numero1}
\frac{1}{\cfl
  V^{n}_{i+1/2}}\big(f^{n}_{i} - \min(f^{n}_{i},f^{n}_{i-1})\big) + \min(f^{n}_{i},f^{n}_{i-1})\geq f^{n}_{i}. 
\end{equation}
 
In second hand we have 
$$
\big(\frac{1}{\cfl
  V^{n}_{i+1/2}}-1\big)\big(f^{n}_{i} - \max(f^{n}_{i},f^{n}_{i-1})\big)
\leq 0, 
$$
then we deduce
\begin{equation}\label{numero2}
\frac{1}{\cfl
  V^{n}_{i+1/2}}\big(f^{n}_{i} - \max(f^{n}_{i},f^{n}_{i-1})\big) +
\max(f^{n}_{i},f^{n}_{i-1})\leq f^{n}_{i}. 
\end{equation}
In last hand, assuming that $f^{n}_{i}\geq 0 \; \forall i$ (the positivity constraint is discussed later) we deduce the obvious relation 
\begin{equation}\label{numero3}
\begin{array}{c @{\geq} c}
 \min(f^{n}_{i},f^{n}_{i-1})\ds\frac{V^{n}_{i-1/2}}{V^{n}_{i+1/2}} + \ds\frac{f^{n}_{i}}{\cfl V^{n}_{i+1/2}} & \ds\frac{f^{n}_{i}}{\cfl V^{n}_{i+1/2}} \\
     &  f^{n}_{i}. 
 \end{array}
 \end{equation}
 Owing to \eqref{numero1}-\eqref{numero3} and the previous notations, one deduce this first non empty interval where the adequate flux will be chosen
 \begin{equation}\label{first_bound}
 f^{n}_{i}\in [b^{n}_{i+1/2},  \mathscr B^{n}_{i+1/2}] \neq \emptyset. 
 \end{equation}
 
 For the case where the velocity is locally negative mean $ V^{n}_{i-1/2}<0$ and  $V^{n}_{i+1/2}<0$, we perform the same reasoning 
 and obtain the non empty interval  $f^{n}_{i}\in [b^{n}_{i-1/2},  \mathscr B^{n}_{i-1/2}] \neq \emptyset.$

 \subsubsection{Consistency constraints}
 
 For the consistency of the scheme, we write the very definition which is that the numerical flux between cell $i$ and cell $i+1$ belongs 
 necessary to the interval defined by the numerical solutions $f^{n}_{i}$ and $f^{n}_{i+1}$. So the consistency constraint is written as follows
  \begin{equation}\label{second_bound}
   m^{n}_{i+1/2}\leq f^{n}_{i+1/2} \leq  M^{n}_{i+1/2} , \;\;  \forall \, i  \text{ in the grid}. 
  \end{equation}  
  
  \subsubsection{Positivity constraints}
  From a positive initial solution, we need to impose a condition on the numerical fluxes in order to ensure the positivity of the numerical 
  approximation of the solution given by \eqref{discret_new}. For this, we are looking for conditions so that  
  \begin{eqnarray*}
  f^{n+1}_{i}\geq 0  & \Longrightarrow & \cfl (V^{n}_{i+1/2}f^{n}_{i+1/2} - V^{n}_{i-1/2}f^{n}_{i-1/2}) \leq f^{n}_{i} \\
  & \Longrightarrow & f^{n}_{i+1/2} \leq \ds  f^{n}_{i-1/2}\ds\frac{V^{n}_{i-1/2}}{V^{n}_{i+1/2}} + \ds\frac{f^{n}_{i}}{\cfl  V^{n}_{i+1/2}},
 \end{eqnarray*}
 so, obtaining the later inequality is done by imposing the flux to satisfy 
 \begin{equation}\label{third_bound}
 f^{n}_{i+1/2} \leq \ds m^{n}_{i-1/2}\ds\frac{V^{n}_{i-1/2}}{V^{n}_{i+1/2}} + \ds\frac{f^{n}_{i}}{\cfl  V^{n}_{i+1/2}},
 \end{equation}  
 because by the consistency constraints we know that $m^{n}_{i-1/2} \leq f^{n}_{i-1/2} $.
 
 \begin{proposition}\label{propo1}
Assuming the {\it{CFL}} condition \eqref{cfl} be satisfied. Then for any $i$ the interval $[\mu^{n}_{i+1/2}, \mathscr M^n_{i+1/2} ]$ is 
non empty. By choosing the fluxes $f^{n}_{i+1/2} \in [\mu^{n}_{i+1/2}, \mathscr M^n_{i+1/2} ]$ for any $i$, then the following assertions 
hold:
\begin{itemize}
\item[1)] The scheme \eqref{discret_new} is consistent with \eqref{eq:1D}.
\item[2)] The scheme \eqref{discret_new} remains positive for positive initial solution: mean if $f^{n}_{i}\geq0$ for any $i$ then $f^{n+1}_{i}\geq0$ too.
\item[3)] The scheme satisfies: if $V^{n}_{i} \geq 0$ then $m^{n}_{i-1/2} \leq f^{n+1}_{i} \leq M^{n}_{i-1/2}$, while if $V^{n}_{i} \leq0 $ then 
$m^{n}_{i+1/2} \leq f^{n+1}_{i} \leq M^{n}_{i+1/2}$.   
\end{itemize} 
 \end{proposition}
 
 \noindent
{\bf Proof.}
The proof of the proposition is essentially based on the gathering of the results \eqref{numero1}--\eqref{numero3} and \eqref{second_bound}--\eqref{third_bound}. 
For more details on the proof, one can refer to \cite{bibGLT}. \\

\begin{remark}
In the previous reasoning, we exclude the case where the velocities on the interfaces cell are of different sign. So in the case where $V^{n}_{i+1/2}<0$ and $V^{n}_{i-1/2}>0$ then it is obvious that the positivity of the numerical solution $f^{n+1}_{i}$ is obtained without any time step condition. Nevertheless, the case where $V^{n}_{i+1/2}>0$ and $V^{n}_{i-1/2}<0$ mean the possible empty of the cell from the two sides. So we choose the upwind fluxes and the numerical solution becomes constant in the cell $i$.
\end{remark}   

\begin{remark}
In this above presentation, the numerical flux is designed for general conservation laws, while the one in~\cite{bibGLT} is more suitable for the transport equations.
\end{remark}

For the anti-dissipative strategy, we define the flux $f^{n}_{i+1/2}$ by solving the following minimization problem: 
 \begin{equation*}
\begin{array}{lll} 
\text{To minimize} \; \vert f^{n}_{i+1/2}-f^{n}_{i+1}\vert \\
 \text{under the constraint}\; \;  
 f^n_{i+1/2} \in [\mu^n_{i+1/2},\mathscr M_{i+1/2}^n]
\end{array}
\end{equation*}
which solution is given for instance in the case $V_{i-1/2}^n > 0$, $V_i^n > 0$ and $V_{i+1/2}^n > 0$ by
\begin{equation}\left\{\begin{array}{ll}\label{sol_mini1}  
 f^{n}_{i+1/2} = \mu^n_{i+1/2} \qquad  &\text{if} \;\;
 f^{n}_{i+1}\leq \mu^n_{i+1/2},\\ 
f^{n}_{i+1/2} = f^{n}_{i+1} \qquad&  \text{if} \;\;
\mu^n_{i+1/2} \leq f^{n}_{i+1} \leq \mathscr M_{i+1/2}^n,\\ 
f^{n}_{i+1/2} = \mathscr M_{i+1/2}^n \qquad&  \text{if}\;
\;f^{n}_{i+1}\geq \mathscr M_{i+1/2}^n.\\  
\end{array}\right.
\end{equation}

 This kind of anti-dissipative method is very suitable for discontinuous initial solution, which has been shown in~\cite{bibdespres}. However,  it is not suitable for smooth solution. Indeed, it turns the very smooth solution into a series of step functions with respect to time evolution~\cite{bibdespres, bibGLT}. The objective of remain part is to find an alternative method such that it keeps the shock near discontinuities while has high accuracy in smooth regions. To the end, we propose the following hybrid method.

 \subsection{A hybrid method}
 We denote $f_{i+1/2}^{\text{A}}$ the flux computed by the  anti-dissipative method and $f_{i+1/2}^{\text{W}}$ the flux computed by a high accurate method, for instance the fifth order method~\cite{bibJS}. So our desirable flux $f_{i+1/2}^{\text{H}}$ by the hybrid method will just be a convex combination of  $f_{i+1/2}^{\text{A}}$ and  $f_{i+1/2}^{\text{W}}$, {\it i.e.}
\begin{equation}
  f_{i+1/2}^{\text{H}} = \omega_{i+1/2}^{\text{A}} f_{i+1/2}^{\text{A}} + \omega_{i+1/2}^{\text{W}} f_{i+1/2}^{\text{W}}
\end{equation}
where $\omega_{i+1/2}^{\text{A}} + \omega_{i+1/2}^{\text{W}} = 1$, $\omega_{i+1/2}^{\text{A}}, \omega_{i+1/2}^{\text{W}} \geq 0$. Moreover, it is desirable:
\begin{itemize}
  \item $\omega_{i+1/2}^{\text{A}} = \mathcal{O}(1)$ and $\omega_{i+1/2}^{\text{W}} = {o}(1)$, near discontinuities,
  \item $\omega_{i+1/2}^{\text{A}} = \mathcal{O}(\Delta x^M)$ and $\omega_{i+1/2}^{\text{W}} = \mathcal{O}(1)$, in smooth regions,
\end{itemize}
where $M$ is a large enough number that causes $f_{i+1/2}^{\text{W}}$ to be the dominant term in smooth regions. 

\begin{figure}
  \begin{center}
    \includegraphics[width=7cm]{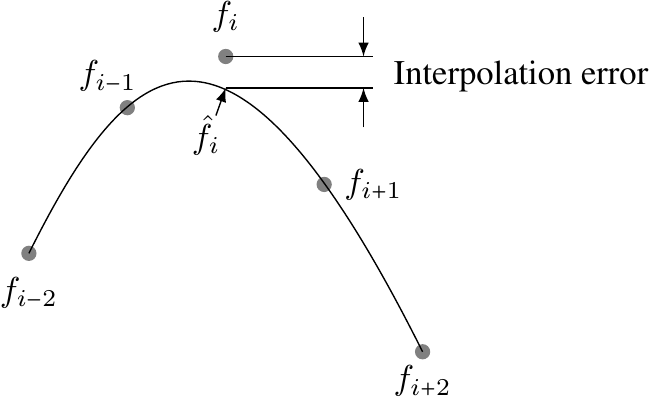}
\caption{\label{fig:interpolation}Interpolation error between the node value $f_i$ and interpolated value $\hat{f}_i$.}
 \end{center}
\end{figure}

According to these considerations, we then need  to identify the smoothness of solution. We consider a similar smoothness measurement, which was first proposed in~\cite{bibDarian}
\begin{equation}\label{eq:smooth_indicator1}
  e_i = \frac{|\hat{f}_i-f_i|}{D},
\end{equation}
where $D$ is a scaling value given as
\begin{equation}\label{eq:smooth_indicator2}
  D = \max_i f_i -  \min_i f_i.
\end{equation}
Figure~\ref{fig:interpolation} shows the interpolated value, {\it i.e.} $\hat{f}_i$ is obtained by the fourth-order interpolation. Note that the point $f_i$ itself, as shown in Figure~\ref{fig:interpolation}, is not included in the interpolation. 

\begin{equation}\label{eq:smooth_indicator3}
  \hat{f}_i = \frac{1}{6}(-f_{i-2} + 4f_{i-1} + 4f_{i+1} - f_{i+2}).
\end{equation}

Then we choose the  smoothness measurement at the cell interface by an upwind way
\begin{equation}
  e_{i+1/2} = \left\{
  \begin{array}{ll}
    e_i,& \text{ if } V_{i+1/2} \geq 0,\\[3mm]
    e_{i+1},& \text{ else}.
  \end{array}
  \right.
\end{equation}
Finally, the weights have forms
\begin{equation}
  \omega_{i+1/2}^{\text{A}} = 1 - \exp(-e_{i+1/2}^2/c_x),\quad \omega_{i+1/2}^{\text{W}} = \exp(-e_{i+1/2}^2/c_x),
\end{equation}
where 
\begin{equation} 
  c_x = \left(\frac{\Delta x}{L_x}\right)^\alpha.
\label{eq:parameter_cx}
\end{equation}
The parameter $\alpha$ is a constant, and a discussion of its value will be given in the section~\ref{sec:parameter_cx}. We now present the order analysis of the weights. In smooth regions, using a Taylor expansion analysis, we have:
\begin{equation*}
  e_{i+1/2} = \mathcal{O}(\Delta x^4).
\end{equation*}
Now using the Taylor expansion of the exponential function and assuming $\alpha \leq 1$, one obtains:
\begin{equation*}
  \omega_{i+1/2}^{\text{A}} = 1 - \exp(-e_{i+1/2}^2/c_x)  \approx  e_{i+1/2}^2/c_x  = \mathcal{O}(\Delta x^{8-\alpha}),
\end{equation*}
and 
\begin{equation*}
  \omega_{i+1/2}^{\text{W}} = \exp(-e_{i+1/2}^2/c_x)  \approx  1 - e_{i+1/2}^2/c_x  = \mathcal{O}(1).
\end{equation*}
Therefore, we use the fifth order WENO flux in smooth regions, {\it i.e.}
\begin{equation*}
  f_{i+1/2}^{\text{H}} \approx f_{i+1/2}^{\text{W}}.
\end{equation*}

For a discontinuity, since the size of the discontinuity does not change as $\Delta x\to0$, one can conclude
\begin{equation*}
  e_{i+1/2} = \mathcal{O}(1).
\end{equation*}
Therefore, applying the same analysis used above leads to
\begin{equation*}
  \omega_{i+1/2}^{\text{A}} =  \mathcal{O}(1).
\end{equation*}
Thus, at discontinuities the anti-dissipative flux is activated.

\begin{remark}
In the original method~\cite{bibDarian}, the parameter $D$ depended on the index $i$ and it was a convex combination of the global and local scales, {\it i.e.}
\begin{equation*}
  D_i=c_s S_g + (1-c_s) S_l,
\end{equation*}
where $S_g=\max_i f_i -  \min_i f_i,\, i\in\{1,\dots,i_{\max}\}$, $S_l = \max_i f_i -  \min_i f_i,\, i\in\{i-2,\dots,i+2\}$. The parameter $c_s$ was chosen as $0.1$ or $0.01$ in~\cite{bibDarian}. This choice can highlight the small jumps in the solution. However, in the paper we focus on the major jumps in the solution. Therefore, the global scale $S_g$ seems more appropriate.
\end{remark}  

\begin{remark}
Again, in the original method~\cite{bibDarian}, the parameter $c_x$ is a constant, which varies from problem to problem. At contrast, we propose the parameter $c_x$ depending on the scaled mesh step, {\it i.e.} $\Delta x/L_x$ and the exponent $\alpha$. Moreover, according to the numerical experience in section~\ref{sec:parameter_cx}, the parameter $\alpha$ can be a fixed number.
\end{remark}  

\section{Numerical results}
\label{sec:detector}
\setcounter{equation}{0}

In this section, we will present several numerical results to depict the behaviors of our hybrid method and its applications in population dynamics.
\subsection{The classical numerical tests}
Here, we first discuss the  parameter $\alpha$ in~\eqref{eq:parameter_cx}, which is the key point in our hybrid method. Then we perform the classical tests with the 1D free transport equation and the 2D rotation equation to verify the convergence of our hybrid method for the regular and irregular initial data. In the sequel, the third order Runge-Kutta method is used for time discretization.
\subsubsection{A discussion of the parameter $\alpha$ in~\eqref{eq:parameter_cx}}
\label{sec:parameter_cx}
To determine the parameter $\alpha$, we use  the free  transport equation
\begin{equation}
 \frac{\partial f}{\partial t}\,\,+\,\,\,  \frac{\partial f}{\partial x}\,\,=\,\,0,\quad x\in[-1,1],\quad t\geq0,
 \label{eq:test1D}
\end{equation}
with a very oscillatory  initial condition, given by~\cite{bibJS}, 
\begin{equation}
 f(0,x)=
 \left\{
 \begin{array}{ll}
  f_1(x) = \frac{1}{6}[G(x,z-\delta)\,+\,G(x,z-\delta)\,+\,4\,G(x,z)],&\textrm{ if } -0.8\leq x\leq-0.6,\\[3mm]
  f_2(x) = 1,&\textrm{ if } -0.4\leq x\leq-0.2,\\[3mm]
  f_3(x) = 1-|10(x-0.1)|, &\textrm{ if } 0\leq x\leq0.2,\\[3mm]
  f_4(x) =  \frac{1}{6}[F(x,z-\delta)\,+\,F(x,z-\delta)\,+\,4\,F(x,z)],&\textrm{ if } 0.4\leq x\leq0.6,\\[3mm]
  0,&\text{ otherwise}.
 \end{array}
 \right.
  \label{eq:discontinuous_initial2}
\end{equation}
where $G(x,z)=\exp(-\beta(x-z)^2)$, $F(x,a)=\{\max((1-\alpha^2(x-a)^2)^{1/2},0)\}$ with  $\alpha=0.5$, $z=-0.7$, $\delta=0.005$, $\alpha=10$ and $\beta=(\log2)/36\delta^2$. Moreover, the periodic boundary conditions is considered.

\begin{figure}
\begin{center}
 \begin{tabular}{cc}
  \includegraphics[width=5cm]{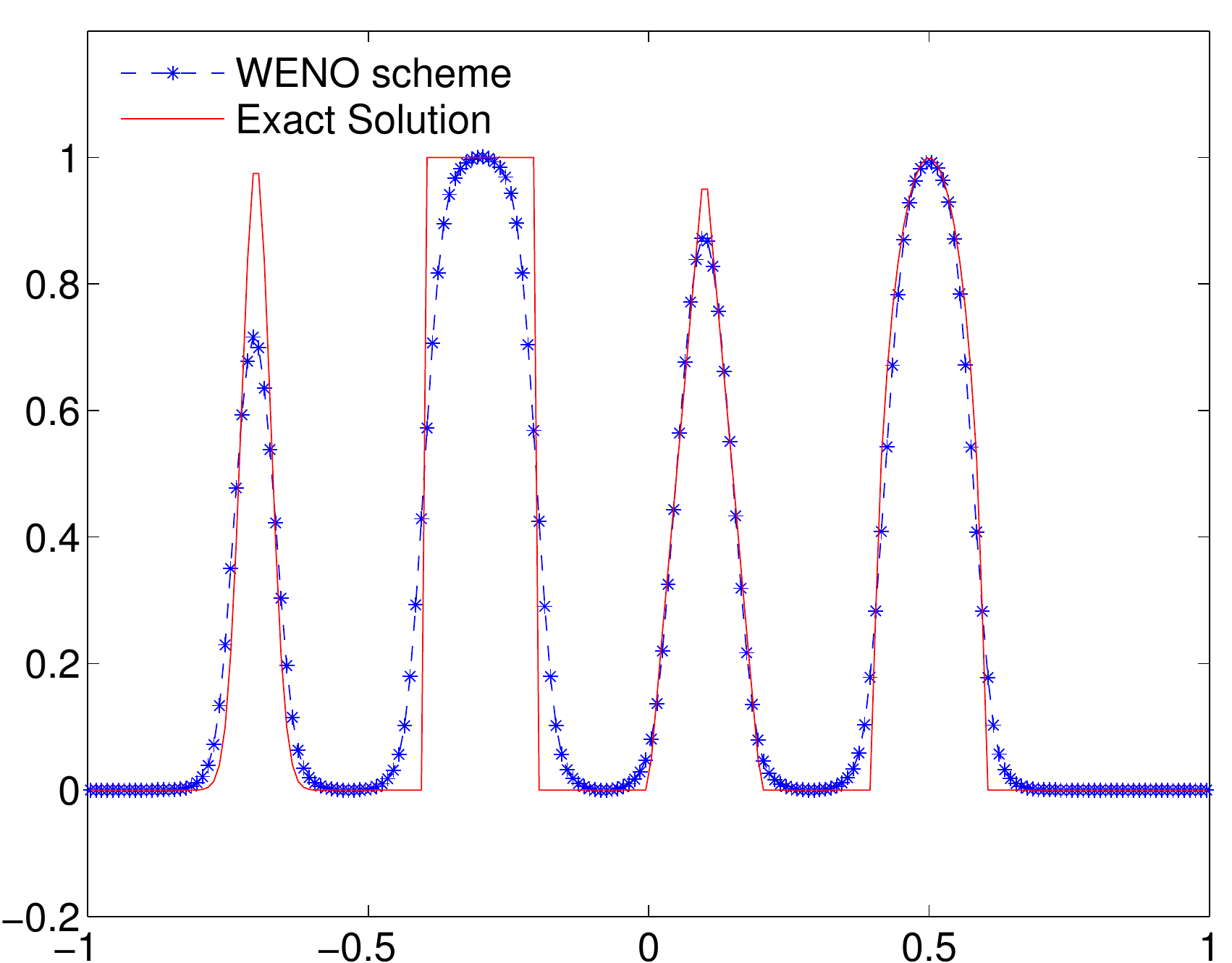} &    \includegraphics[width=5cm]{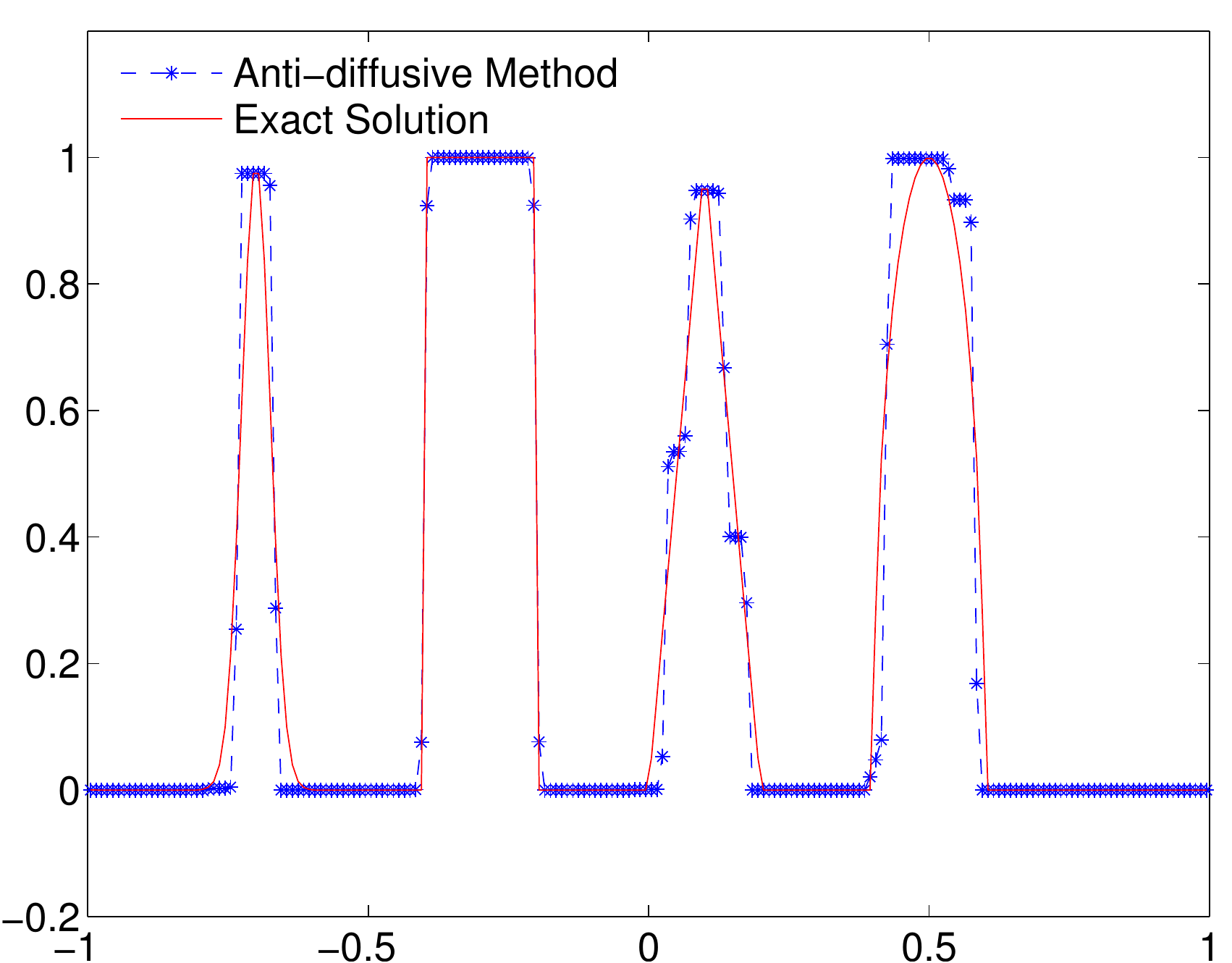}\\
  (a)                  &    (b) \\
   \includegraphics[width=5cm]{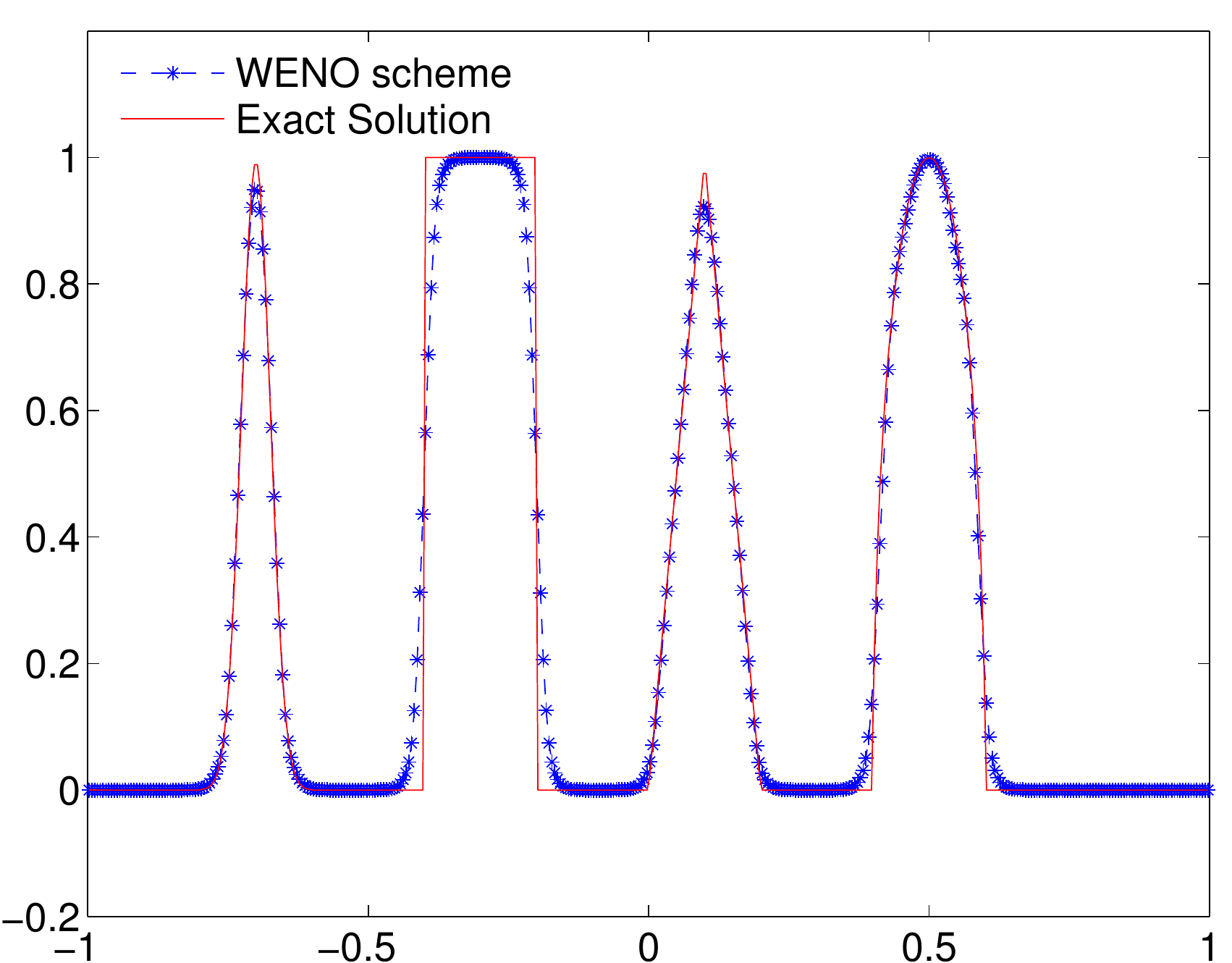} &    \includegraphics[width=5cm]{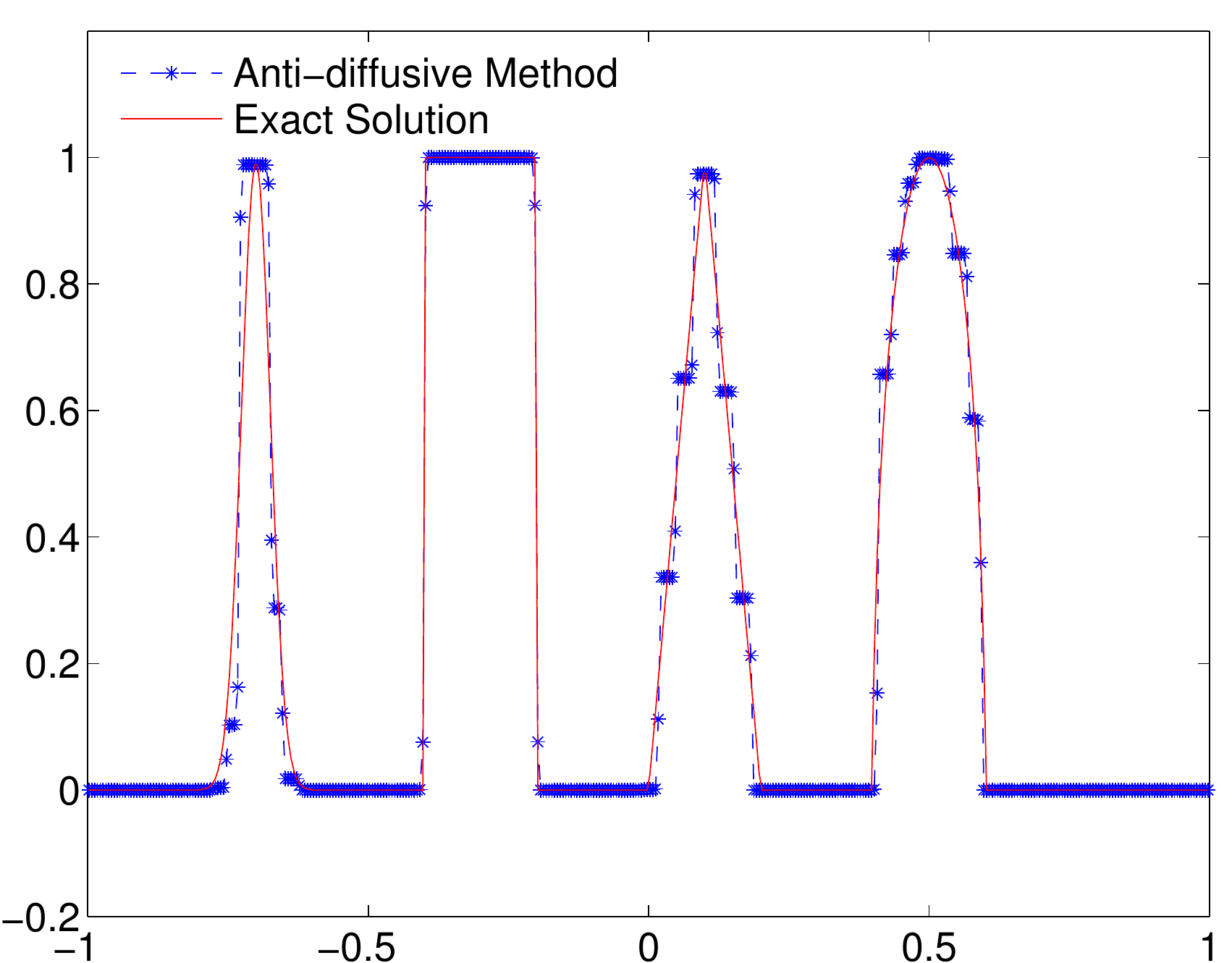} \\
  (c)                  &    (d)  
 \end{tabular}
\caption{\label{fig:1Dtransport}Numerical solutions of the WENO method and the anti-diffusive method for 1D transport equation : {\it Plot solutions of the linear equation~\eqref{eq:test1D} with initial data~\eqref{eq:discontinuous_initial2}. $n_x=200$ in the first row and $n_x=400$ in the second row, while $\Delta t$ is chosen such that CFL$=0.2$. The  final time is $T_{\text{end}}=8$.}}
 \end{center}
\end{figure}

The numerical results is presented in Figure~\ref{fig:1Dtransport}. We first see that the WENO method is well adapted for the smooth regions of solution, while it becomes significantly diffuse near the steep regions. For the refined mesh, we can always observed this diffusion near the step function clearly. The anti-diffusive method has perfect performance for the step function, however it destroys the very smooth function. Indeed, it turns the smooth function into a sawtooth profile, which can not be regularized by refining mesh.

\begin{figure}
\begin{center}
 \begin{tabular}{ccc}
  \includegraphics[width=4.5cm]{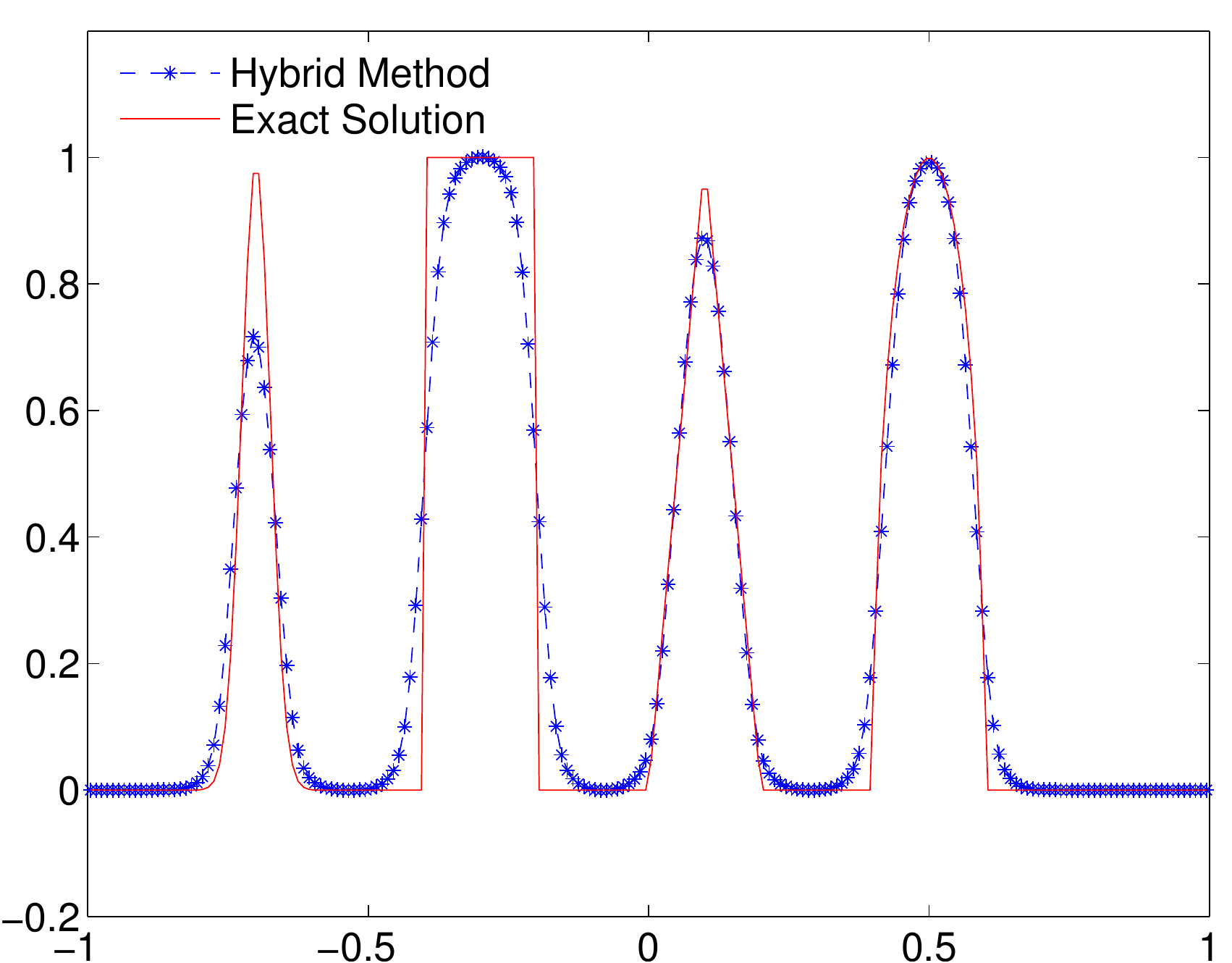} &    \includegraphics[width=4.5cm]{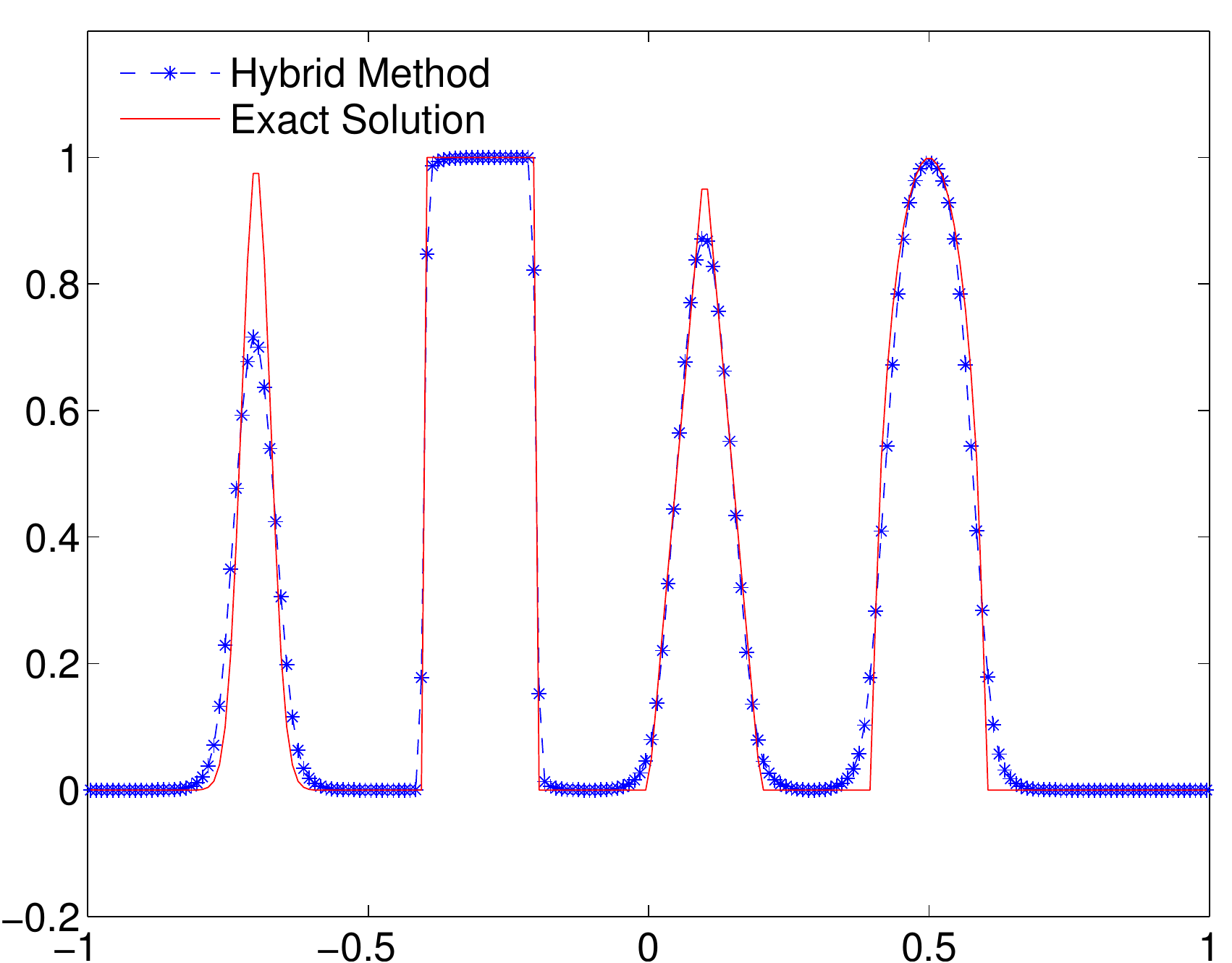} &   \includegraphics[width=4.5cm]{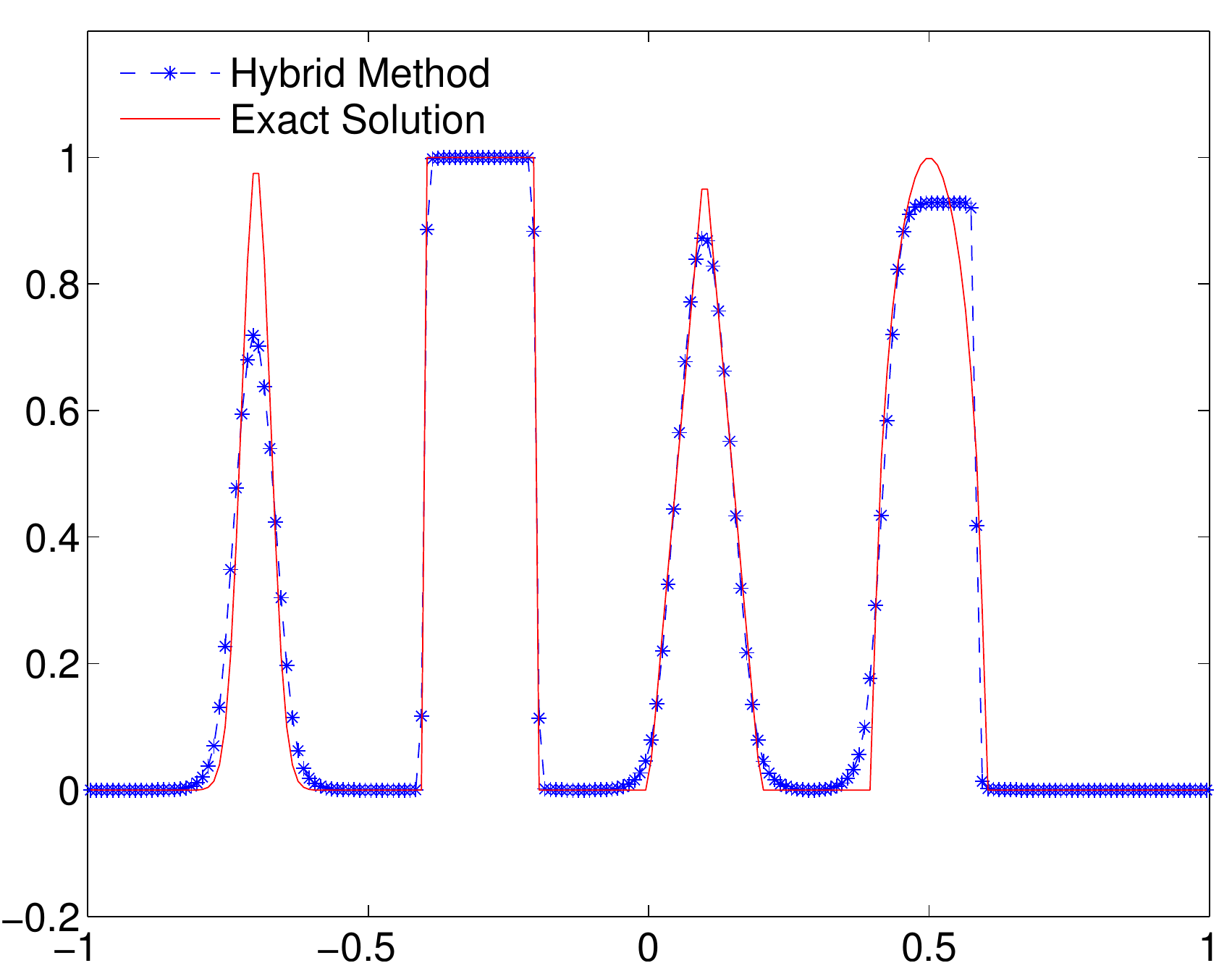}\\
  (a) $\alpha=0.65$                 &    (b) $\alpha=0.75$  &  (c) $\alpha=0.85$ \\
   \includegraphics[width=4.5cm]{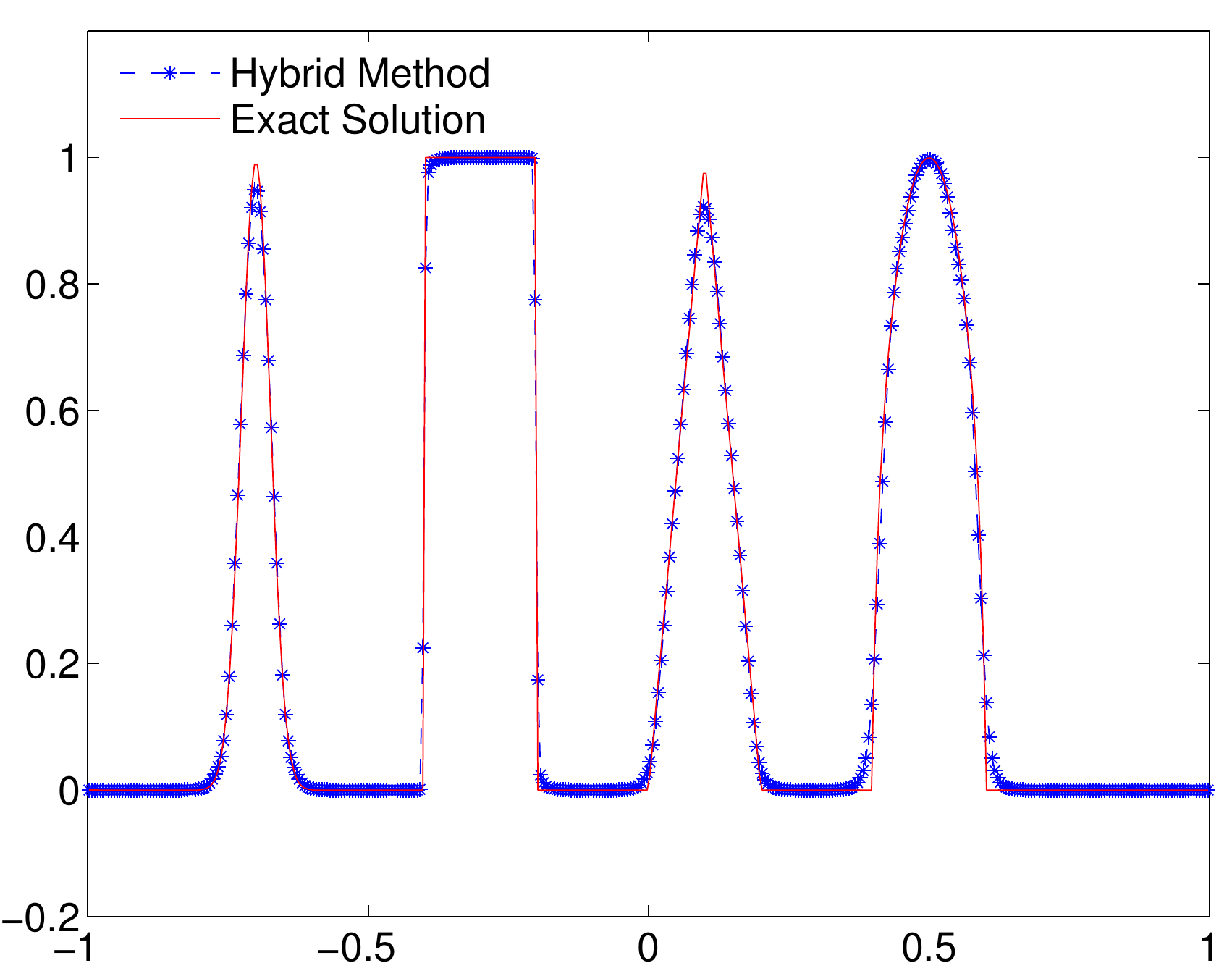} &    \includegraphics[width=4.5cm]{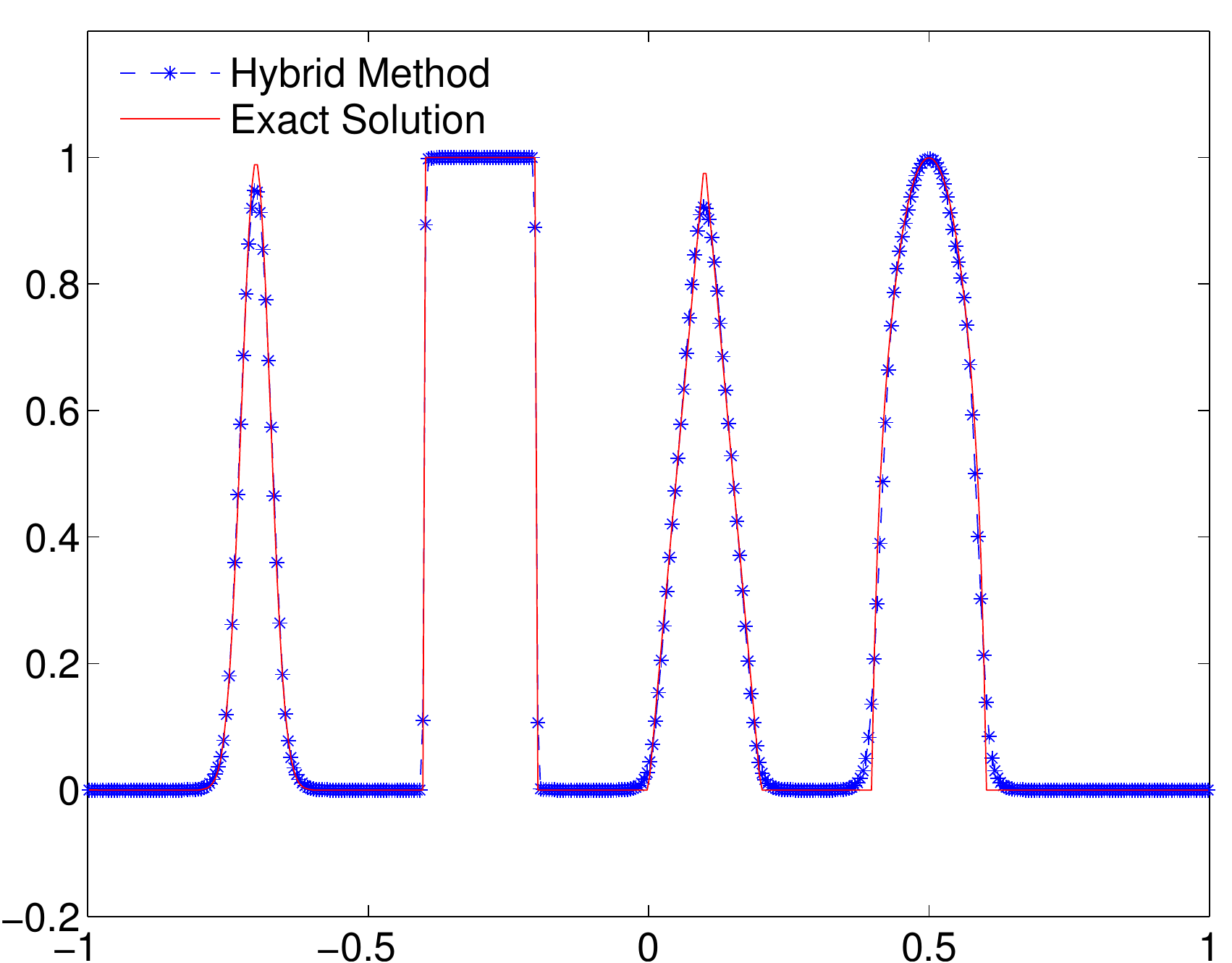} &   \includegraphics[width=4.5cm]{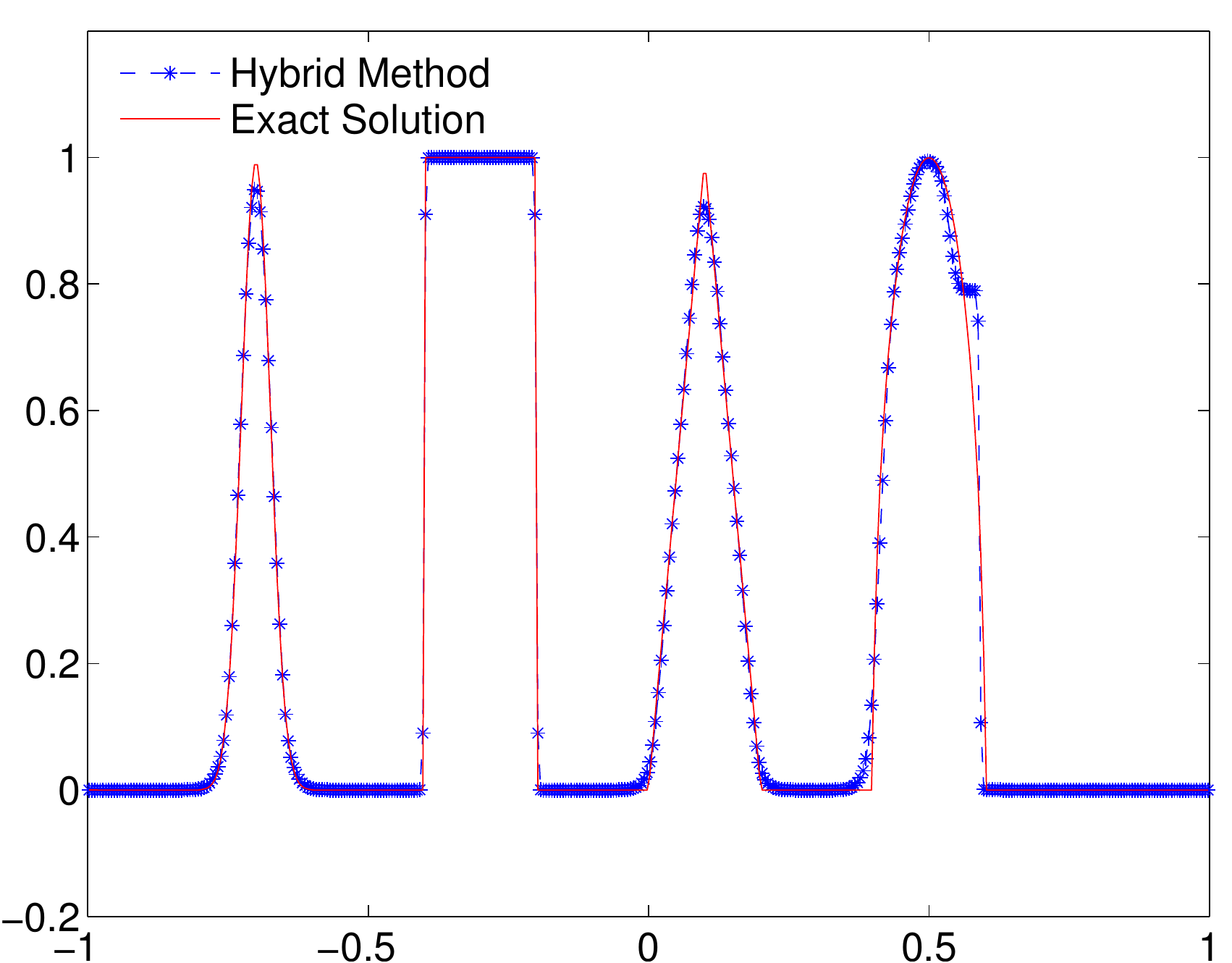}\\
  (d)  $\alpha=0.65$                 &    (e) $\alpha=0.75$  &  (f) $\alpha=0.85$ 
 \end{tabular}
\caption{\label{fig:1Dtransport_alpha}Comparison of different values of the parameter $\alpha$ in~\eqref{eq:parameter_cx} of the hybrid method for  1D transport equation  : {\it Plot solutions of the linear equation~\eqref{eq:test1D} with initial data~\eqref{eq:discontinuous_initial2}. $n_x=200$ in the first row and $n_x=400$ in the second row, while $\Delta t$ is chosen such that CFL$=0.2$. The  final time is $T_{\text{end}}=8$.}}
 \end{center}
\end{figure}

Now we expect that the hybrid method can combine the advantages of both the WENO method and the anti-diffusive method. However, it remains to choose an appropriate parameter $\alpha$ in~\eqref{eq:parameter_cx}. In Figure~\ref{fig:1Dtransport_alpha}, we present the numerical solutions of the hybrid methods with the different values of the parameter $\alpha$, {\it i.e.} $\alpha=0.65,\,0.75,\,0.85$. We first observe that the hybrid method approximates well the functions $f_1$ and $f_3$ of~\eqref{eq:discontinuous_initial2} for all the values of $\alpha$. However, there are significantly differences for the function  $f_2$ and $f_4$. Indeed, with $\alpha=0.65$, the hybrid method can not preserve the  steep profile of the step function $f_2$; while with  $\alpha=0.85$, the sawtooth profile appears again in the function $f_4$. By refining the mesh, we see that the approximations of the functions $f_2$ and $f_4$ are improved, but they are still not perfect. Finally, with $\alpha=0.75$, we are capable to capture well all configurations of the solution~\eqref{eq:discontinuous_initial2}. Therefore, we will use the value of the parameter $\alpha=0.75$ in the sequel of this paper.

\subsubsection{Convergence rate for the 1D transport equation}
\label{sec:1DTransport}

 Let us first consider a smooth solution, where the  initial condition is chosen as
\begin{equation*}
 f(0,x)=\sin\left(\pi x\right),\quad x\in[-1,1].
\end{equation*} 
The numerical error for  different methods is  presented in Table~\ref{tab:1D1}. We first notice that the anti-diffusive method has only first order of convergence rate. At contrast,  the WENO method and the hybrid method have both the fifth order convergence rate. Thus they are much more precise than the  anti-diffusive method. Moreover, the hybrid method does not disturb at all the precision of numerical results  in smooth solution case.

\begin{table}[h]
\begin{center}
 \begin{tabular}{|c|c|c|c|c|c|c|}
  \hline
  $n_x$&\multicolumn{2}{c|}{200}  & \multicolumn{2}{c|}{400}   & \multicolumn{2}{c|}{800} \\
  \cline{2-7}
  & $\|\cdot\|_1$  & $r$  & $\|\cdot\|_1$  & $r$  & $\|\cdot\|_1$  & $r$ \\
  \hline    
  Anti-diffusive method&    4.31e-2 &  1.07   & 2.06e-2   & 1.07   & 1.06e-2   & 0.96     \\
  \hline
  WENO scheme &     1.14e-7 &  5.00  &   3.57e-9     &  5.00 &   1.12e-10  &  4.99   \\
  \hline
  Hybrid method&   1.14e-7 &  5.00  &   3.57e-9     &  5.00 &   1.13e-10  &  4.99   \\
  \hline
 \end{tabular}
\end{center}
\caption{\label{tab:1D1}1D transport equation : {\it Error  in $L_1$-norm and order of convergence $r$. The  final time is $T_{\text{end}}=8$.}}
\end{table}

We next  consider a step function as follows
\begin{equation}
 f(0,x)=
 \left\{
 \begin{array}{ll}
  1,&\textrm{ for } -1\leq x\leq0,\\[3mm]
  0,&\text{otherwise}.
 \end{array}
 \right.
  \label{eq:discontinuous_initial1}
\end{equation}

Comparisons of the methods are now summarized in Table~\ref{tab:1D2}. We observe that the anti-diffusive method has the first order convergence rate for irregular solution, while the WENO method can not achieve  the first order convergence rate. The hybrid method is less precise than the  anti-diffusive method, however it is  more precise of one degree than the WENO method. It is interesting to notice that the hybrid method has a convergence rate even more than one for irregular solution. At last, we see that the total variation for the hybrid method is comparable with the WENO method, which means that the hybrid method is also an essentially non-oscillatory method.

\begin{table}[h]
\begin{center}
 \begin{tabular}{|c|c|c|c|c|c|c|}
  \hline
  $n_x$&\multicolumn{2}{c|}{200}  & \multicolumn{2}{c|}{400}   & \multicolumn{2}{c|}{800} \\
  \cline{2-7}
  & $\|\cdot\|_1$  & $r$  & $\|\cdot\|_1$  & $r$  & $\|\cdot\|_1$  & $r$ \\
  \hline
  Anti-diffusive method&    1.48e-3 &  1.00   & 7.39e-4   & 1.00   & 3.69e-4   & 1.00     \\
  \hline  
  WENO scheme &    4.50e-2 &  0.83  &   2.53e-2     &  0.83 &   1.43e-2  &  0.82   \\
  \hline
  Hybrid method&    7.93e-3 &  3.32  &   2.36e-3     &  1.75 &   9.24e-4  &  1.35    \\
  \hline
 \end{tabular}
 \vskip 3mm
 (a) Error between exact solution and approximated solution
  \vskip 3mm
   \begin{tabular}{|c|c|c|c|}
  \hline
  $n_x$&200&400&800\\
  \hline
  Anti-diffusive method&  -1.33e-15      &  -1.55e-15  &  -1.78e-15    \\
  \hline
  WENO scheme &  9.71e-7      &  9.52e-7      & 7.75e-7    \\
  \hline
  Hybrid method&  8.51e-7      &  5.28e-6  &  1.16e-6  \\
  \hline
 \end{tabular}
  \vskip 3mm
 (b) Error of total variation
  \vskip 3mm
 \caption{ \label{tab:1D2}1D transport equation : {\it  Comparison of different methods for the linear equation~\eqref{eq:test1D} with initial data~\eqref{eq:discontinuous_initial1}. (a) Error in $L_1$ norm and $r$ is the order of accuracy (b) Error on the total variation. The  final time is $T_{\text{end}}=8$.}}
 \end{center}
\end{table}

Finally, the hybrid method does not significantly increase the computational cost. Indeed, in average the WENO scheme and the ADM method take 74 operations and 27 operations respectively. While the hybrid method consists of these two methods (103 operations) and in addition the computation of smoothness indicator (27 operations). So in total the hybrid method takes approximately 130 operations, about 1.76 times of the WENO scheme or 4.81 times of the ADM method. However, considering its important features, we think the additional computational cost of our hybrid method is acceptable.

\subsubsection{2D rotation equation}
\label{sec:2DRotation}
We can directly extend the hybrid method to two dimensional cases. In fact, the numerical flux can be computed dimension by dimension. For example, in the $x$ direction, the numerical flux can be computed by
\begin{equation}
  f_{i+1/2,j}^{\text{H}} = \omega_{i+1/2,j}^{\text{A}} f_{i+1/2,j}^{\text{A}} + \omega_{i+1/2,j}^{\text{W}} f_{i+1/2,j}^{\text{W}}
\end{equation}
where $\omega_{i+1/2,j}^{\text{A}} + \omega_{i+1/2,j}^{\text{W}} = 1$, $\omega_{i+1/2,j}^{\text{A}}, \omega_{i+1/2,j}^{\text{W}} \geq 0$. Then we choose the  smoothness measurement at the cell interface by an upwind way
\begin{equation}
  e_{i+1/2,j} = \left\{
  \begin{array}{ll}
    e^x_{i,j},& \text{ if } V_{i+1/2,j} \geq 0,\\[3mm]
    e^x_{i+1,j},& \text{ else},
  \end{array}
  \right.
\end{equation}
where the smoothness measurement $e^x_{i,j}$ is the same as in~\eqref{eq:smooth_indicator1}-\eqref{eq:smooth_indicator3} for $j$ fixed.
The weights have forms
\begin{equation}
  \omega_{i+1/2,j}^{\text{A}} = 1 - \exp(-e_{i+1/2,j}^2/c_x),\quad \omega_{i+1/2,j}^{\text{W}} = \exp(-e_{i+1/2,j}^2/c_x),
\end{equation}
where 
\begin{equation} 
  c_x = \left(\frac{\Delta x}{L_x}\right)^\alpha.
\end{equation}
Similarly, we can compute the numerical flux $f_{i,j+1/2}^{\text{H}} $ in $y$ direction.

To test the efficiency of the hybrid method, we use the famous Zalesak' disk test. The governed equation is the 2D rotation equation
\begin{equation}
 \frac{\partial f}{\partial t}\,\,+\,\,\,  y\frac{\partial f}{\partial x}\,\, - \,\,\,  x\frac{\partial f}{\partial y}\,\,=\,\,0,\quad x\times y\in[-1,1]\times[-1,1],\quad t\geq0.
 \label{eq:test2D}
\end{equation}

The initial solution is~\cite{bibLeveque}
\begin{equation}
 f(0,x)=
 \left\{
 \begin{array}{ll}
  1,&\textrm{ if } \sqrt{x^2 + (y-0.5)^2}\leq r_0 \,\&\&\, (|x|\geq0.025 \,||\, y\geq0.75),\\[3mm]
  1-\frac{1}{r_0}\sqrt{x^2 + (y+0.5)^2},&\textrm{ if }  \sqrt{x^2 + (y+0.5)^2}\leq r_0,\\[3mm]
  \frac{1 + \cos{\pi} \frac{1}{r_0}\sqrt{(x+0.5)^2 + y^2}}{4}, &\textrm{ if } \sqrt{(x+0.5)^2 + y^2}\leq r_0,\\[3mm]
  0,&\text{ otherwise},
 \end{array}
 \right.
  \label{eq:discontinuous_initial2D}
\end{equation}
where the radius $r_0=0.3$.

\begin{figure}
\begin{center}
 \begin{tabular}{cc}
  \includegraphics[width=7cm]{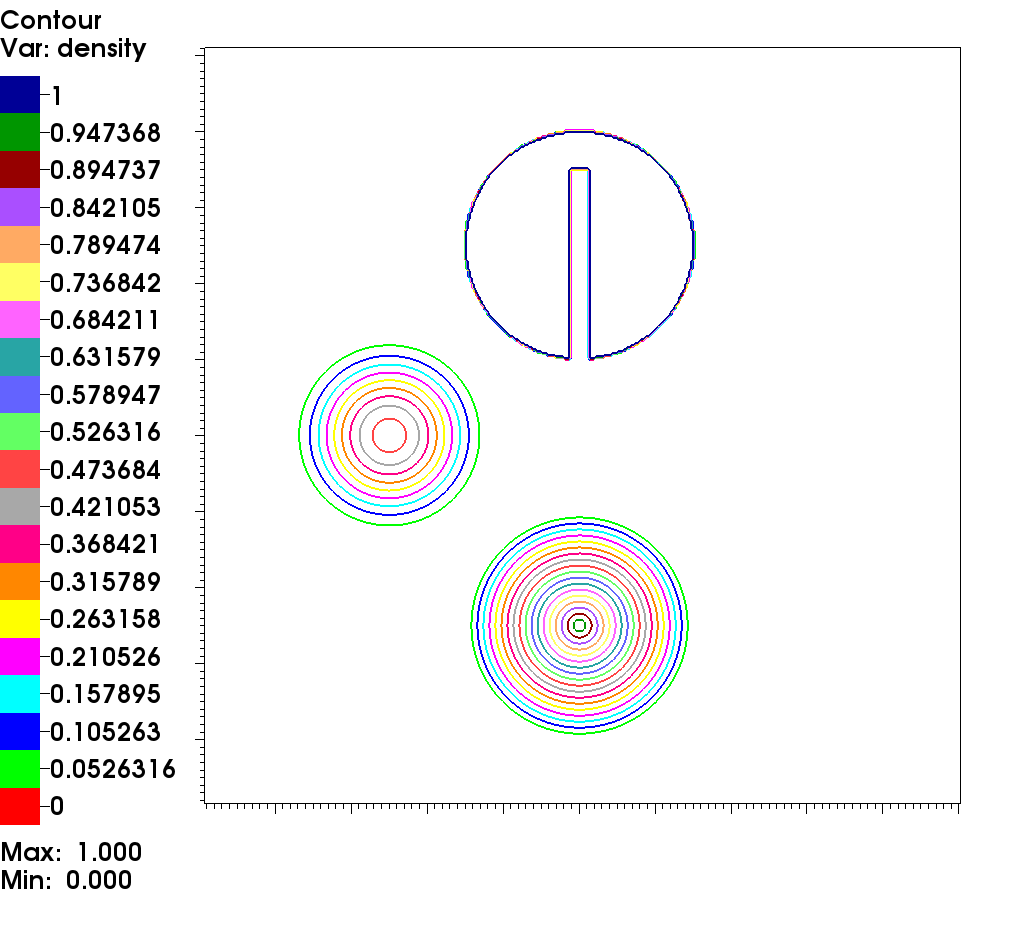}  &   \includegraphics[width=7cm]{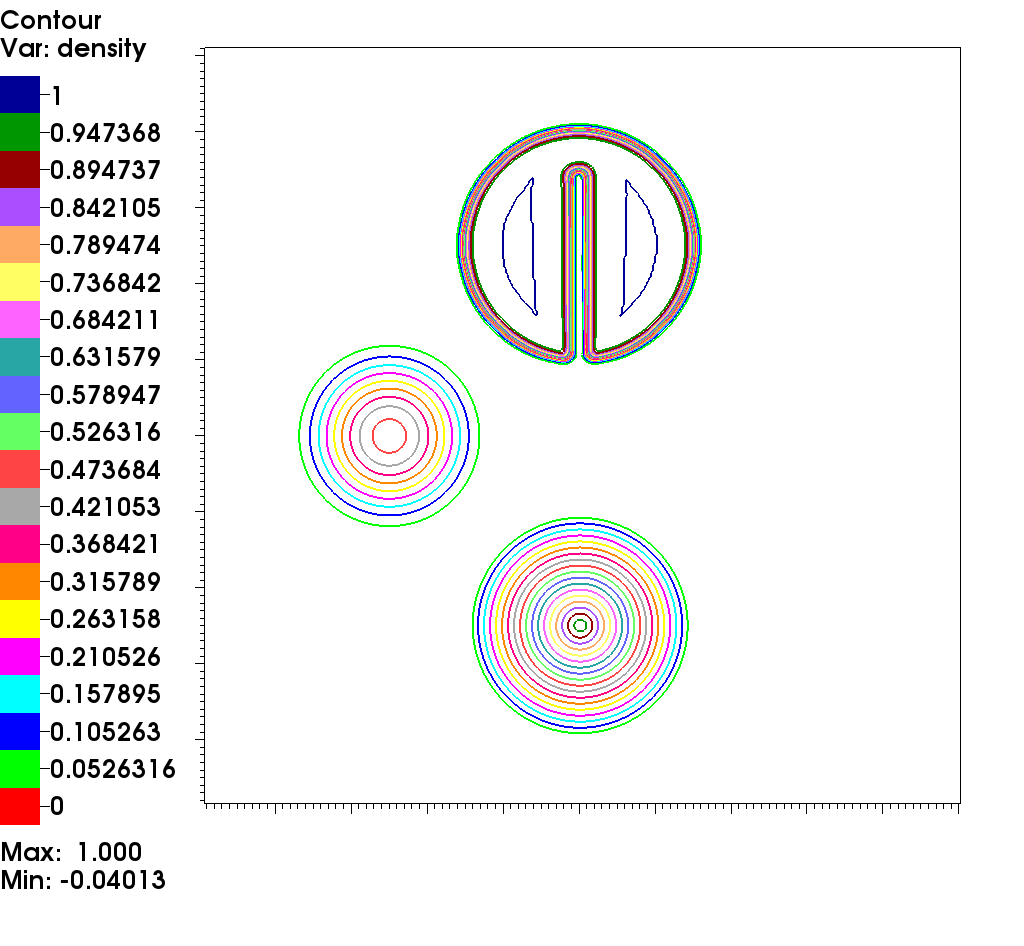}\\
  (a)  Reference solution                &    (b)WENO method \\
   \includegraphics[width=7cm]{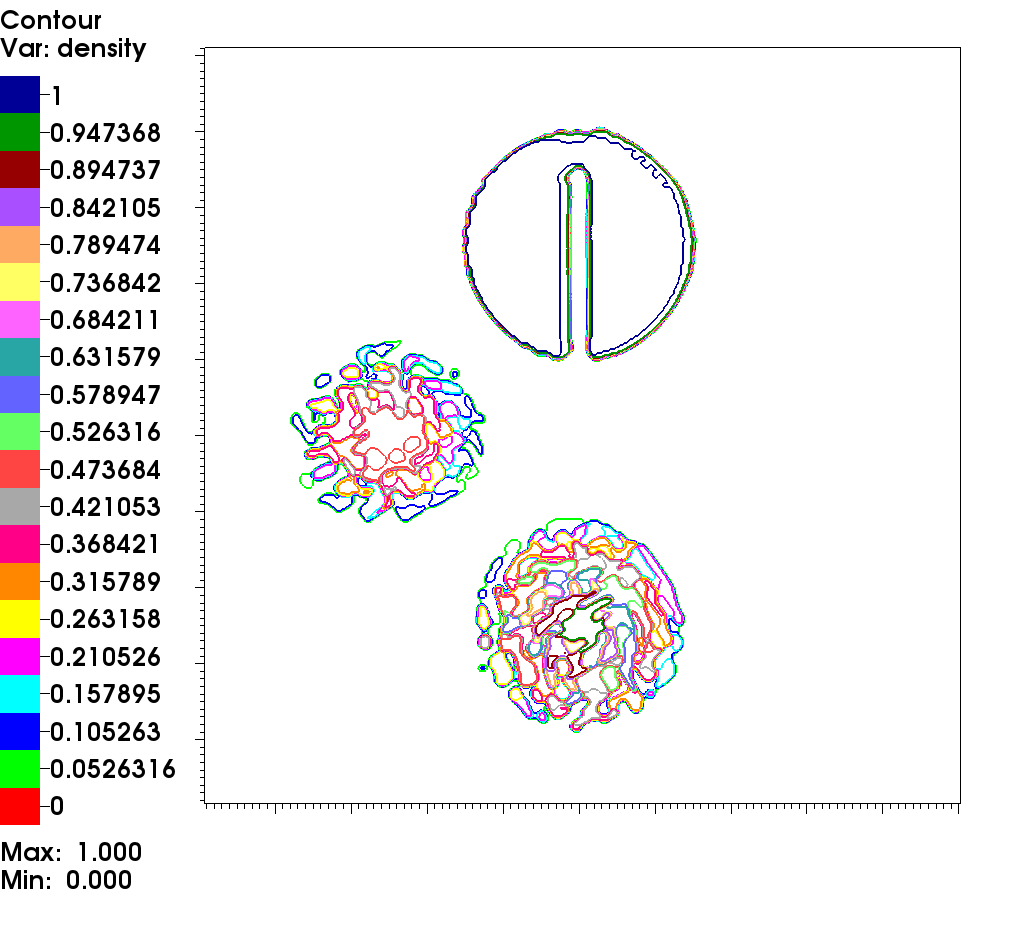} &   \includegraphics[width=7cm]{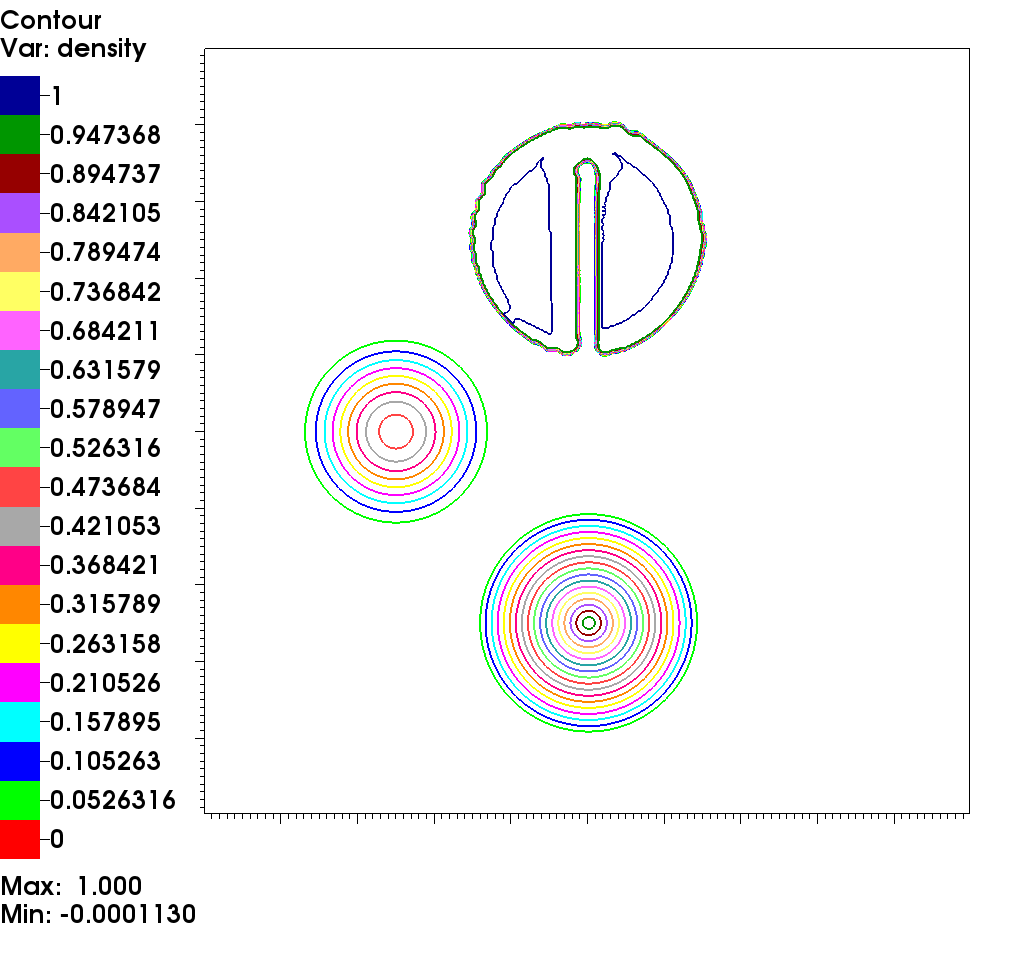}\\
  (c)  Anti-diffusive method                & (d) Hybrid method 
 \end{tabular}
\caption{\label{fig:2Drotation}2D rotation equation : {\it Plot solutions of the linear equation~\eqref{eq:test2D} with initial data~\eqref{eq:discontinuous_initial2D}. Mesh size is $n_x\times n_y=400\times400$, while $\Delta t$ is chosen such that CFL$=0.2$. The  final time is $T_{\text{end}}=2\pi$.}}
 \end{center}
\end{figure}

We illustrate the numerical results in Figure~\ref{fig:2Drotation}. We first notice that the initial condition consists of the Zalesak' disk, the conical body and the peak of the hump. Then after one period, the WENO method preserves well the the conical body and the peak of the hump, however a clear diffusion appears in the  Zalesak' disk. The anti-diffusive method keeps well the shape of  the Zalesak' disk, but it destroys completely the others two objects. Finally, we observe that the hybrid method performs well for all these three objects. Therefore, the hybrid method is suitable for both the smooth solution and the irregular solution for the advection equations.

\subsection{Biological numerical tests}
In this subsection, we will focus on  numerical simulations for the polymerization/depolymerization type model. The objective here is to highlight the good performance of our hybrid method on the long term asymptotic behavior of the solution.

\subsubsection{Space-homogeneous polymerization/depolymerization type model}
Here we are interested on the numerical behavior of the standard Lifshitz-Slyozov equations \cite{bibLS} that we refer as homogeneous in space in comparison to the system \eqref{eq1}--\eqref{eq2}.

 This standard Lifshitz-Slyozov equations  can be interpreted in the point of view of population dynamics as a model describing a population of cells evolving only by nutrients uptake (without birth, death and division) where the nutrients are characterized by their concentration $c(t)$ which fulfills a mass preservation equation as follows 
\begin{equation}
\left\{ \begin{array}{ll}\label{lif-sly}
\displaystyle\frac{\partial}{\partial t} f(t,\xi) +\frac{\partial}{\partial \xi} \bigl((\xi^{1/3}c(t)-1)f(t,\xi)\bigr) =0 ,\\
\\
c(t) + \displaystyle\int_{0}^{\infty} \xi f(t,\xi) d\xi  = \rho,\\
\end{array} \right.
\end{equation}      
where $\rho$ is a constant and measures the total initial mass, $\xi$ the cell-size and $f$ the size density repartition.  

\begin{figure}
\begin{center}
 \begin{tabular}{cc}
  \includegraphics[width=6cm]{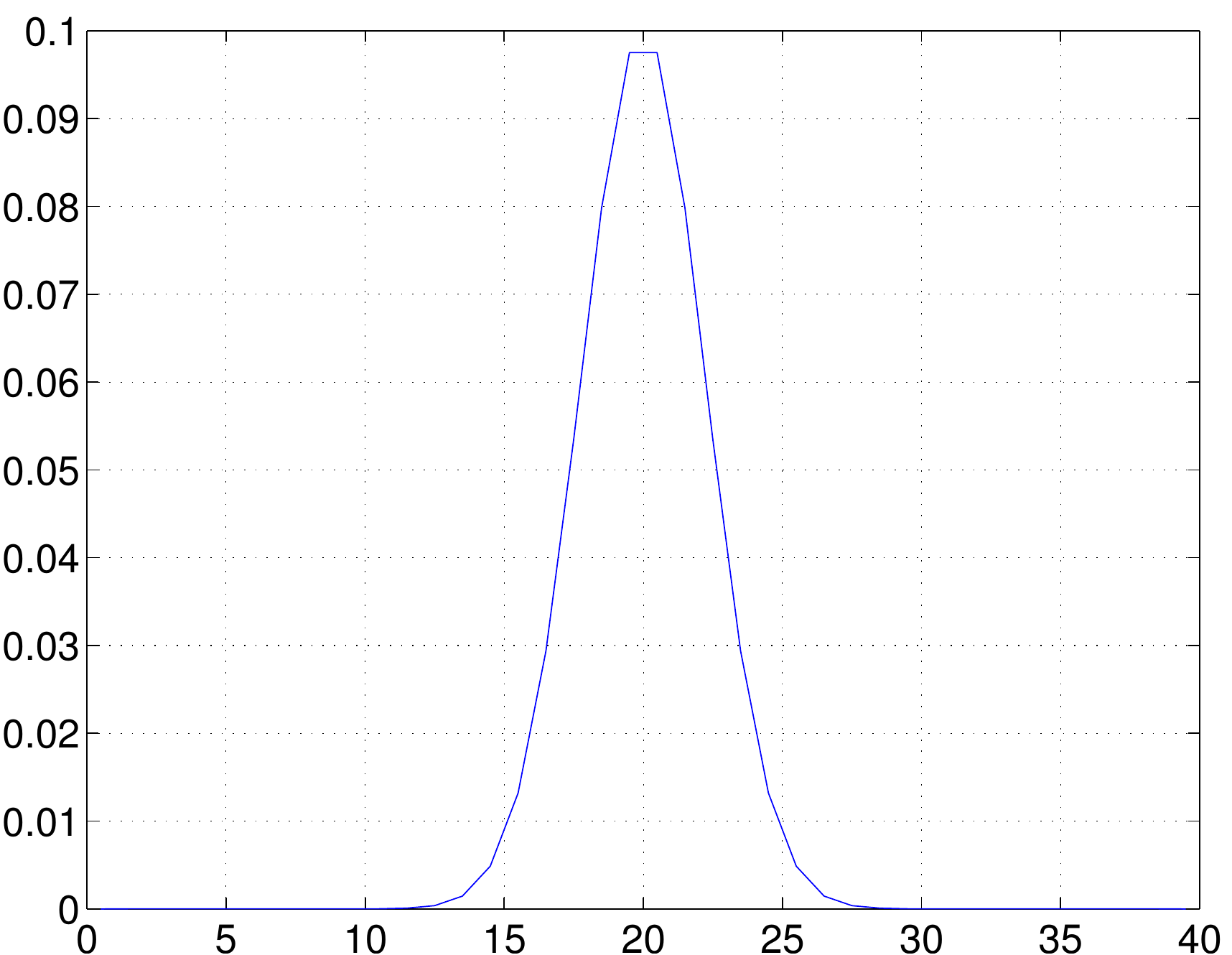}  &   \includegraphics[width=6cm]{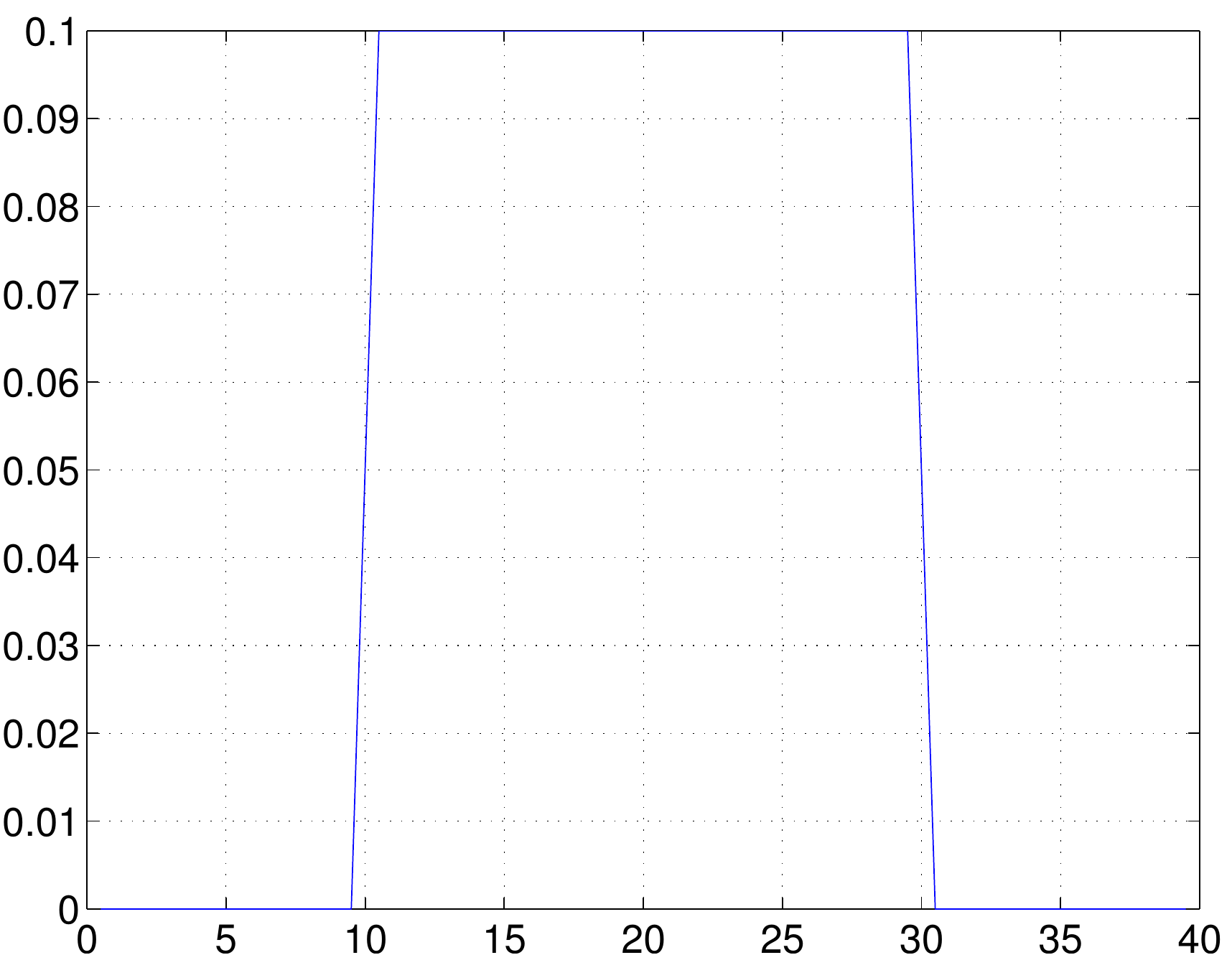}\\
  (a)  Regular initial data                &    (b) Irregular initial data
 \end{tabular}
\caption{\label{fig:growth_homo_init}Polymerization/depolymerization test in homogeneous space: {\it the regular and irregular initial data corresponding to~\eqref{eq:f1d_reg_init} and~\eqref{eq:f1d_irreg_init}.}}
 \end{center}
\end{figure}

Despite its simplistic appearance, this model \eqref{lif-sly} is very intriguing when one is interested on the time asymptotic behavior. An interesting discussion is made in \cite{bibGLT} and the authors highlight specifically the importance of using an anti dissipative numerical scheme in order to avoid numerical artifact (diffusion) and then capture the exact asymptotic profile. In order to test our scheme defined in section \ref{sec:adm} and show its accuracy, we compare the numerical results with those obtained either by a WENO scheme (see \cite{bibCG}) or by the anti-dissipative scheme (see \cite{bibGLT}). Two types of initial distribution functions will be considered: the regular one, which represents the size-density of cells in a normal distribution as in~\eqref{eq:f1d_reg_init},
\begin{equation}\label{eq:f1d_reg_init}
f^{0}(x)=0.1\exp(-0.1(\xi - 20)^2).
\end{equation} 
While the irregular one represents that the  size-density of cells is concentrated around a spot as in~\eqref{eq:f1d_irreg_init},
\begin{equation}\label{eq:f1d_irreg_init}
f^{0}(\xi)=\left\{ \begin{array}{ll}
 1, & \quad \text{if} \quad 10\leq \xi\leq 30, \\
 0, & \quad \text{else}.  \\
 \end{array} \right .
\end{equation} 
These two types of size-density are illustrated in Figure~\ref{fig:growth_homo_init}.

\begin{figure}
\begin{center}
 \begin{tabular}{cc}
  \includegraphics[width=6cm]{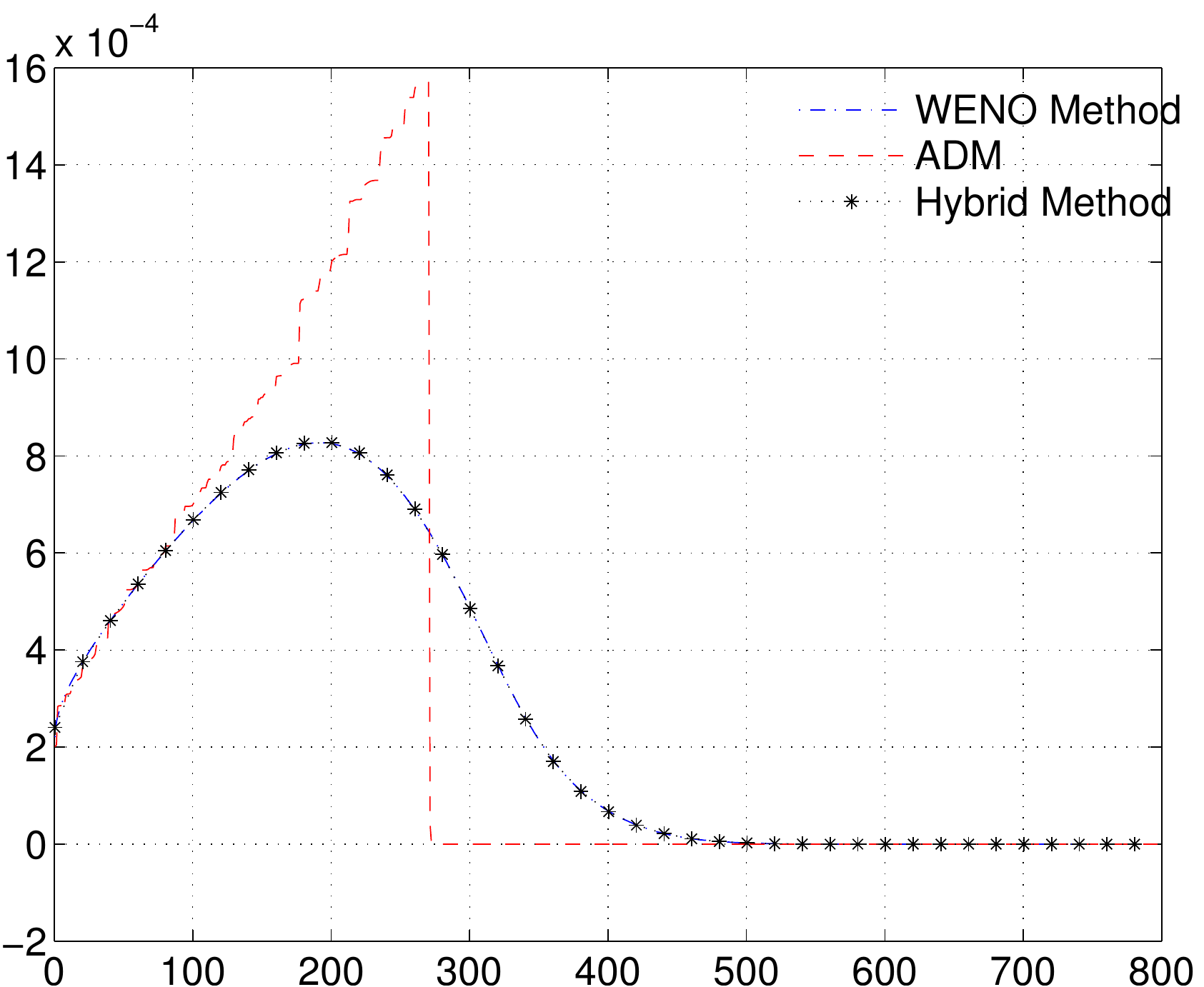}  &   \includegraphics[width=6cm]{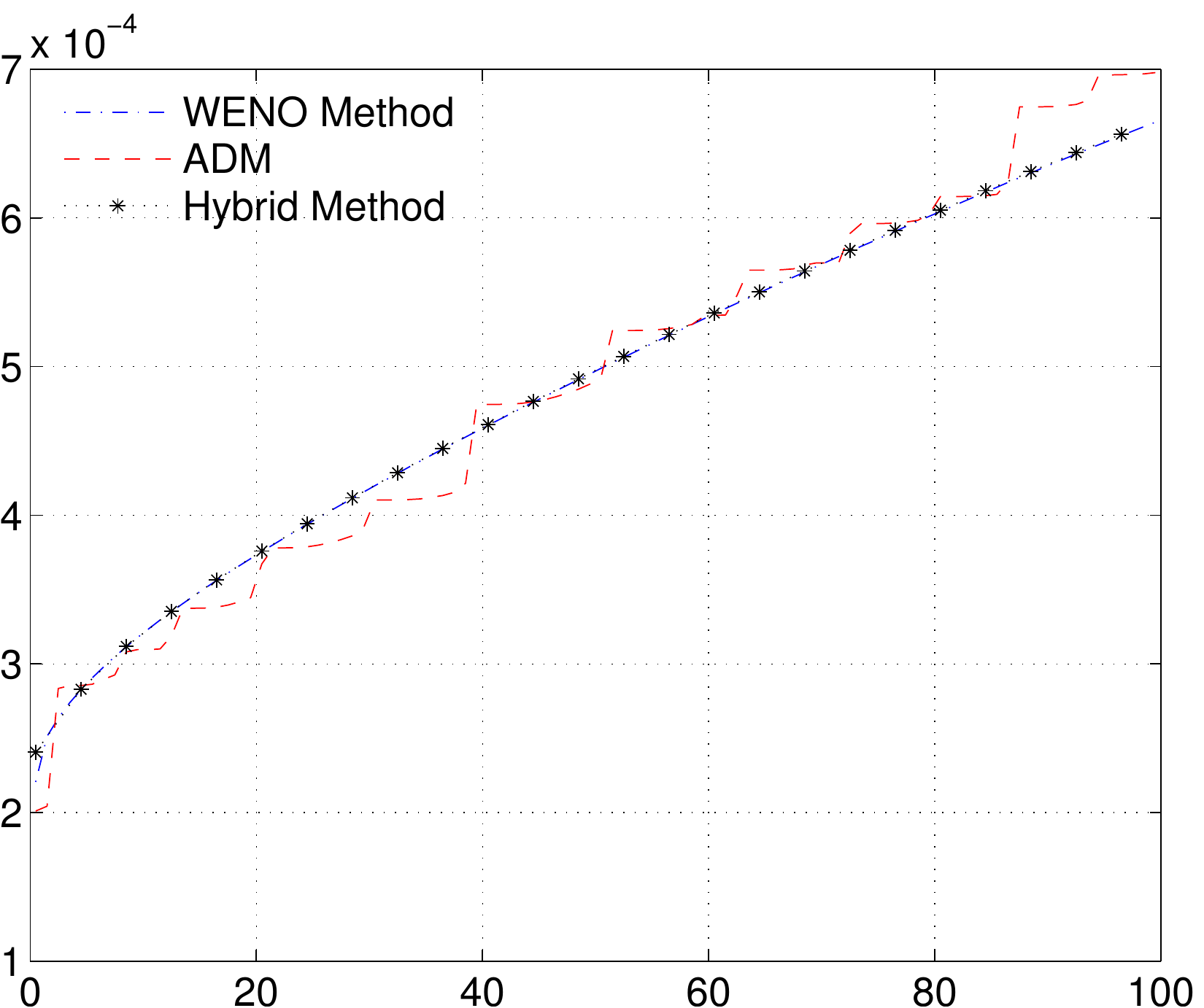}\\
  (a)  $t=1250$                &    (b) Zone, $t=1250$ \\
   \includegraphics[width=6cm]{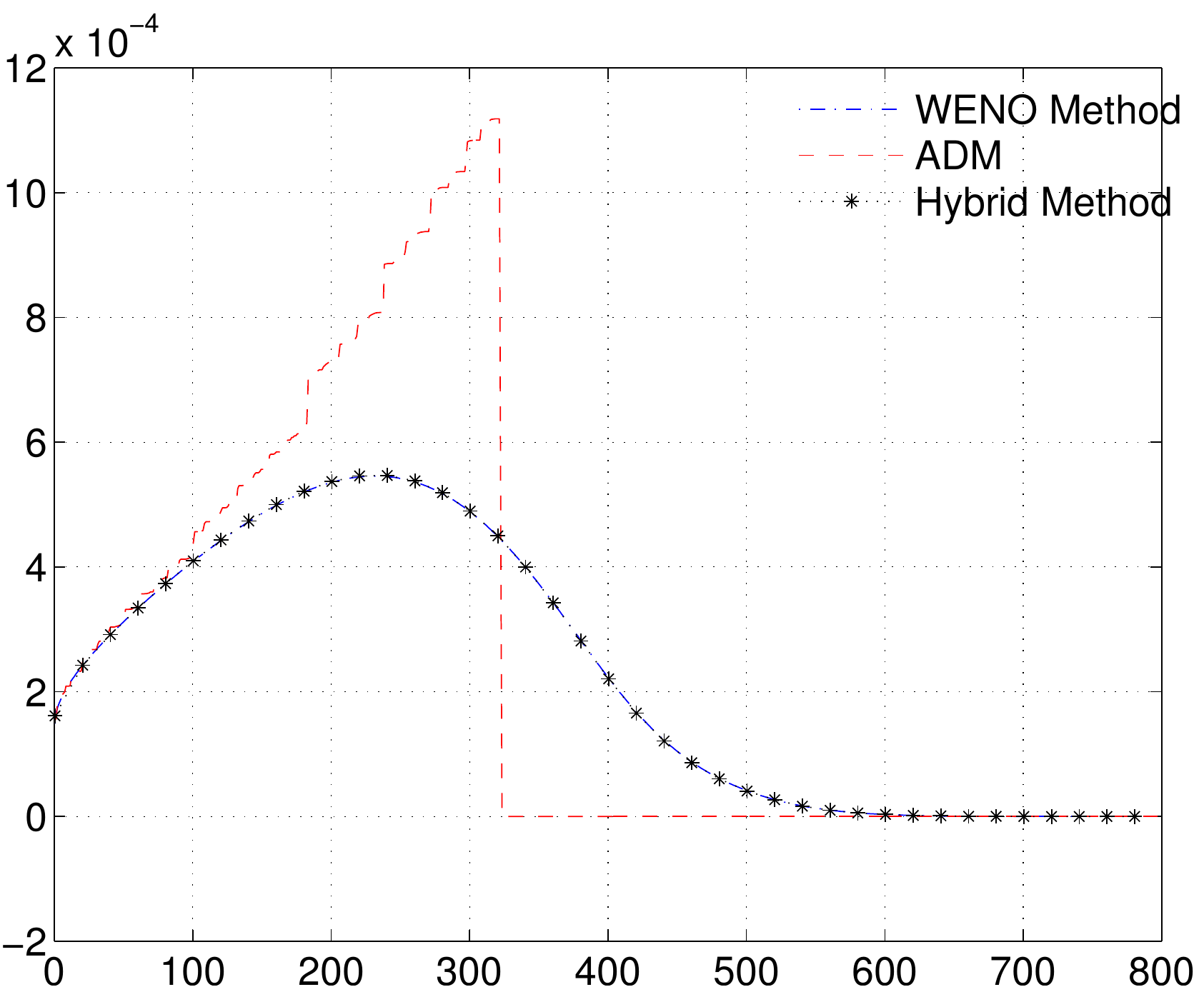}  &   \includegraphics[width=6cm]{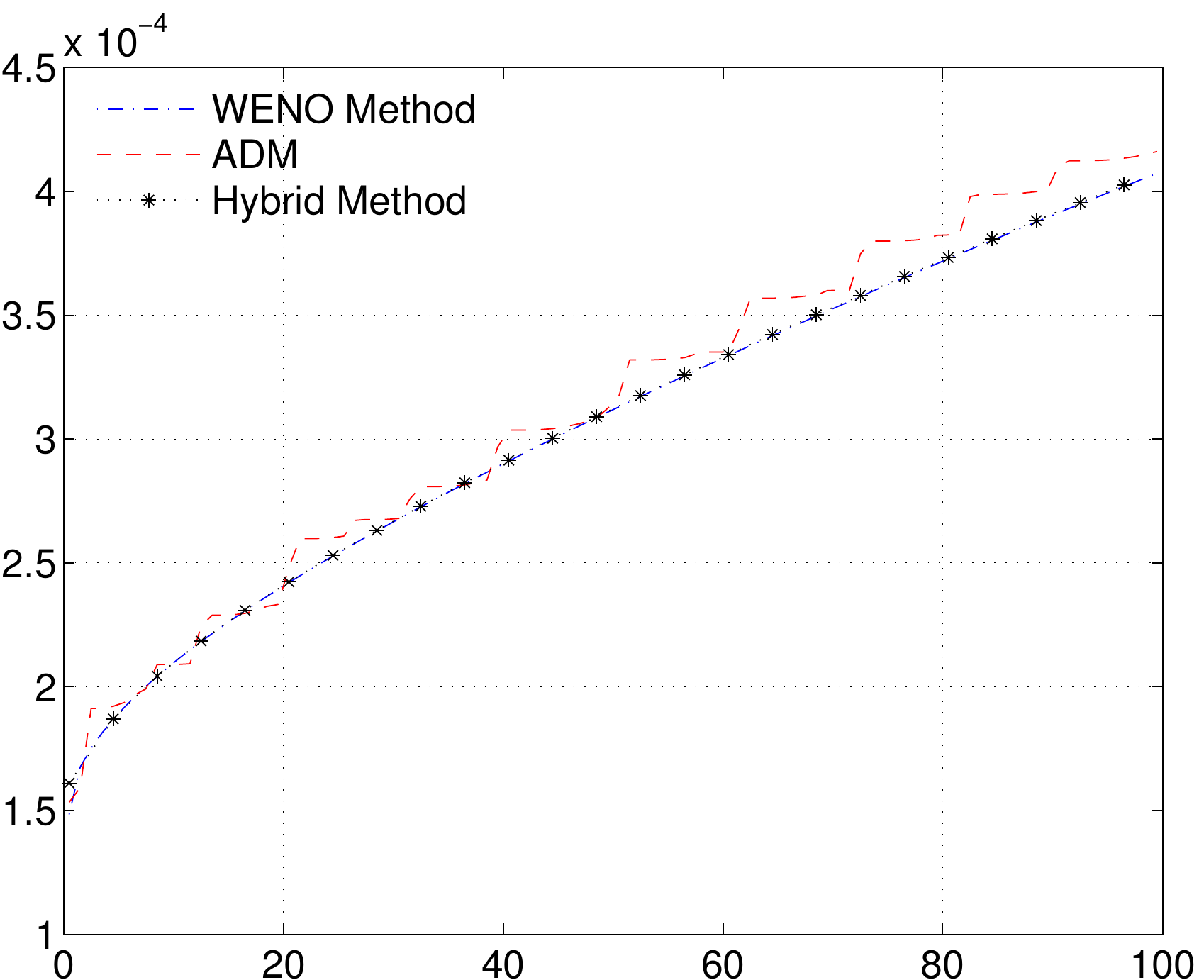}\\
  (c)  $t=1500$                &    (d) Zone, $t=1500$ \\
    \includegraphics[width=6cm]{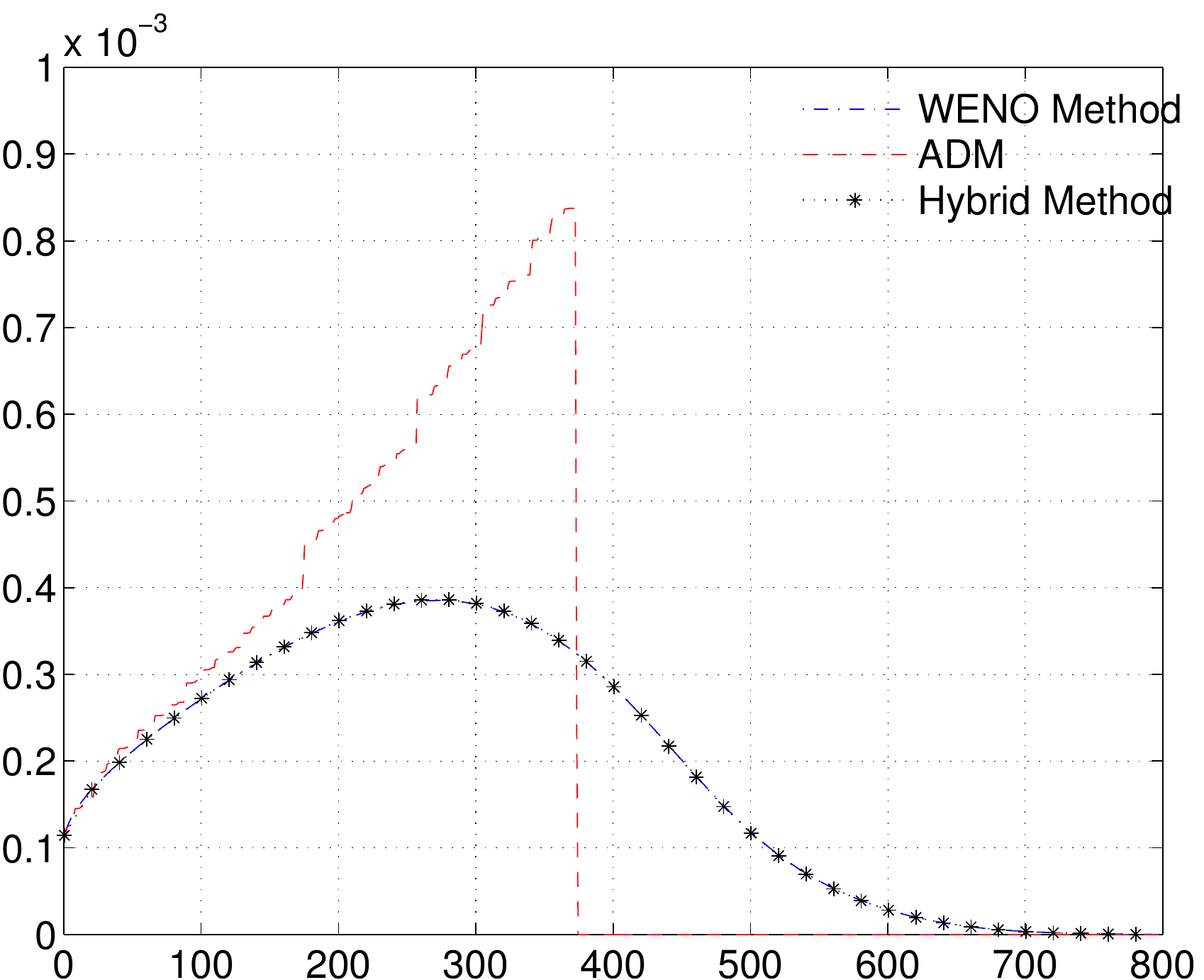}  &   \includegraphics[width=6cm]{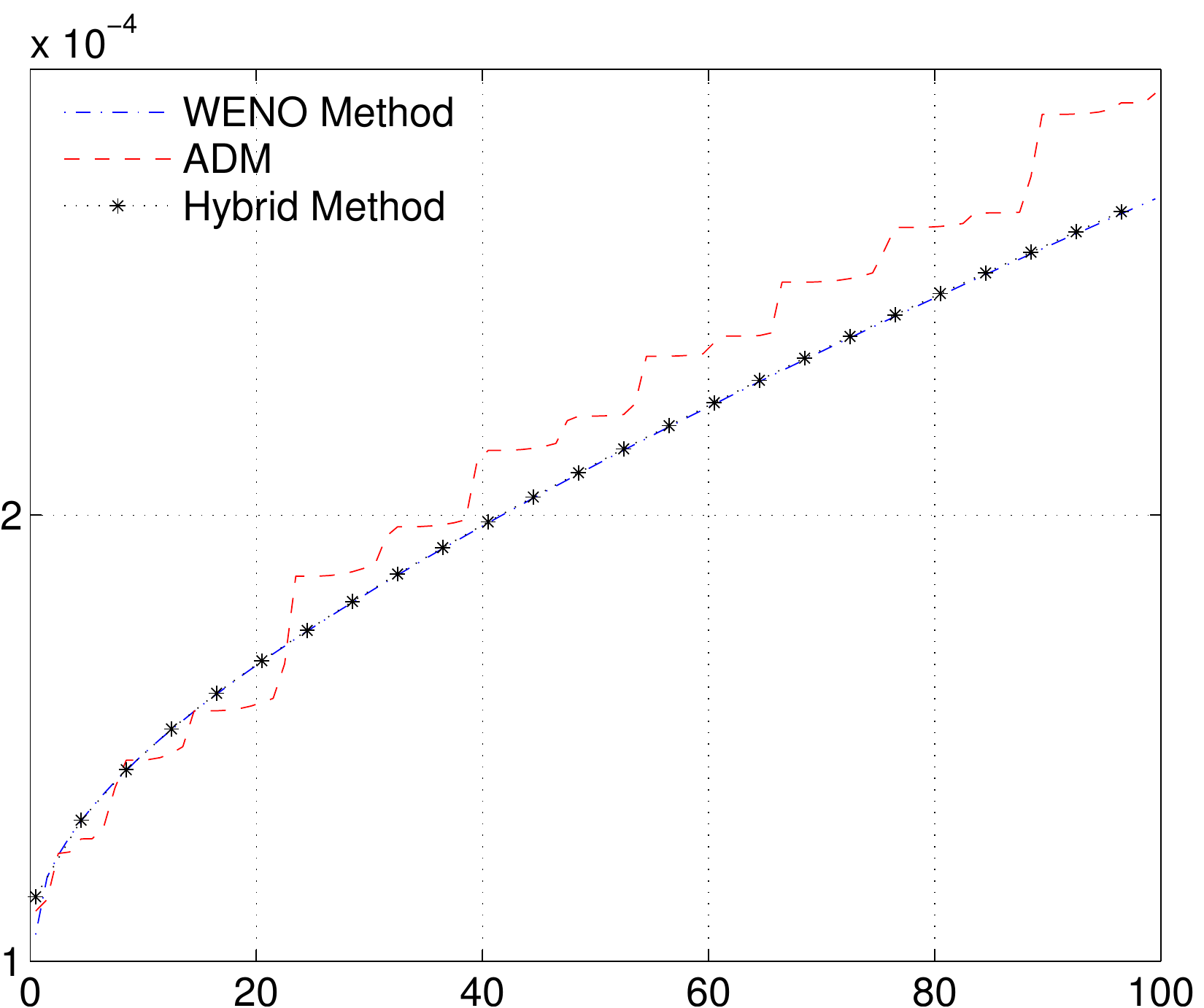}\\
  (e)  $t=1750$                &    (f) Zone, $t=1750$ \\
   \includegraphics[width=6cm]{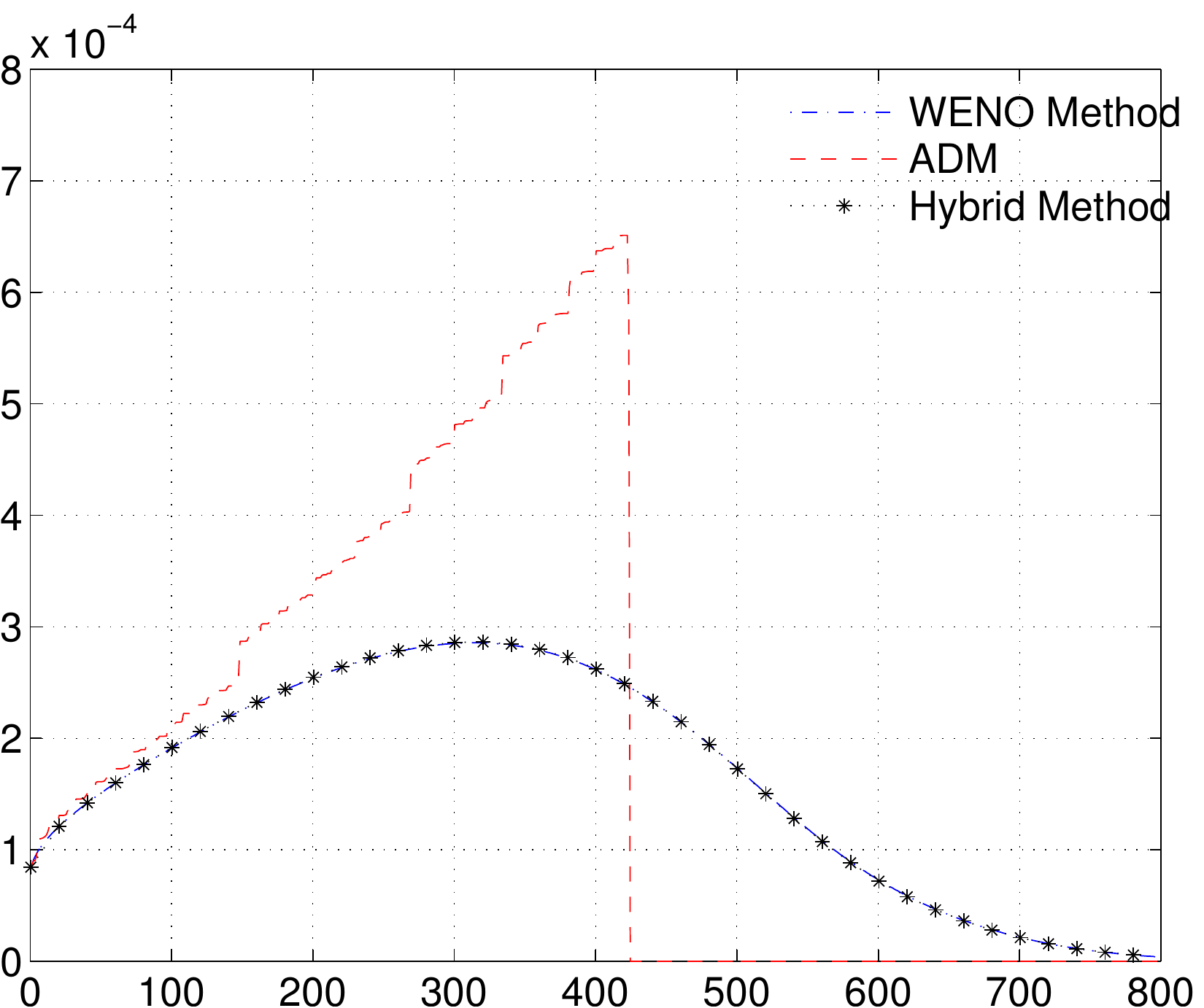}  &   \includegraphics[width=6cm]{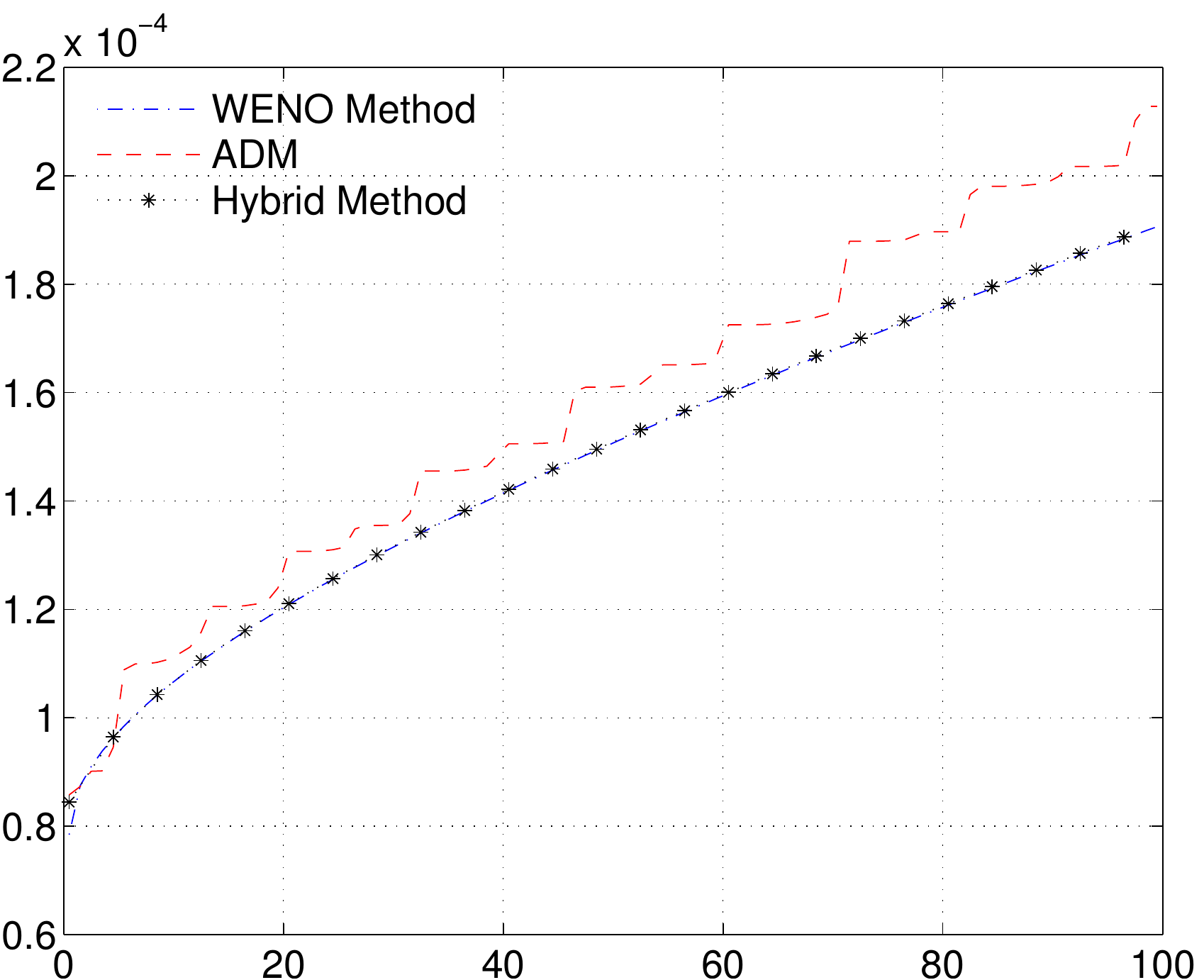}\\
  (g)  $t=2000$                &    (h) Zone, $t=2000$ \\
 \end{tabular}
\caption{\label{fig:growth_homo_reg}Polymerization/depolymerization test in homogeneous space: {\it Plot solutions of the  equation~\eqref{lif-sly} with the regular initial data~\eqref{eq:f1d_reg_init}. Mesh size is $n_x=800$, $\Delta t=1/8$, $CFL\approx0.12$, $\rho=41$.}}
 \end{center}
\end{figure}

We first consider the numerical results with the regular data~\eqref{eq:f1d_reg_init} (see Figure~\ref{fig:growth_homo_reg}). Our hybrid scheme has the same results as the classical WENO scheme. More precisely, a very smooth profile is formed, which moves towards cells of large size. This observation has a good agreement with the physical behavior of transport equation which does not modify the general shape of the initial distribution function. In biological point of view, the observation show the growth of big size cells at the expense of smaller ones. This competition between big size cells and smaller size ones is due to the fact that at any time there exists a critical size $\xi_{crit}=\frac{1}{c(t)^3}$ such that the velocity vanish; then cells of size $\xi>\xi_{crit}$ grow while cells of size $\xi<\xi_{crit}$ shrink. This competition phenomenon is well known under the name ``Ostwald ripening'' ~\cite{bibNieth, bibOstw1, bibOstw2, bibOstw3}. At contrast, the  anti-dissipative method forms a very sharp front, which is caused by the Downwind flux and it can not be significantly improved by mesh refinement. Even in the zone view for small size particles (see the right column of Figure~\ref{fig:growth_homo_reg}), the  anti-dissipative method generates ``stairs" looking like oscillations.

\begin{figure}
\begin{center}
 \begin{tabular}{cc}
  \includegraphics[width=6cm]{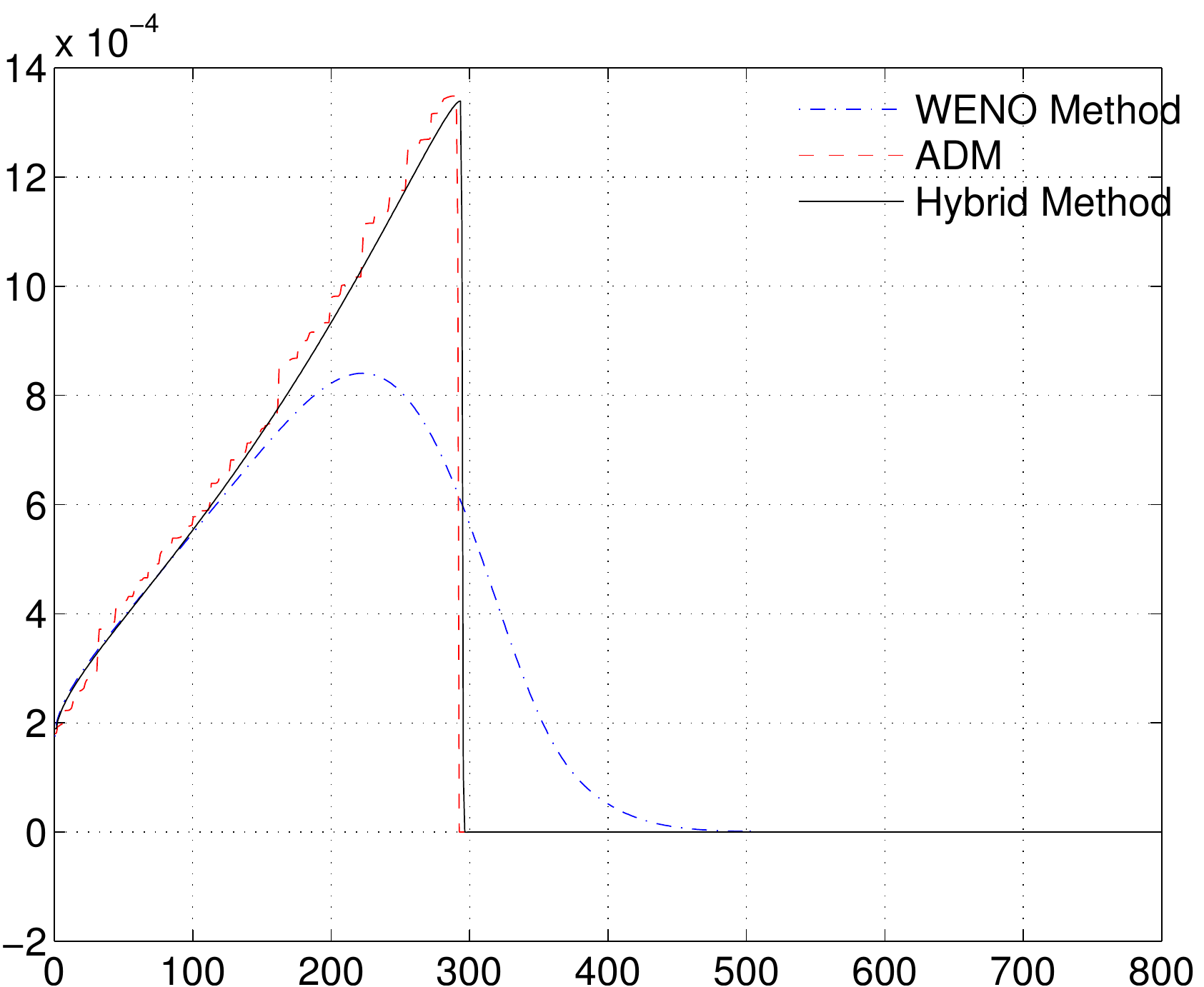}  &   \includegraphics[width=6cm]{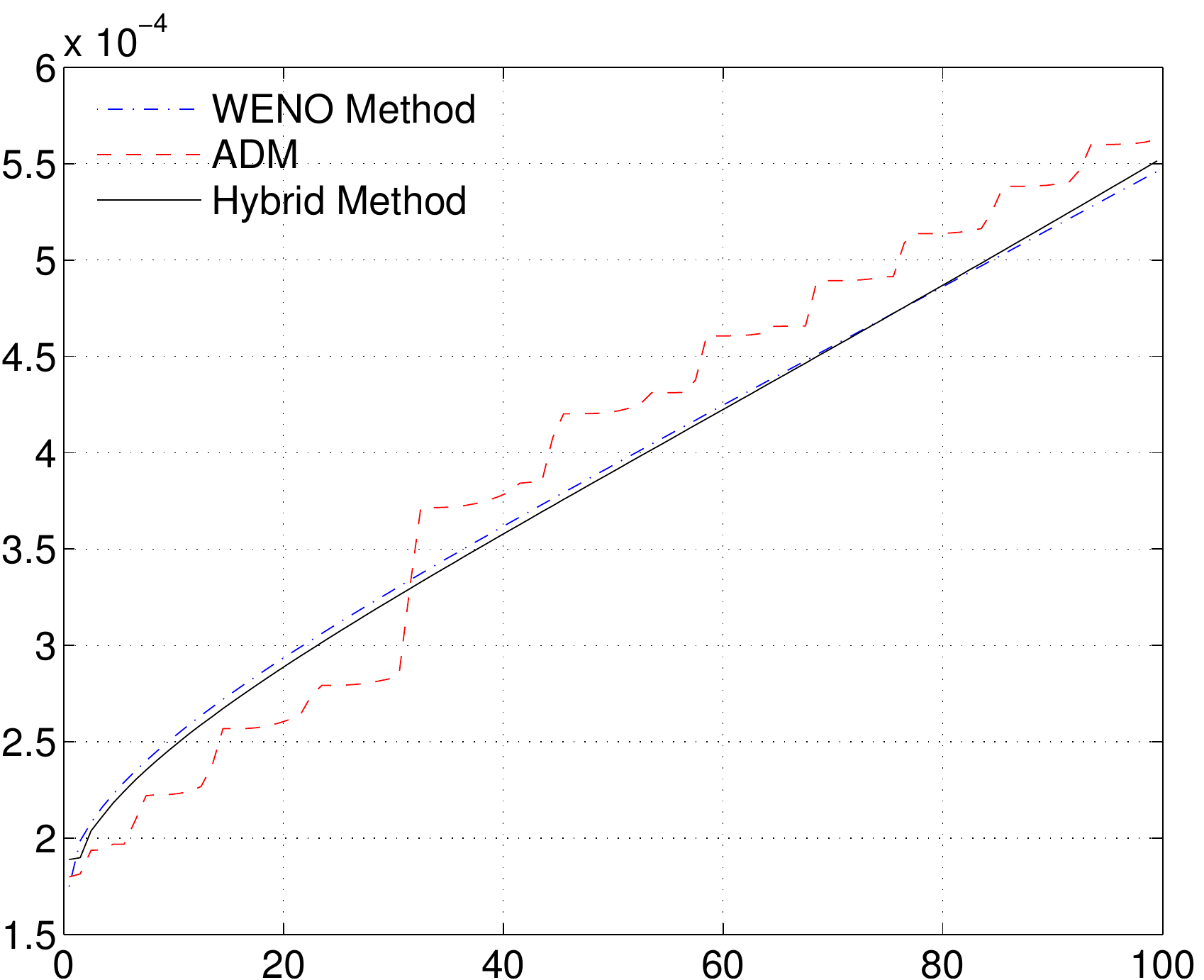}\\
  (a)  $t=1250$                &    (b) Zone, $t=1250$ \\
   \includegraphics[width=6cm]{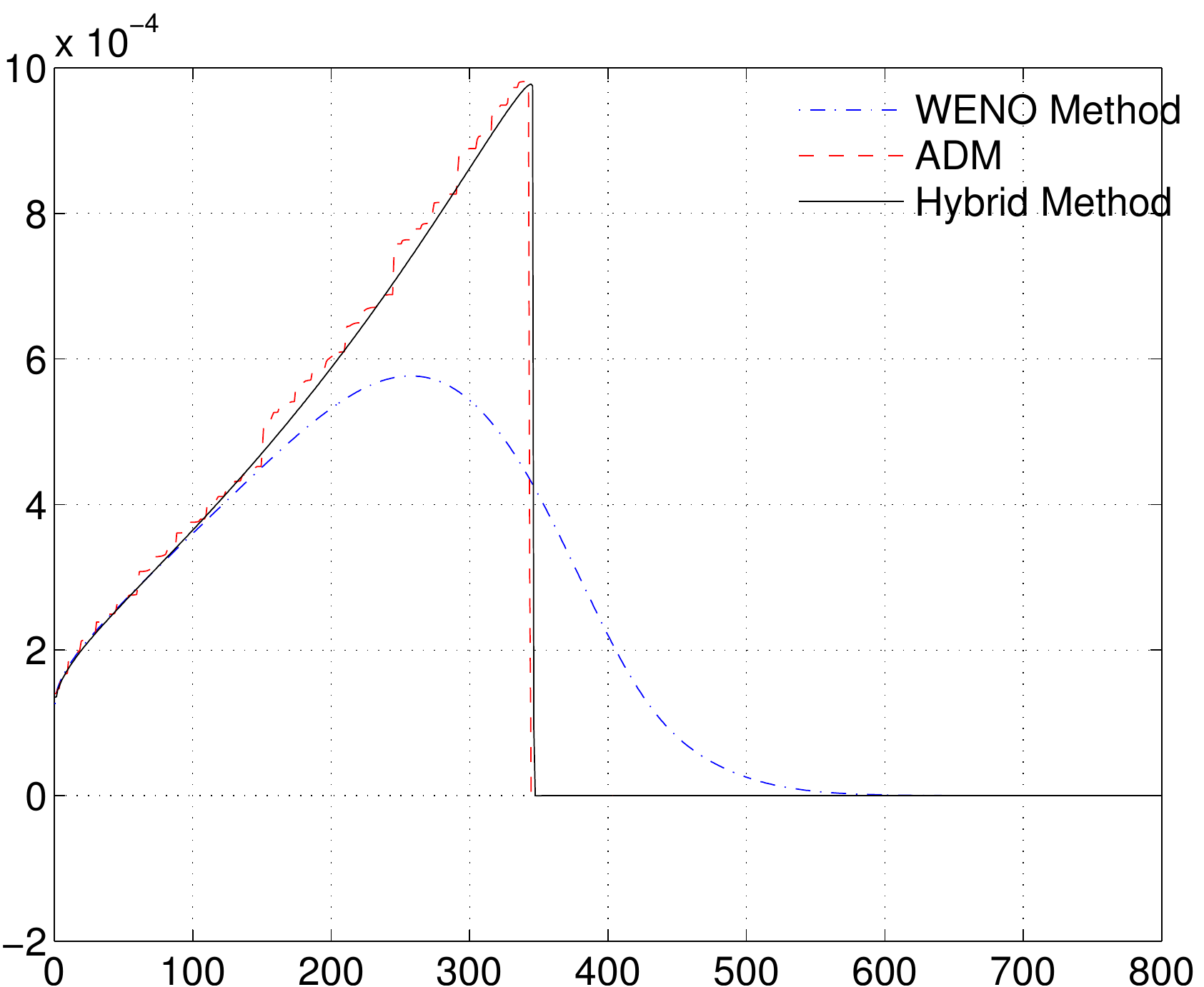}  &   \includegraphics[width=6cm]{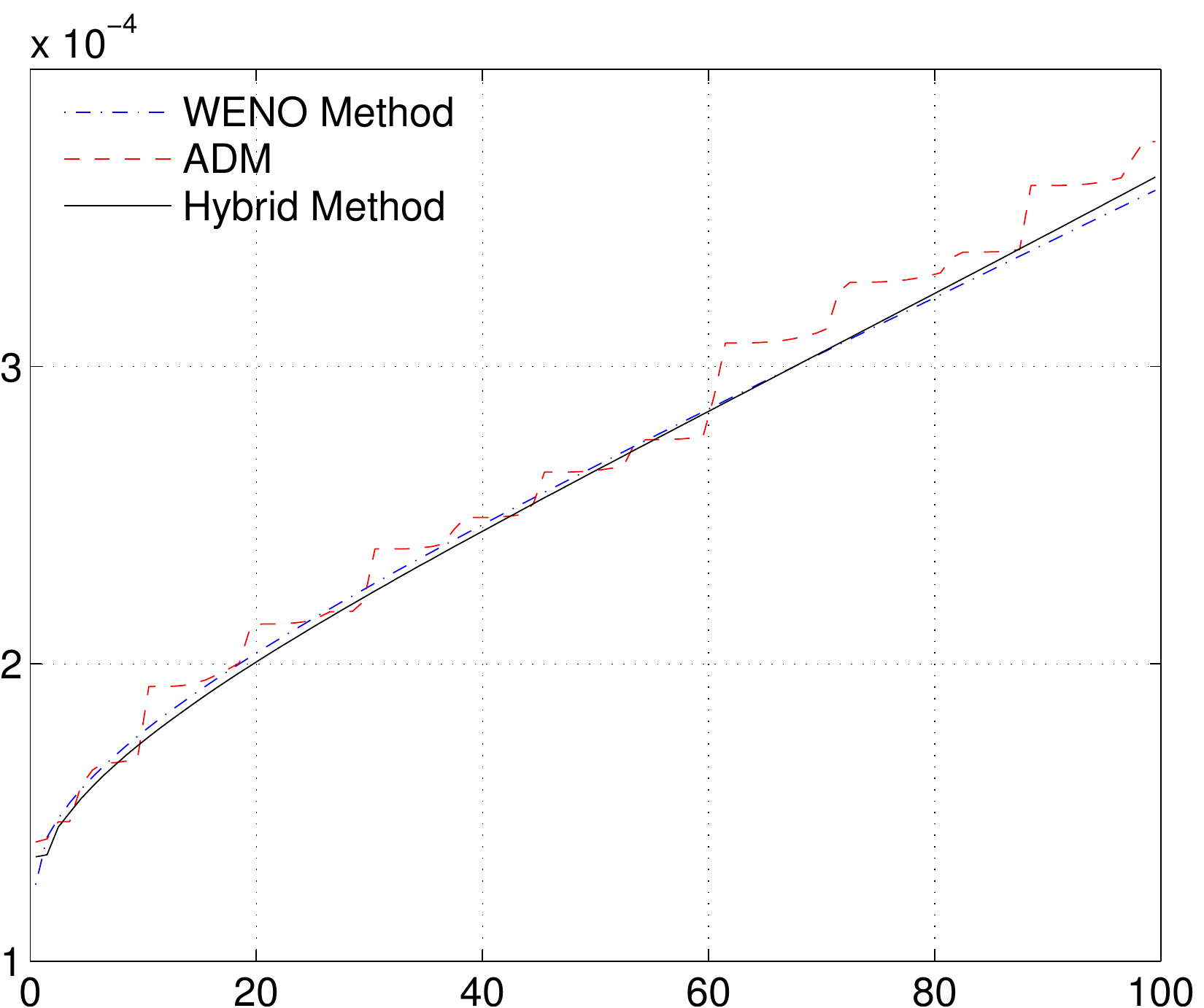}\\
  (c)  $t=1500$                &    (d) Zone, $t=1500$ \\
    \includegraphics[width=6cm]{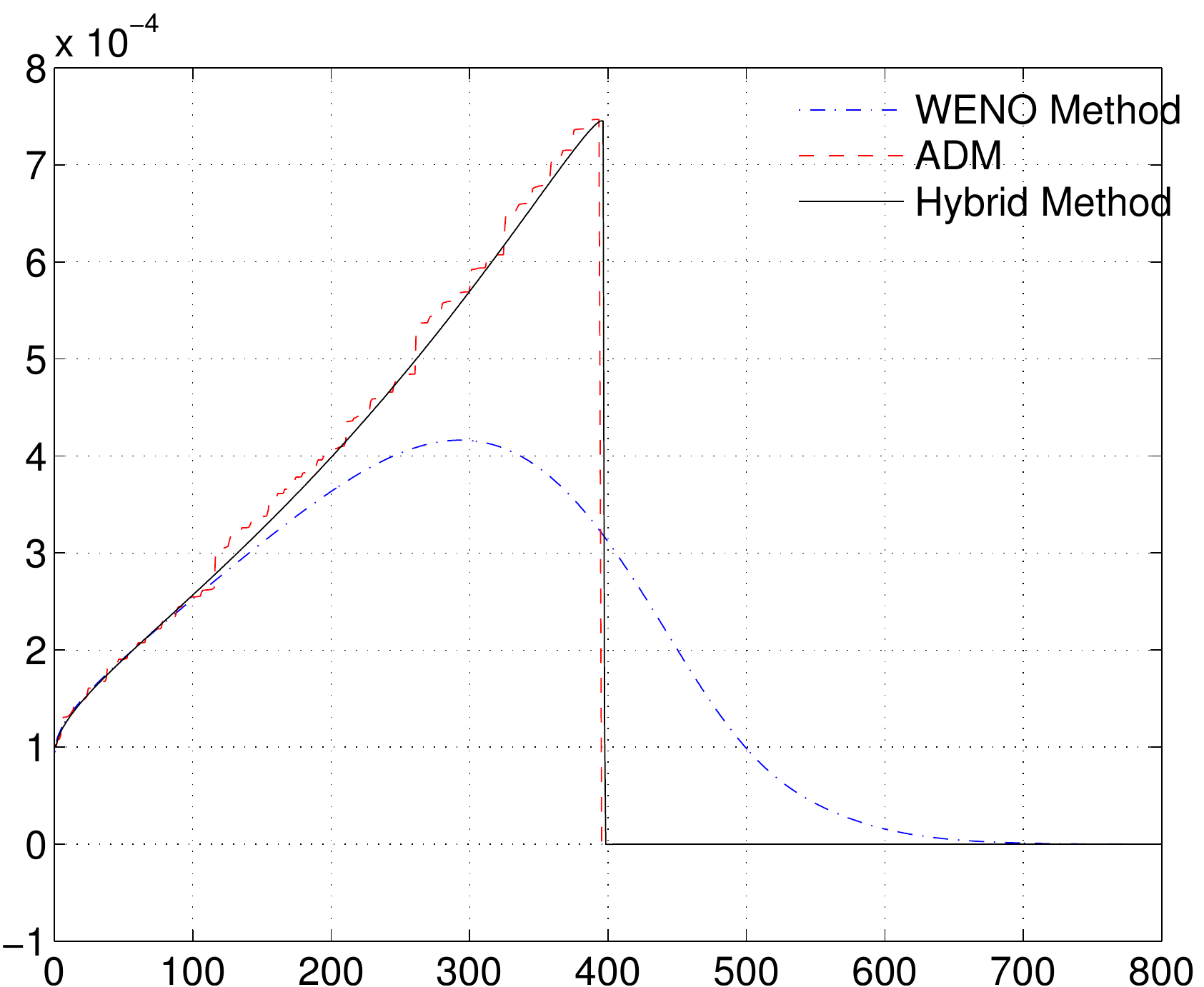}  &   \includegraphics[width=6cm]{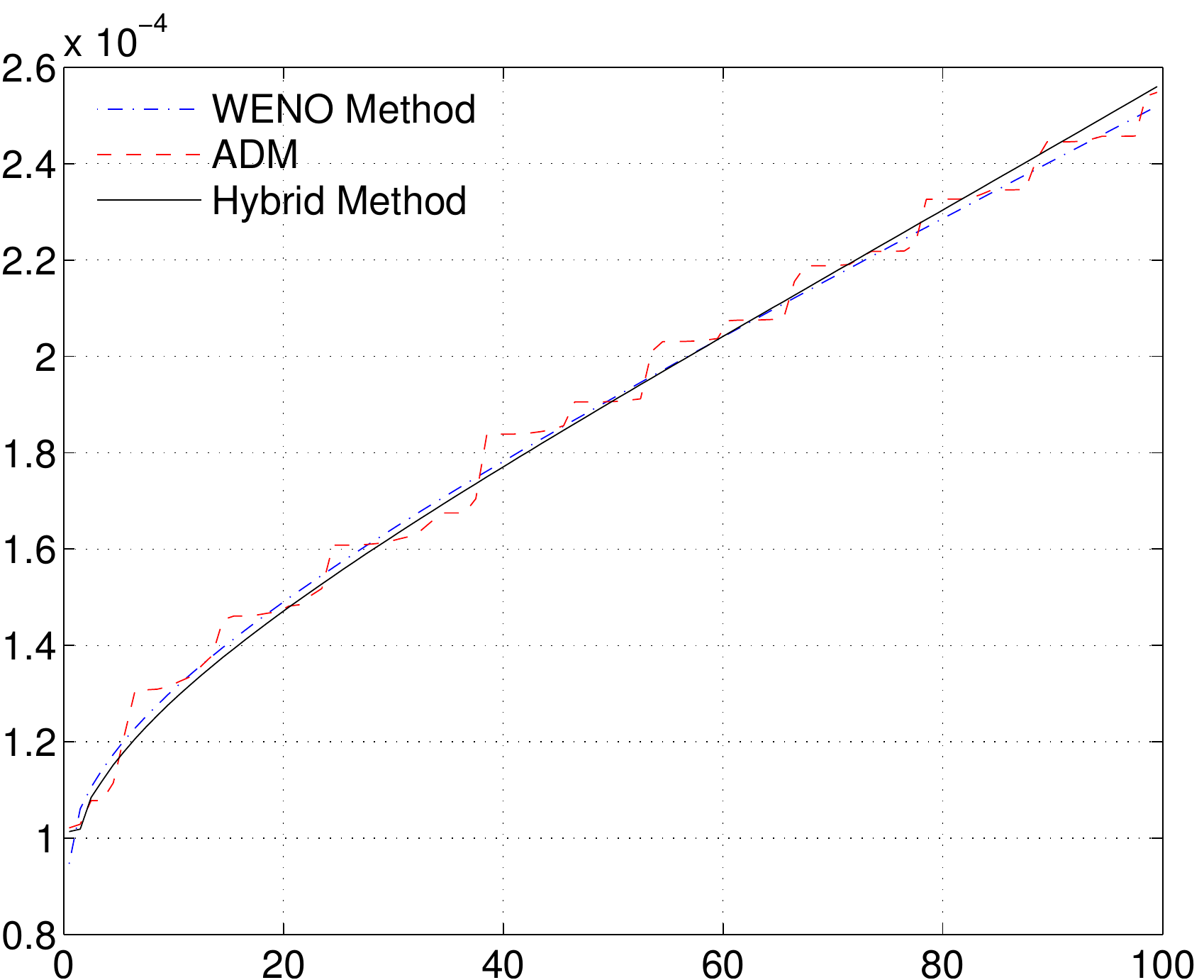}\\
  (e)  $t=1750$                &    (f) Zone, $t=1750$ \\
   \includegraphics[width=6cm]{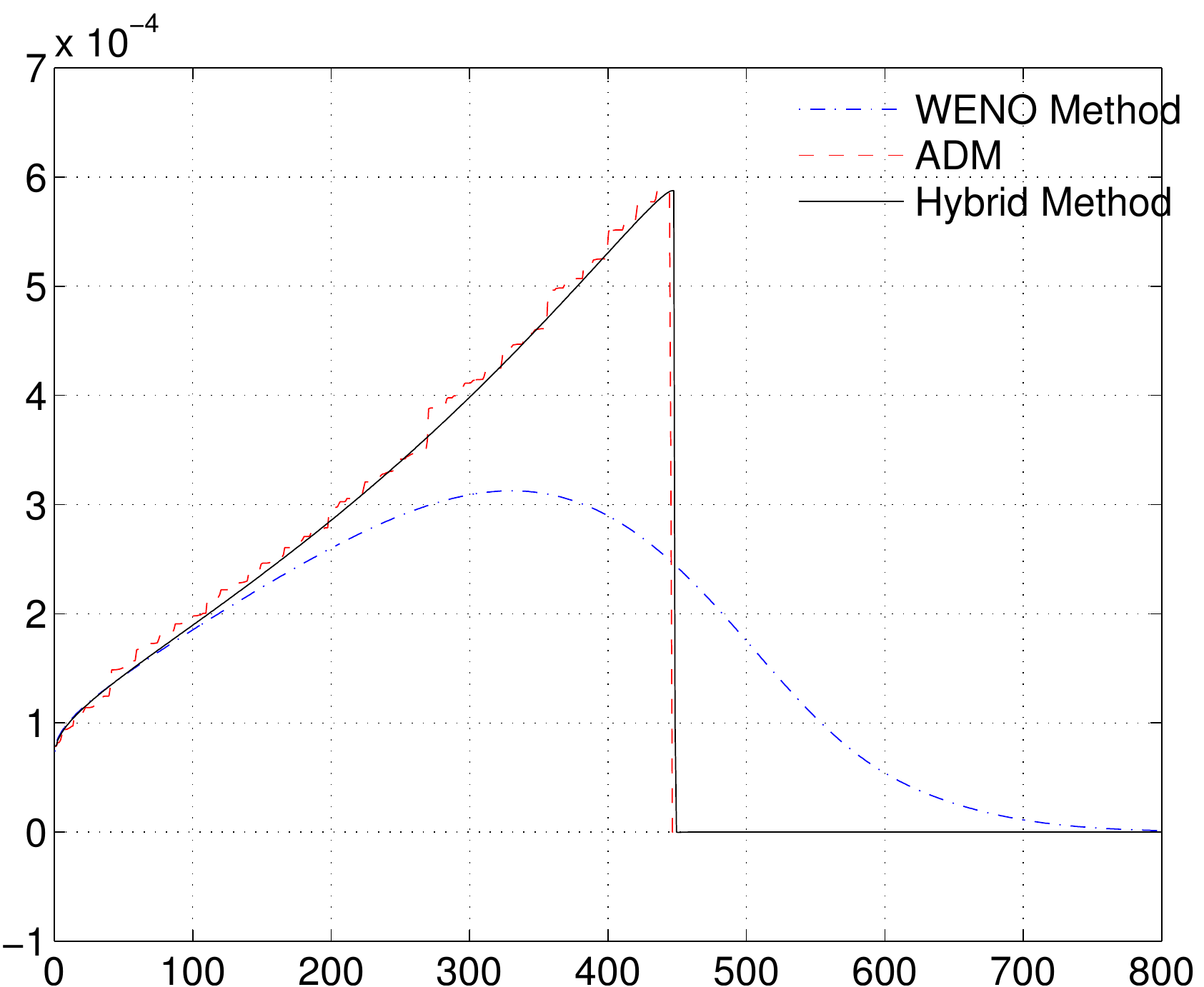}  &   \includegraphics[width=6cm]{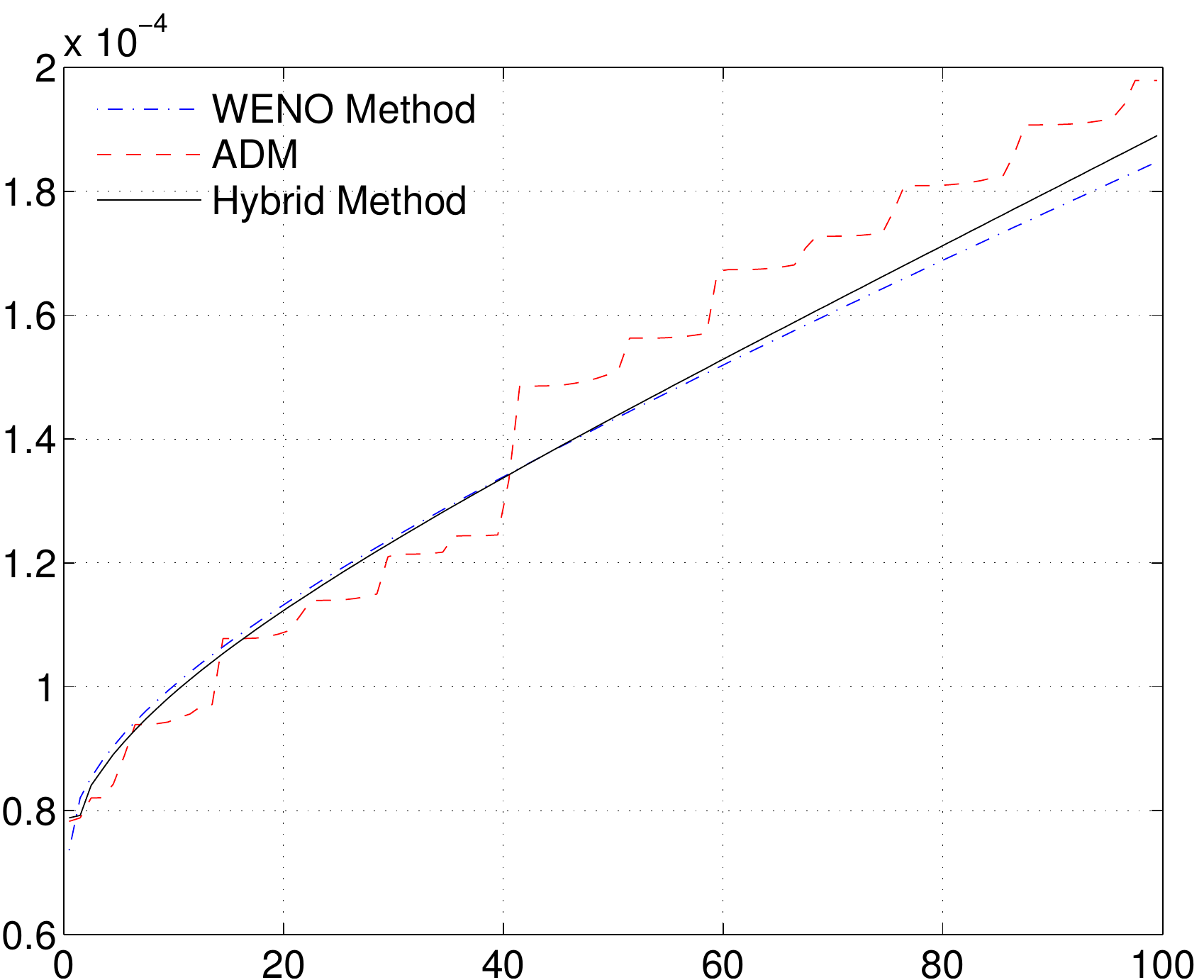}\\
  (g)  $t=2000$                &    (h) Zone, $t=2000$ \\
 \end{tabular}
\caption{\label{fig:growth_homo}Polymerization/depolymerization test in homogeneous space: {\it Plot solutions of the  equation~\eqref{lif-sly} with the irregular initial data~\eqref{eq:f1d_irreg_init}. Mesh size is $n_x=800$, $\Delta t=1/8$, $CFL\approx0.12$, $\rho=41$.}}
 \end{center}
\end{figure}

Then, with  the irregular data~\eqref{eq:f1d_irreg_init}, as pictured at the figure \ref{fig:growth_homo}, we first notice that the irregular data is regularized by the numerical diffusion by using the WENO scheme.  Indeed, I. M. Lifshitz and V. V. Slyozov in their original paper~\cite{bibLS} conjecture that the asymptotic solution is in the same form despite the initial data (compare with the left column of Figure~\ref{fig:growth_homo_reg}). However, as pointed out in~\cite{bibGLT}, the  irregular data does have an influence for long term solution.  Our scheme as the anti-dissipative one preserves well the front propagation at long time (no numerical diffusion) what is physically very important when dealing with transport equation and can be interpreted in biology as the fact that very stiff localized population of cells must remain stiff localized if they are subject to only transport process.  Moreover, our scheme corrects well the ``stairs" looking like oscillations in the numerical solution of the anti-dissipative method (see the right column of Figure~\ref{fig:growth_homo}). 



\subsubsection{Non-space-homogeneous polymerization/depolymerization type model}
Here we are interested on the numerical behavior of the non-homogeneous system \eqref{eq1}-\eqref{eq2} with the choice $a(\xi)=\xi^{1/3}$, $b=1$, that we rewrite as follows 
\begin{equation}
\left\{ \begin{array}{ll}\label{eq:lif-sly1d1d}
\displaystyle\frac{\partial}{\partial t} f(t,x,\xi) +\frac{\partial}{\partial \xi} \bigl((\xi^{1/3}c(t,x)-1)f(t,x,\xi)\bigr) =0 ,\quad\,t\geq0,\,x\in\Omega,\,\xi\geq0,\\
\\
\partial_t\left(c(t,x) + \displaystyle\int_{0}^{\infty} \xi f(t,x,\xi) d\xi\right)  = \Delta_x c(t,x),\quad\,t\geq0,\,x\in\Omega.\\
\end{array} \right.
\end{equation}  
The diffusion equation is endowed with homogeneous Neumann boundary condition
\begin{equation*}
 \partial_\nu c = \nabla c\cdot\nu = 0, \,\text{on }\partial\Omega. 
\end{equation*}
Here we choose $\Omega \subset \mathbb{R}$ means $1D$ in space-variable. The system \eqref{eq:lif-sly1d1d} describes the immersion of a population of cells in a culture of micro-organisms 
(nutrients) that are subjected to diffusion equation while de size density repartition of cells is parametrized by the space position $x$.

\begin{figure}
\begin{center}
 \begin{tabular}{cc}
  \includegraphics[width=7cm]{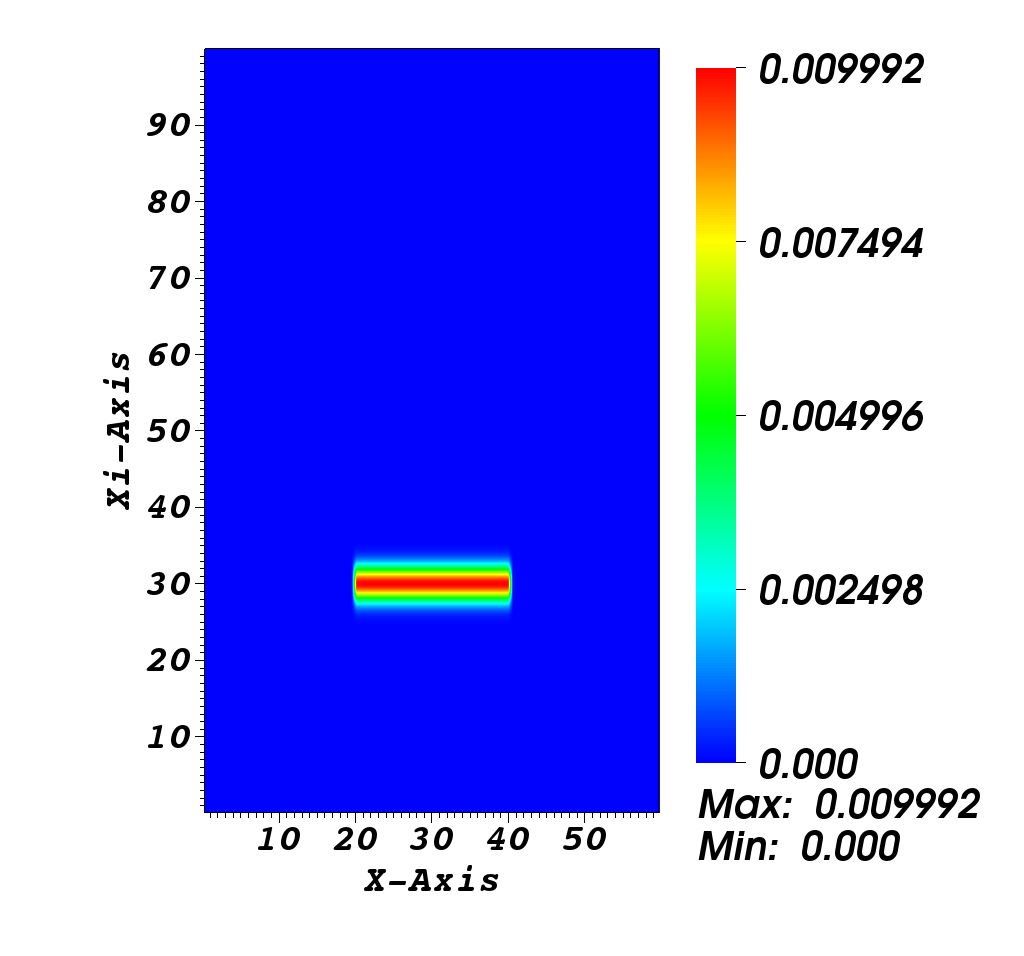}  &   \includegraphics[width=7cm]{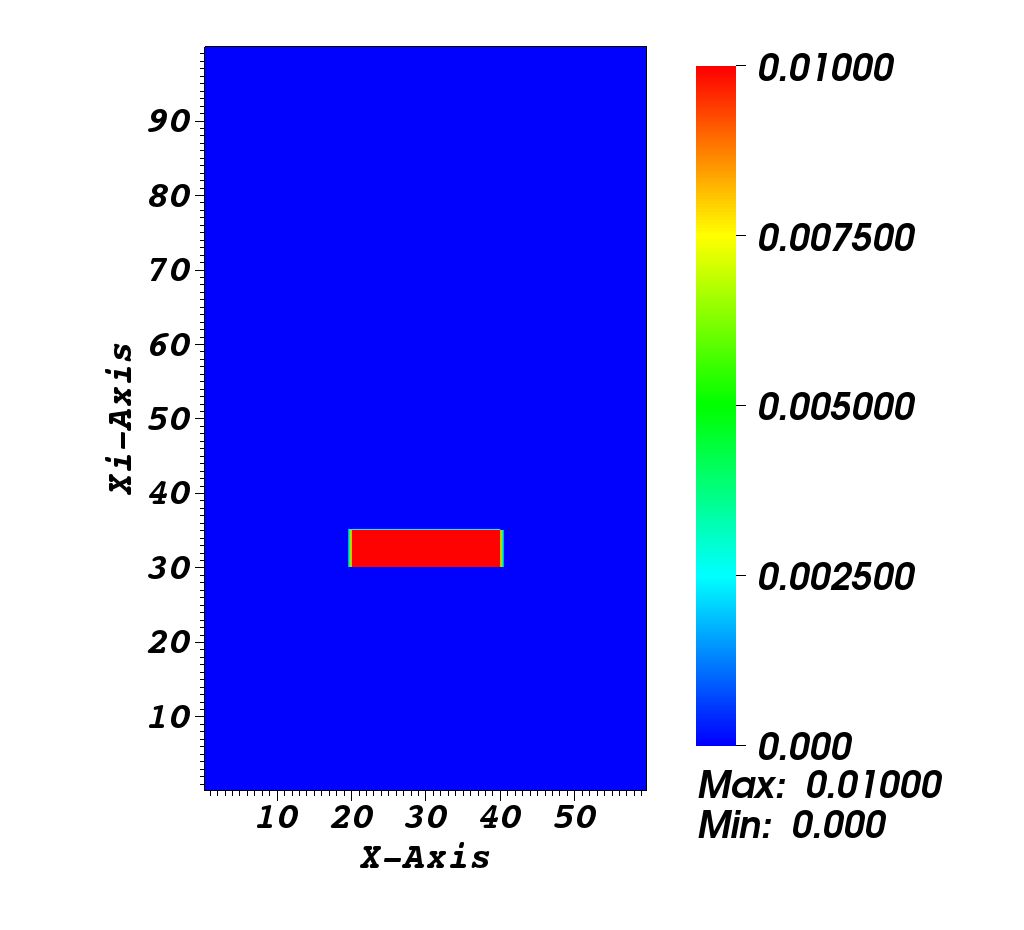} \\
  (a)  Regular initial size distribution function                &    (b) Irregular initial  size distribution function 
 \end{tabular}
 \begin{tabular}{c}
  \includegraphics[width=7cm]{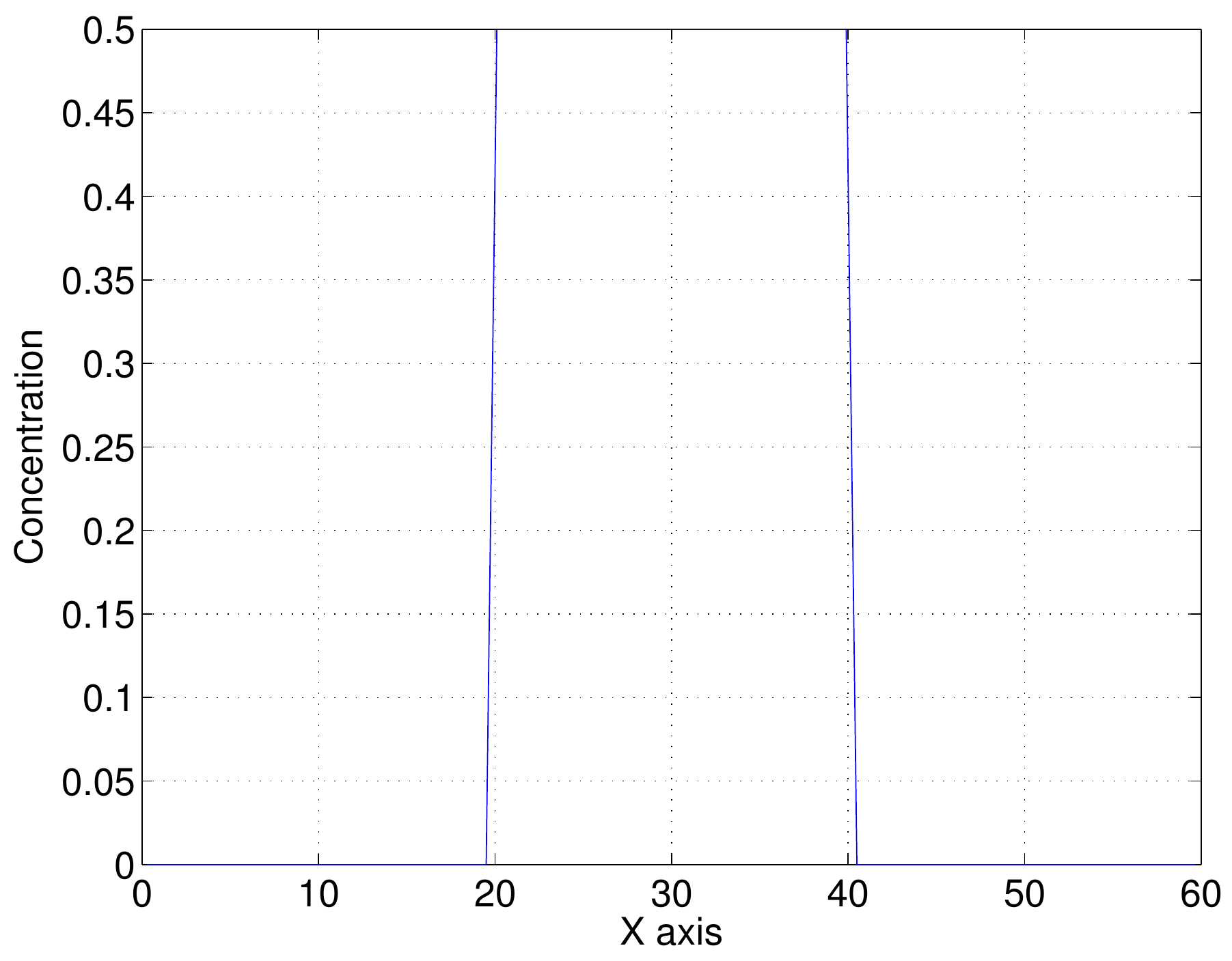}  \\
     (c) Initial concentration
 \end{tabular}
\caption{\label{fig:LS2D_IC}Polymerization/depolymerization test in non-homogeneous space: {\it the  initial data corresponding to~\eqref{eq:lif-sly1d1d_init_concentration}-\eqref{eq:lif-sly1d1d_init_irreg}.} }
 \end{center}
\end{figure}

The numerical simulations are performed in the slab $x\in[0,60]$. The size variable is truncated to $\xi\in[0,100]$. The initial concentration is defined by
\begin{equation}\label{eq:lif-sly1d1d_init_concentration}
c(t=0,x) = 0.5\mathbb{I}_{ x\in[20,40] }.
\end{equation}
As we did in the previous test, two types of size-density are considered: a regular one
 \begin{equation} \label{eq:lif-sly1d1d_init_reg}
   f(t=0,x,\xi) = 0.01\exp(-0.2(\xi-30)^2) \mathbb{I}_{ x\in[20,40] }.
   \end{equation}          
and an irregular one
\begin{equation} \label{eq:lif-sly1d1d_init_irreg}
   f(t=0,x,\xi) = 0.01 \mathbb{I}_{ x\in[20,40] }\times \mathbb{I}_{ \xi\in[30,35] }.
\end{equation}
The initial data is presented in Figure~\ref{fig:LS2D_IC}.

Let us briefly describe the algorithm for the system~\eqref{eq:lif-sly1d1d}. The third order Runge-Kutta method is used as the time discretization for the advection equation of the size-density of macro-particles. In the reconstruction of numerical flux, both the WENO scheme and our hybrid scheme are applied for the purpose of comparison. Then for the diffusion equation of monomers, the crank-Nicolson method is used for time discretization. A classical second-order finite volume method is applied for the lapacian operator, and the integral is approximated by the classical Simpson's rule.

\begin{figure}
\begin{center}
 \begin{tabular}{ccc}
  \includegraphics[width=4.5cm]{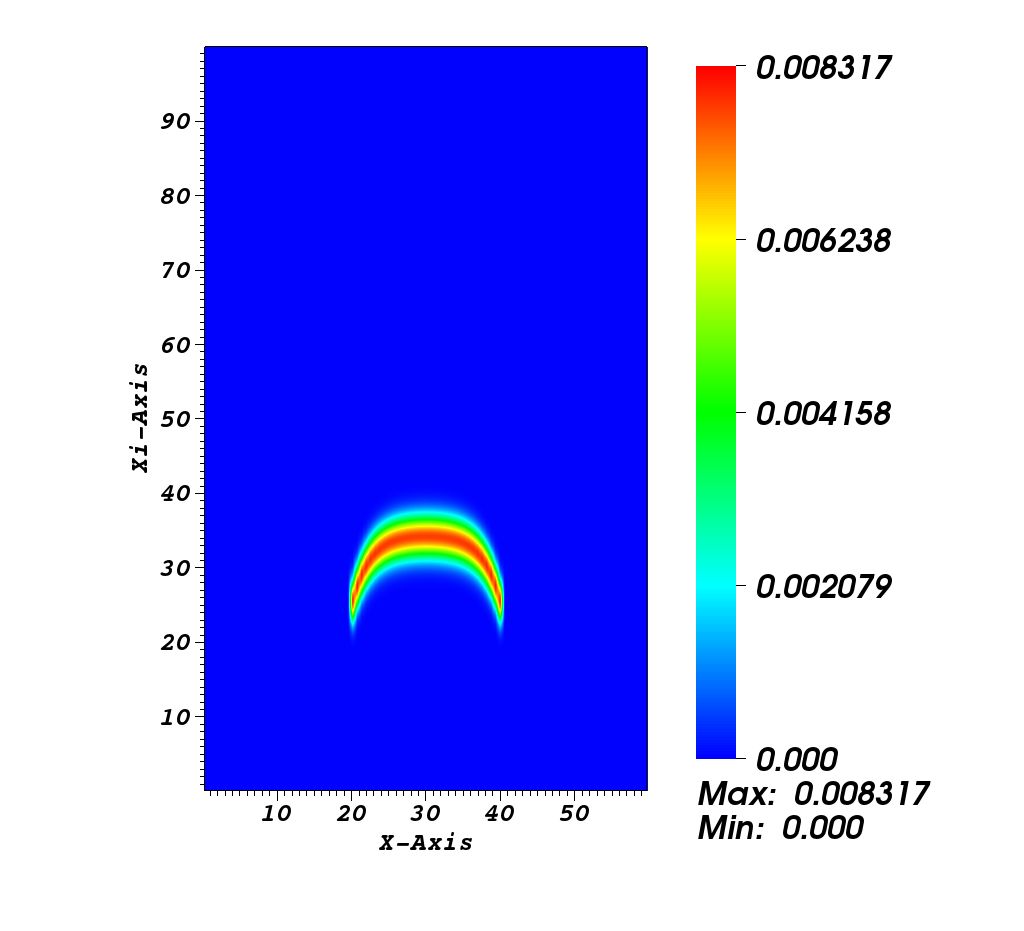}  &   \includegraphics[width=4.5cm]{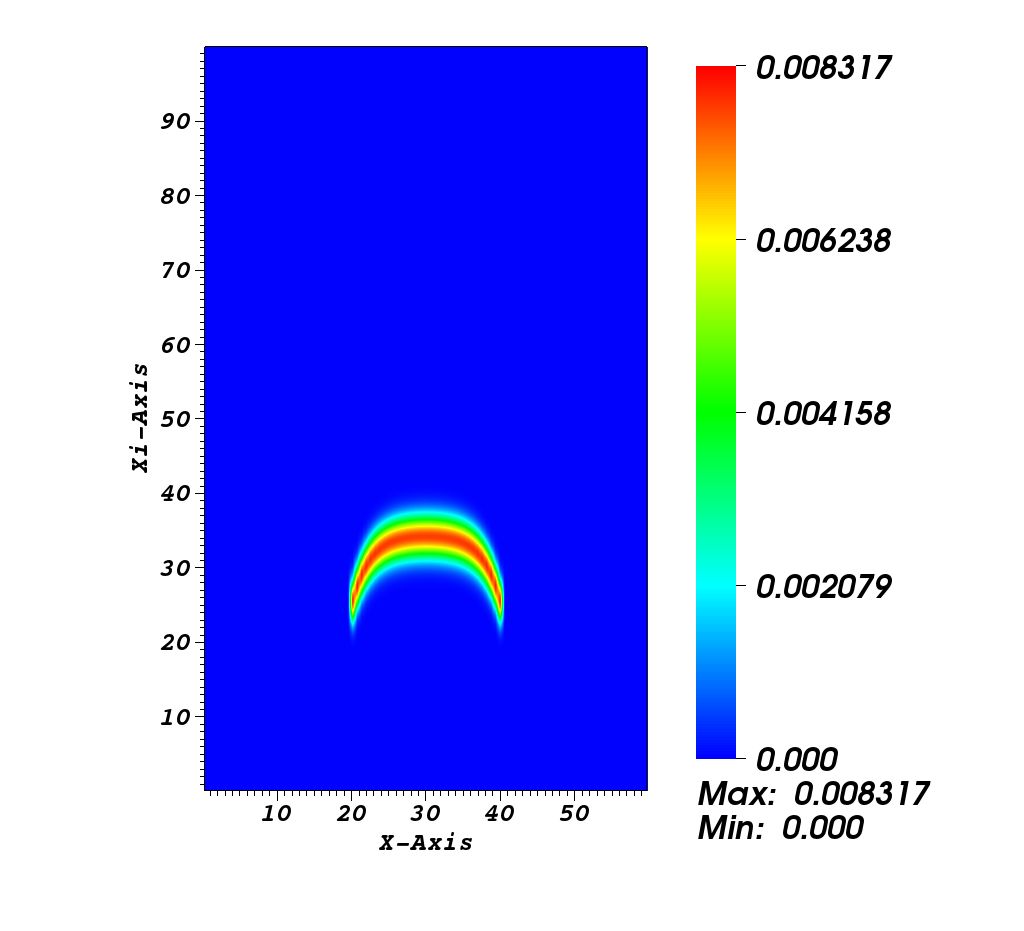}  &  \includegraphics[width=5.3cm]{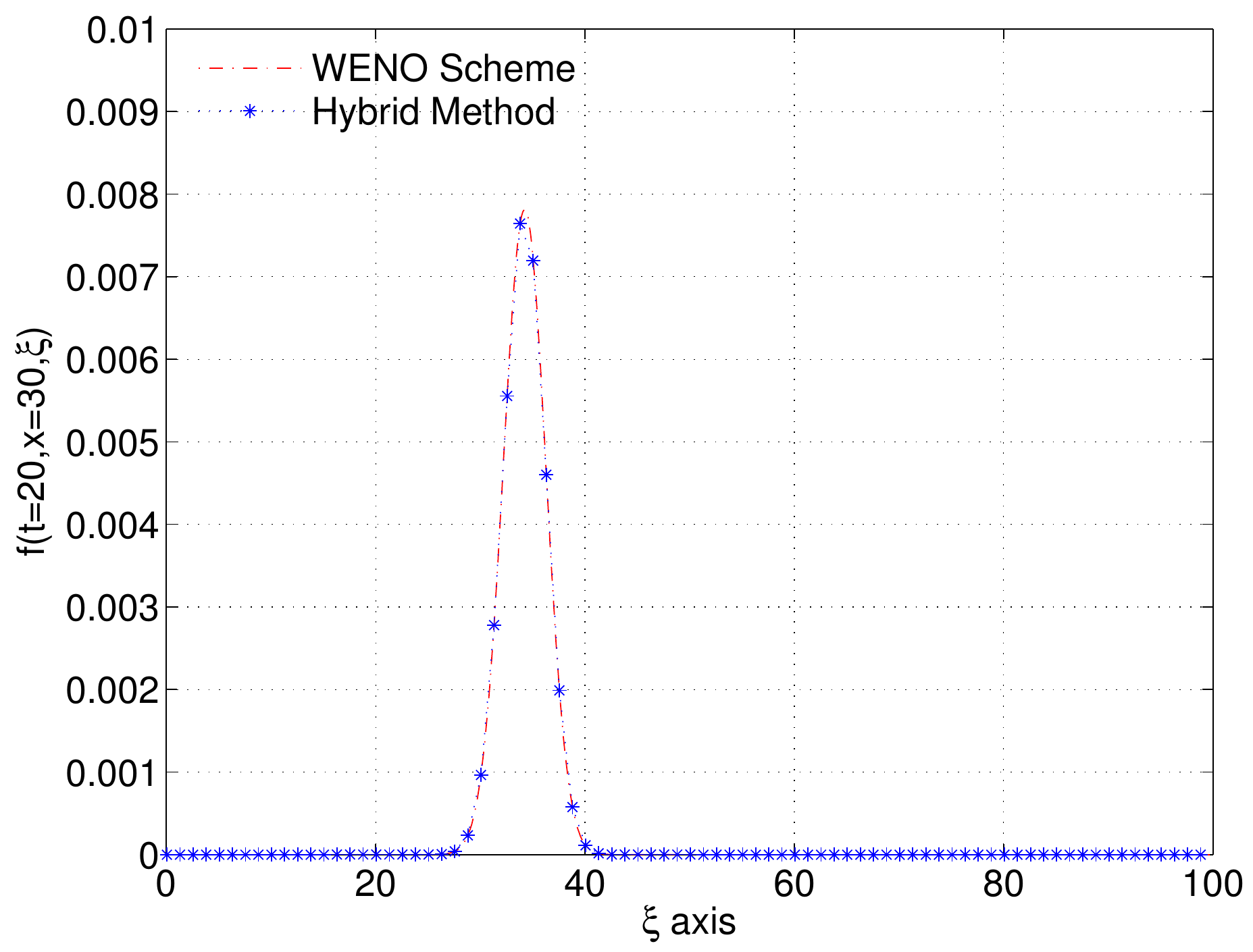}\\
  (a)  WENO Scheme                &    (b) Hybrid Method &    (c) $x=30,t=20$ \\
   \includegraphics[width=4.5cm]{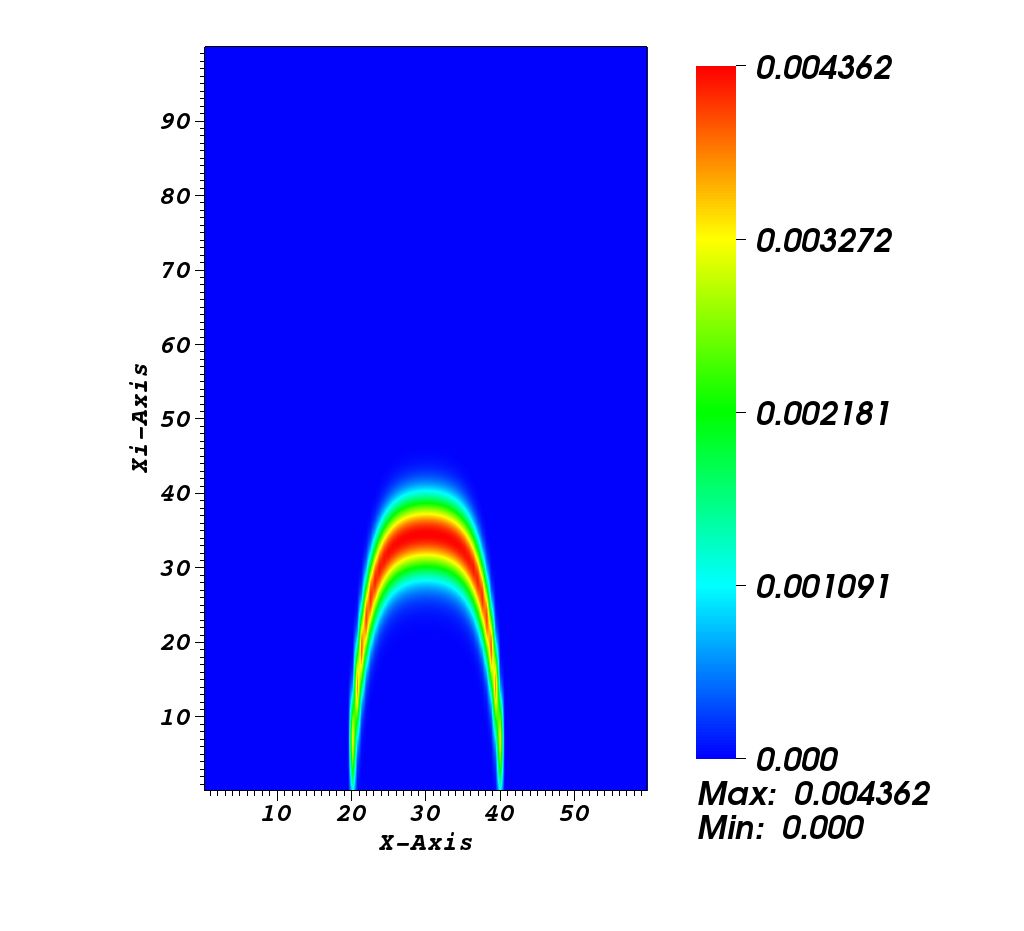}  &   \includegraphics[width=4.5cm]{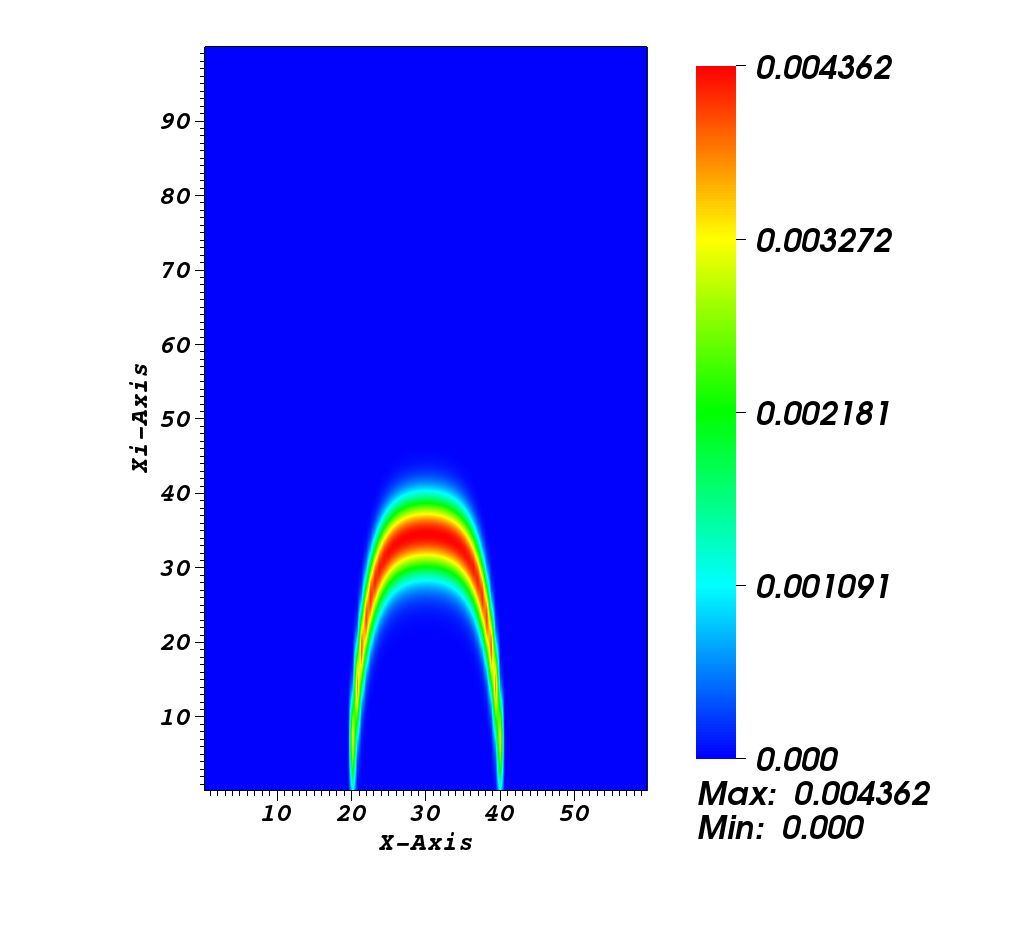}  &  \includegraphics[width=5.3cm]{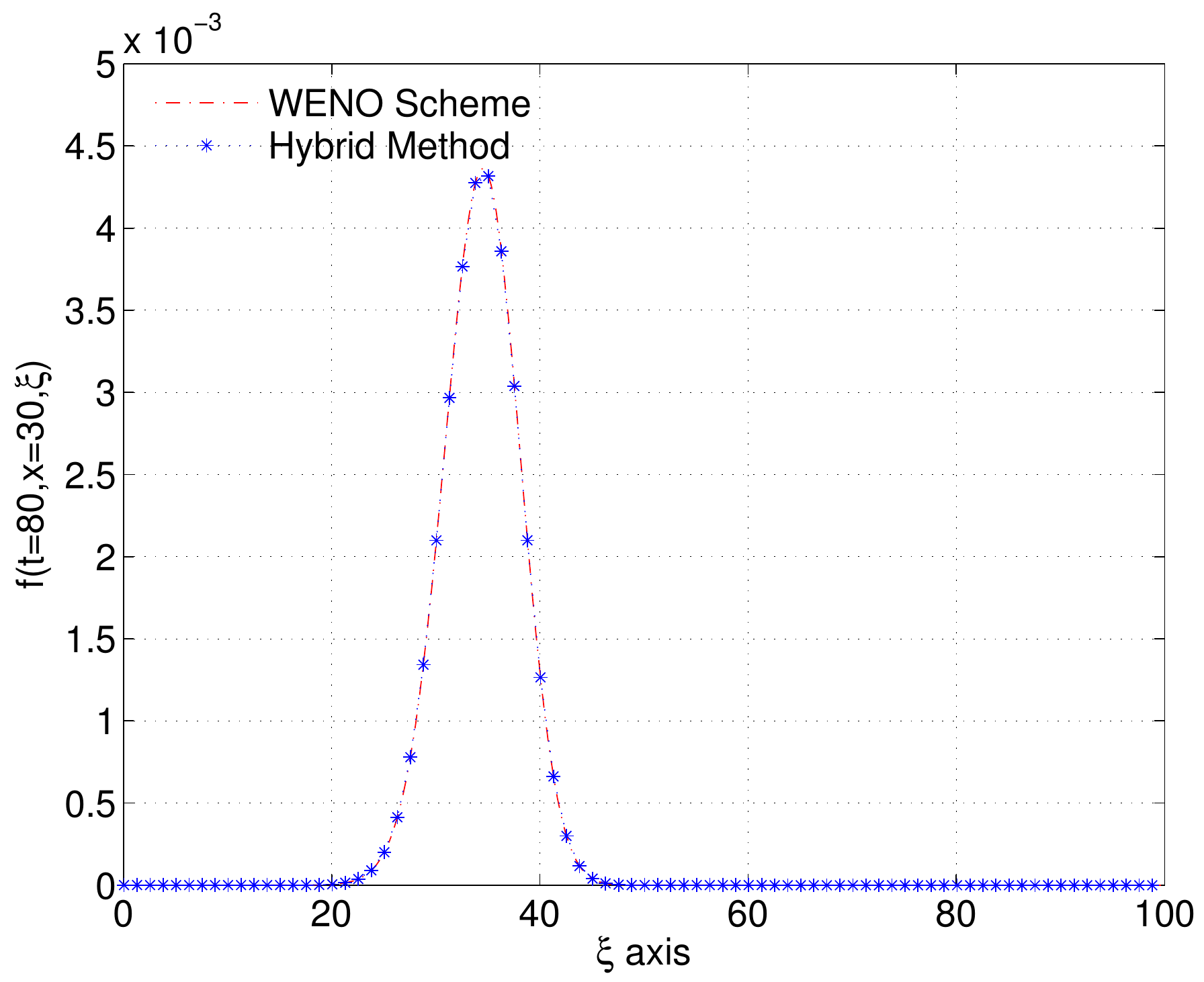}\\
  (d)  WENO Scheme                &    (e) Hybrid Method &    (f) $x=30,t=80$ \\
    \includegraphics[width=4.5cm]{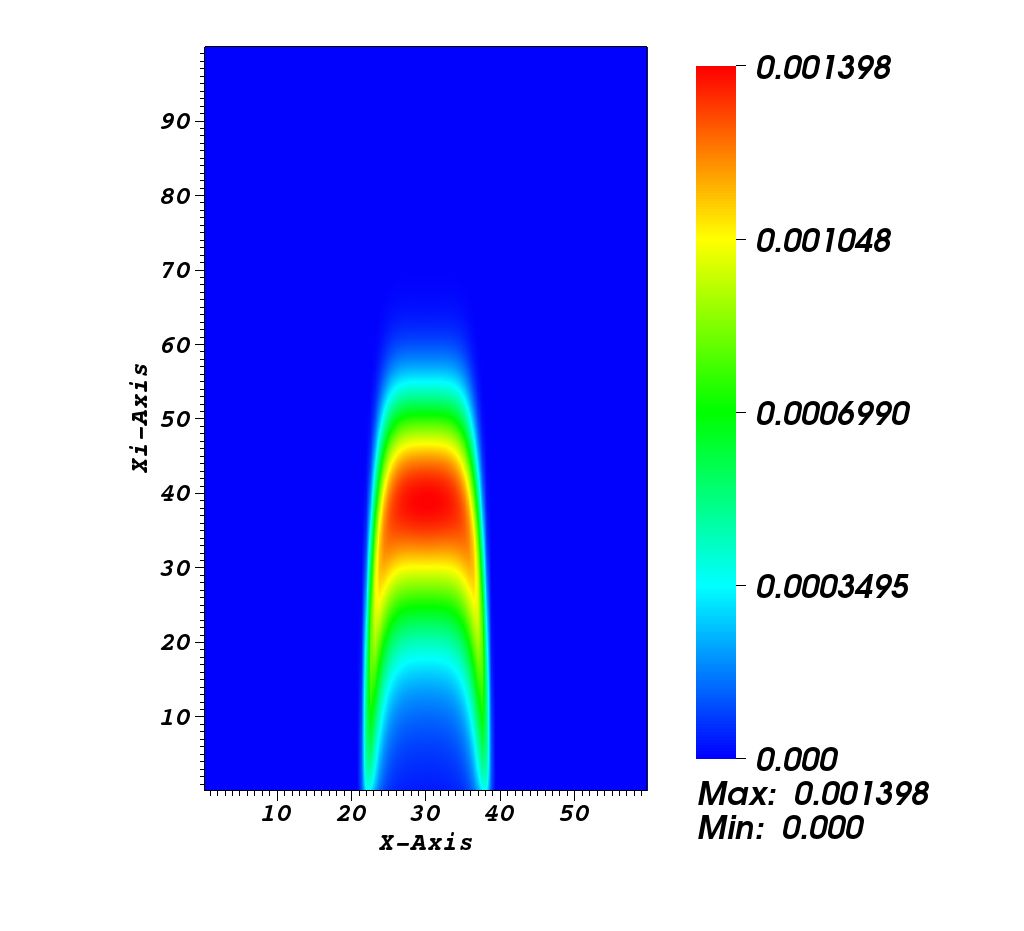}  &   \includegraphics[width=4.45cm]{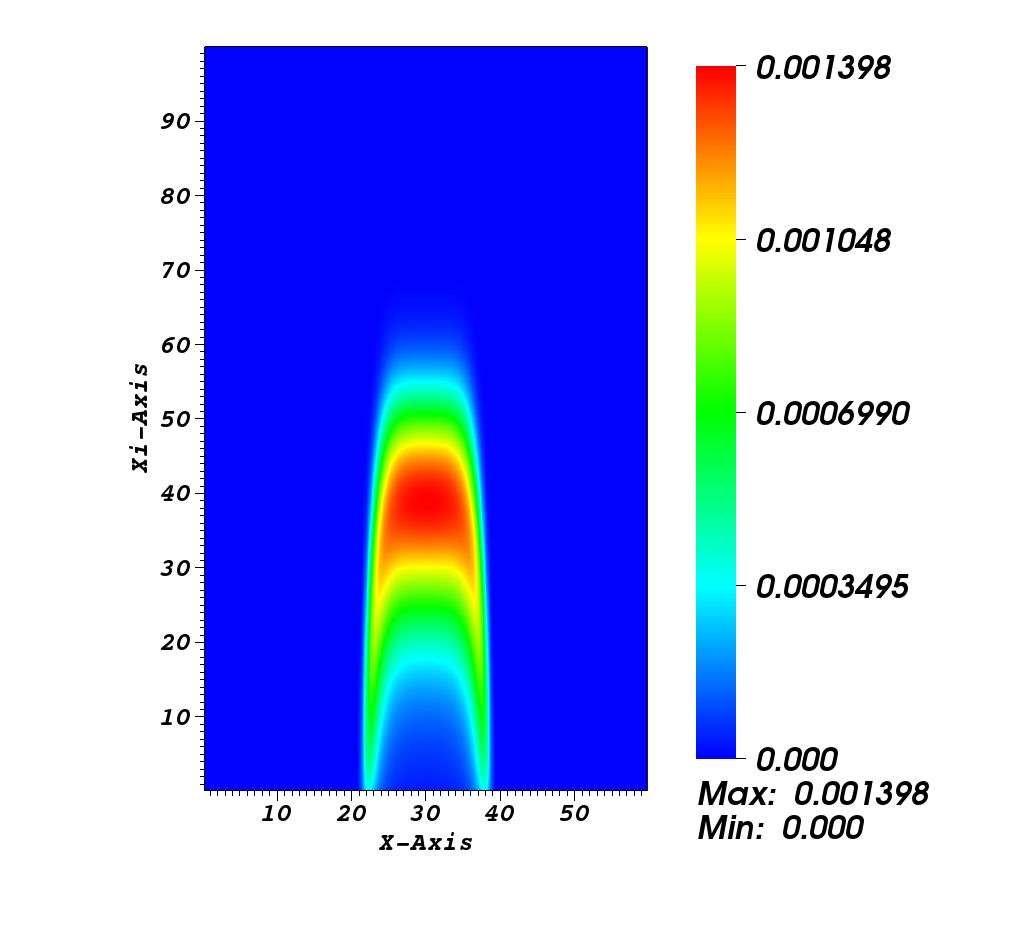}  &  \includegraphics[width=5.3cm]{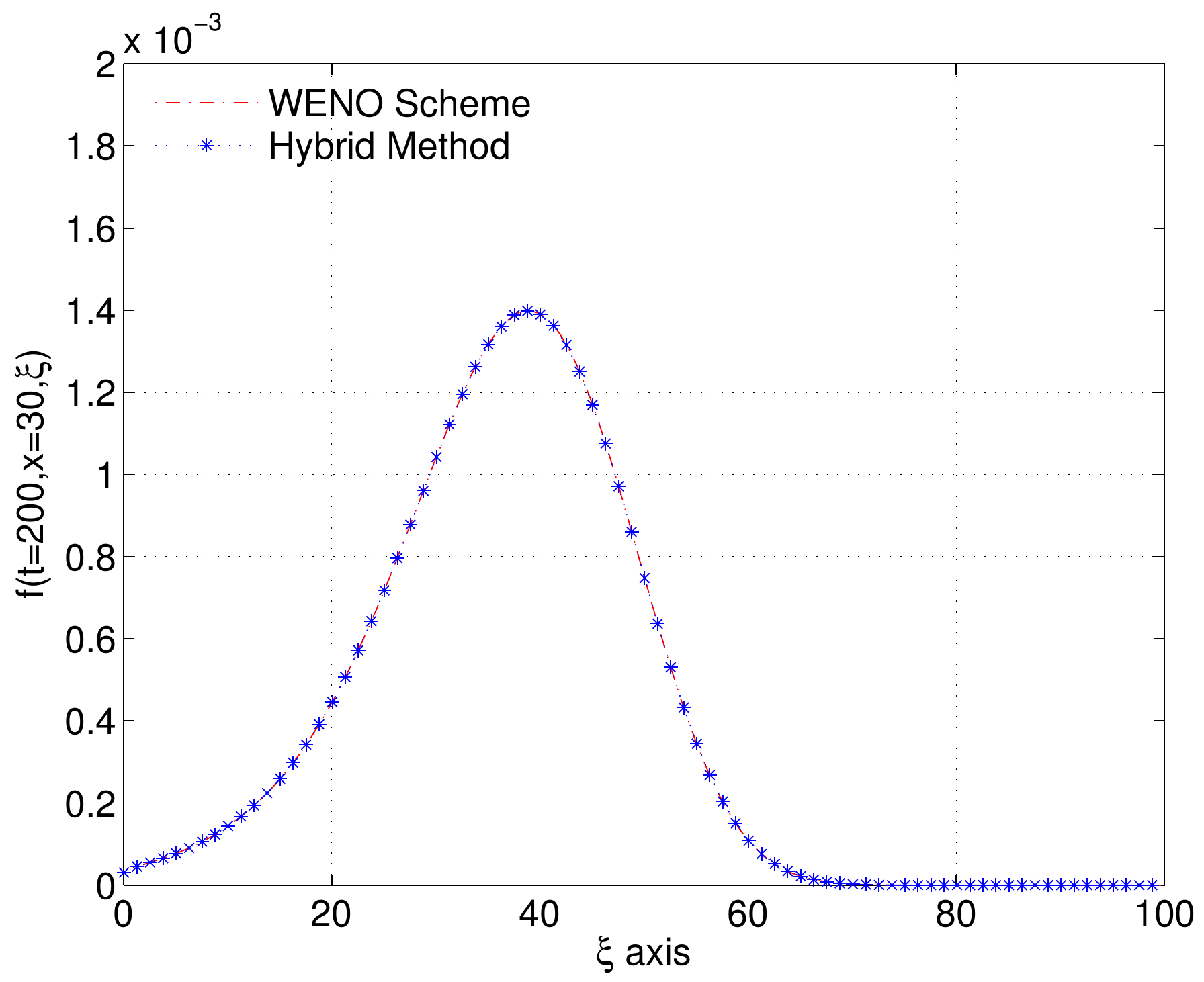}\\
  (g)  WENO Scheme                &    (h) Hybrid Method &    (i) $x=30,t=200$ \\
   \includegraphics[width=4.5cm]{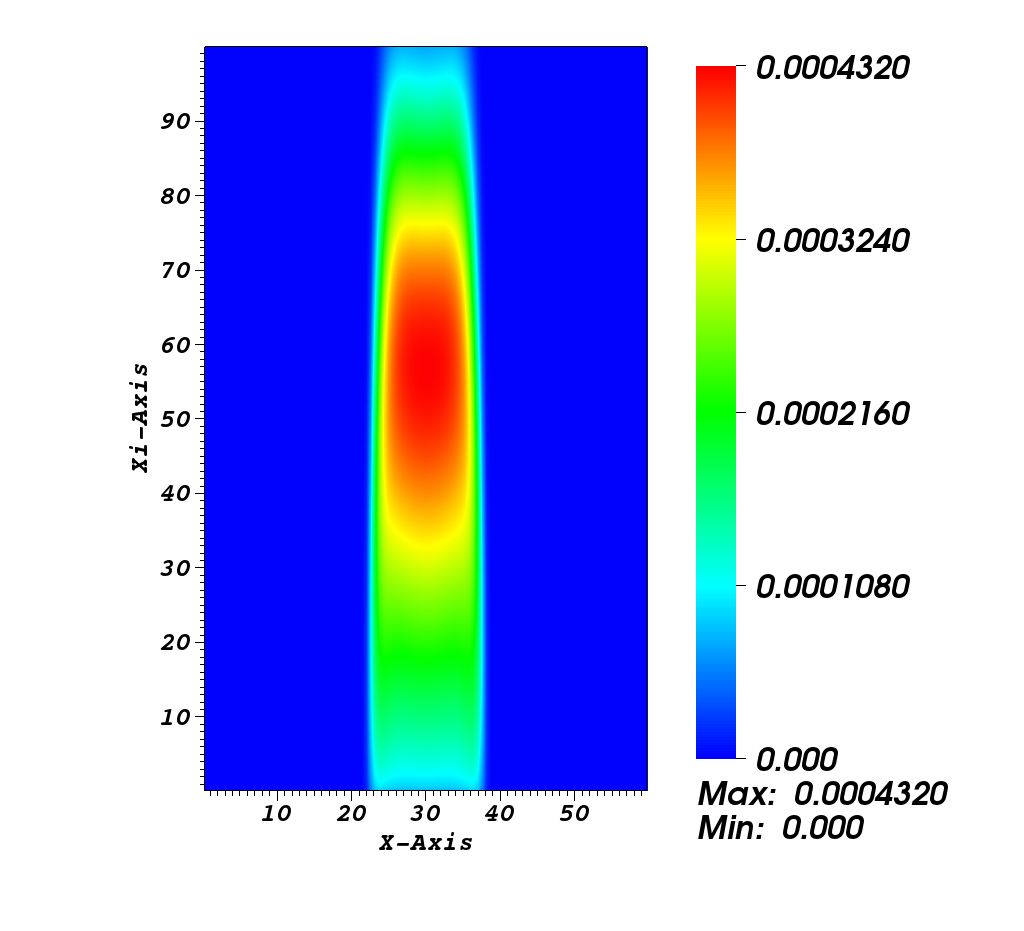}  &   \includegraphics[width=4.5cm]{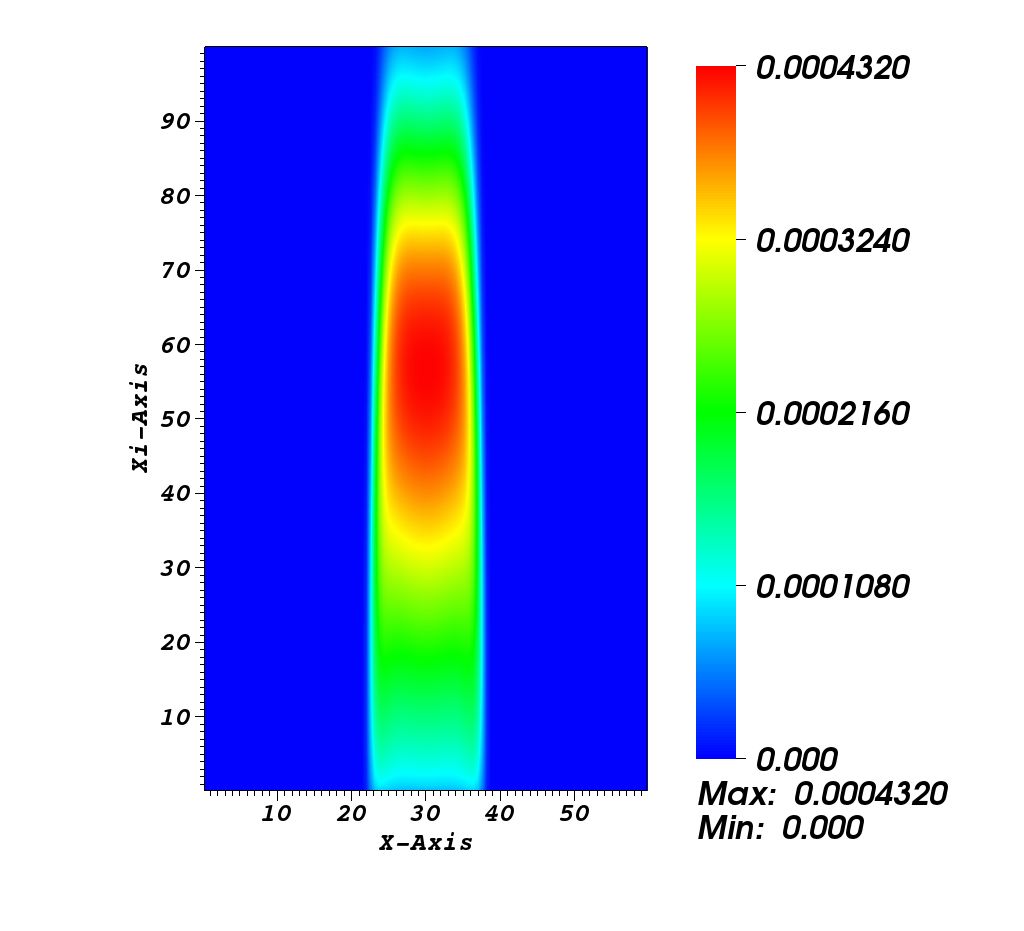}  &  \includegraphics[width=5.3cm]{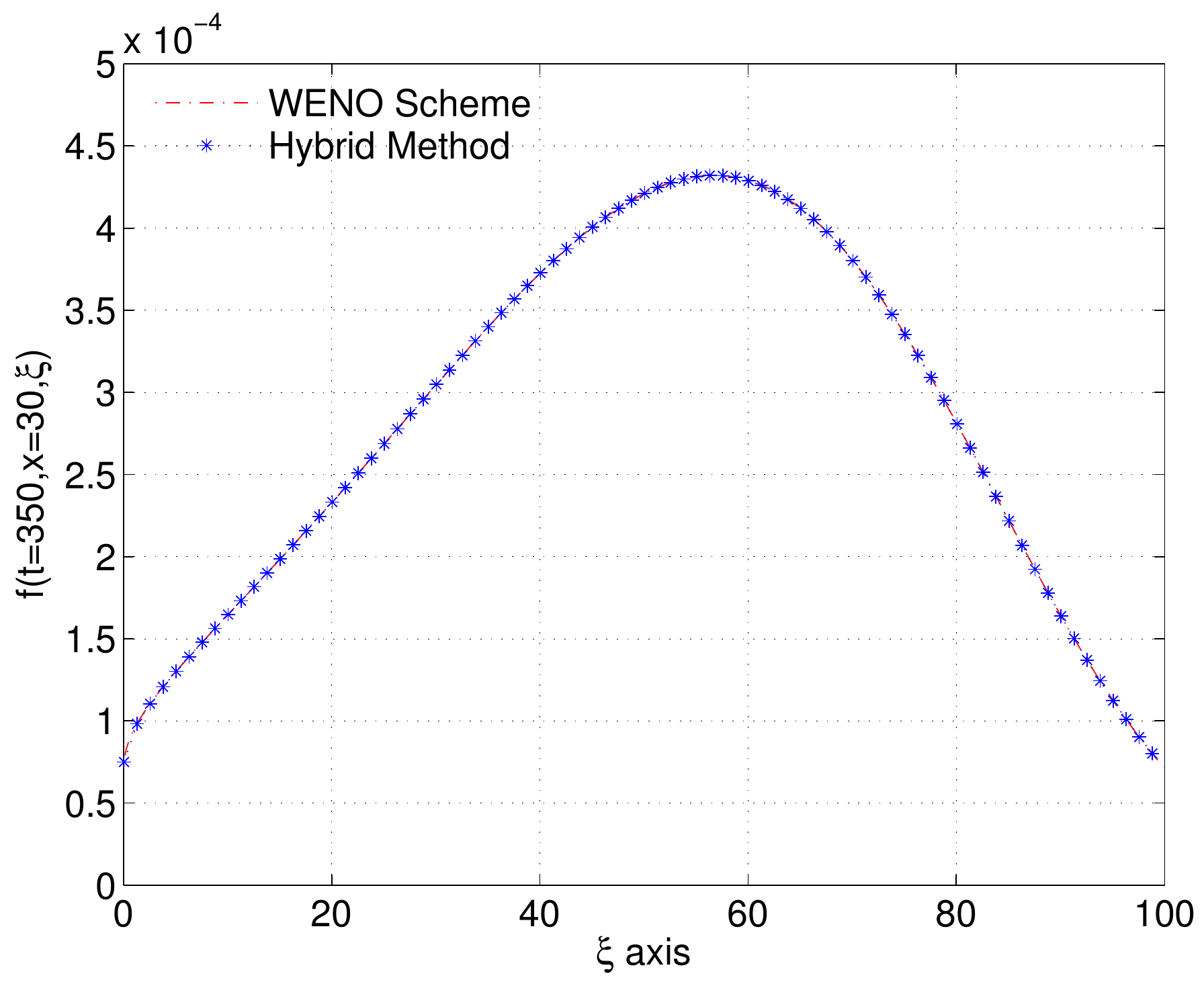}\\
  (j)  WENO Scheme                &    (k) Hybrid Method &    (l) $x=30,t=350$ \\
 \end{tabular}
\caption{\label{fig:LS2D_reg}Polymerization/depolymerization test in non-homogeneous space: {\it Plot the size distribution function of the  equation~\eqref{eq:lif-sly1d1d} with the regular initial data~\eqref{eq:lif-sly1d1d_init_concentration},~\eqref{eq:lif-sly1d1d_init_reg} at different time. Mesh size is $n_x\times n_{\xi}=100\times800$, $\Delta t=0.1$, $CFL\approx0.13$.}}
 \end{center}
\end{figure}

 Firstly, in the case with the regular data~\eqref{eq:lif-sly1d1d_init_reg},  there are almost no difference from the size-density distribution between the two numerical methods (see Figure~\ref{fig:LS2D_reg}). The initial distribution is in a plaque form at $t=0$, then it generates a moon shape caused by the diffusion of the concentration of monomers. In evolution in time, we see that the solution is very smooth in size direction.

\begin{figure}
\begin{center}
 \begin{tabular}{cc}
  \includegraphics[width=7cm]{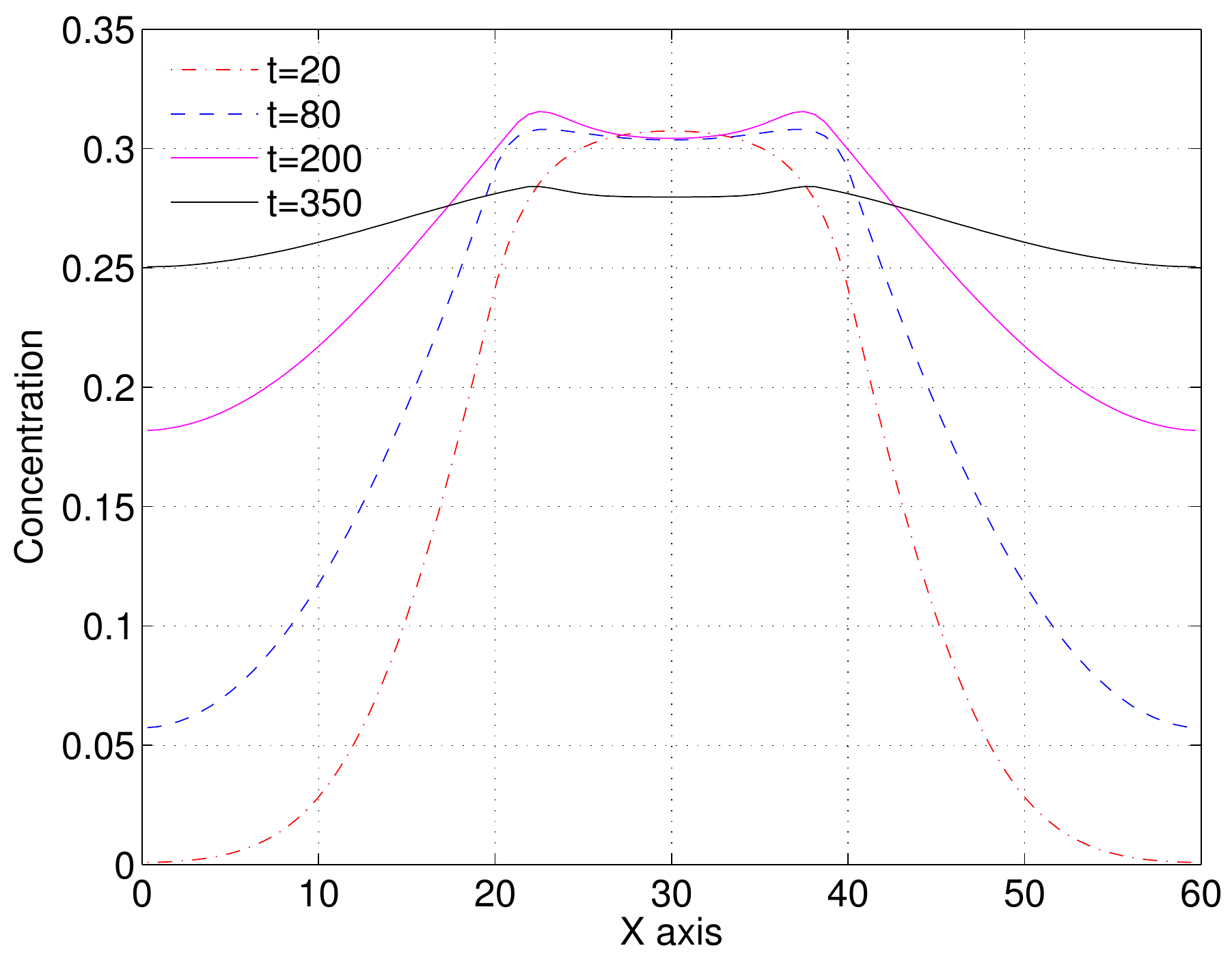}  &   \includegraphics[width=7cm]{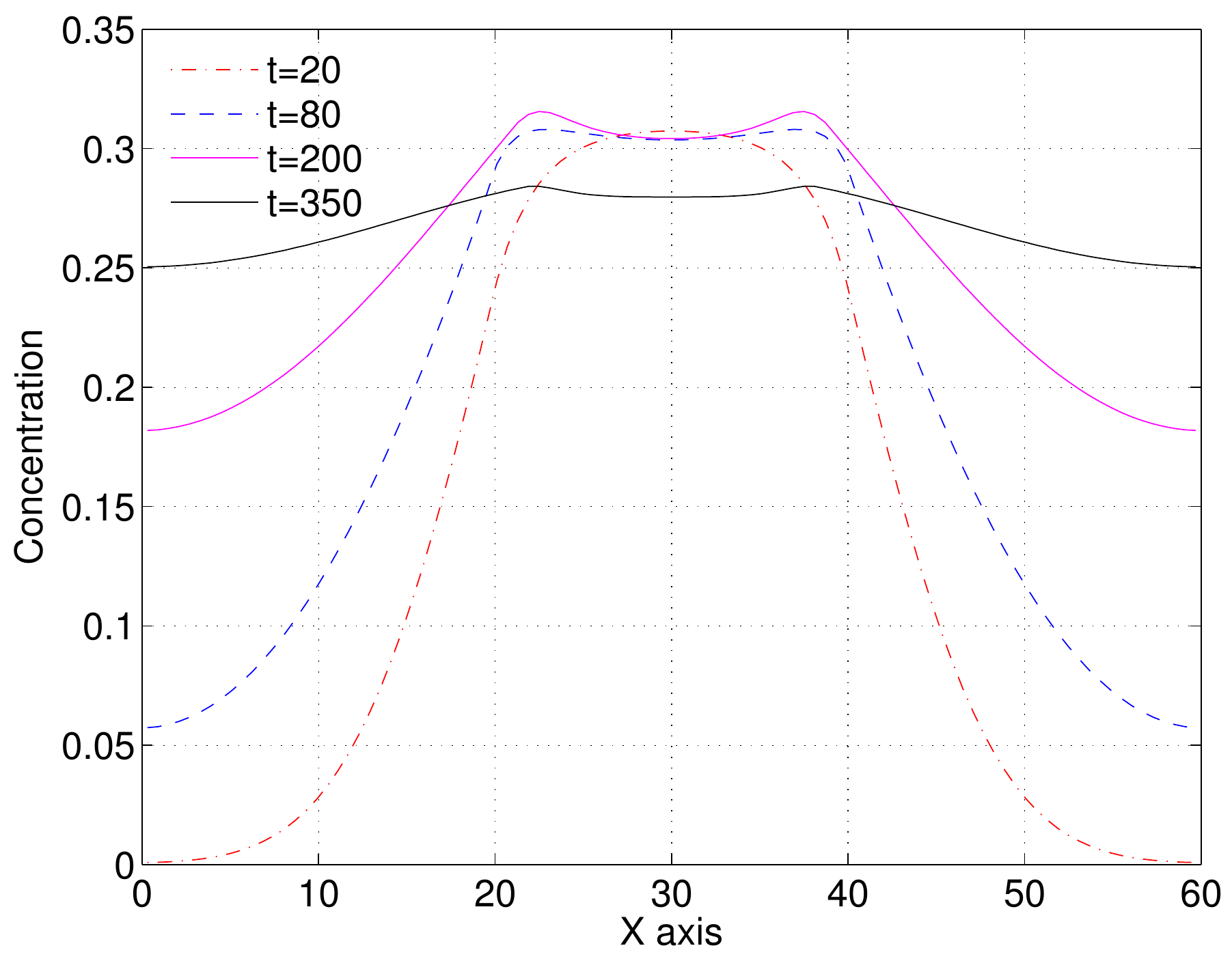} \\
  (a)  WENO scheme                &    (b) Hybrid Method \\
  \includegraphics[width=7cm]{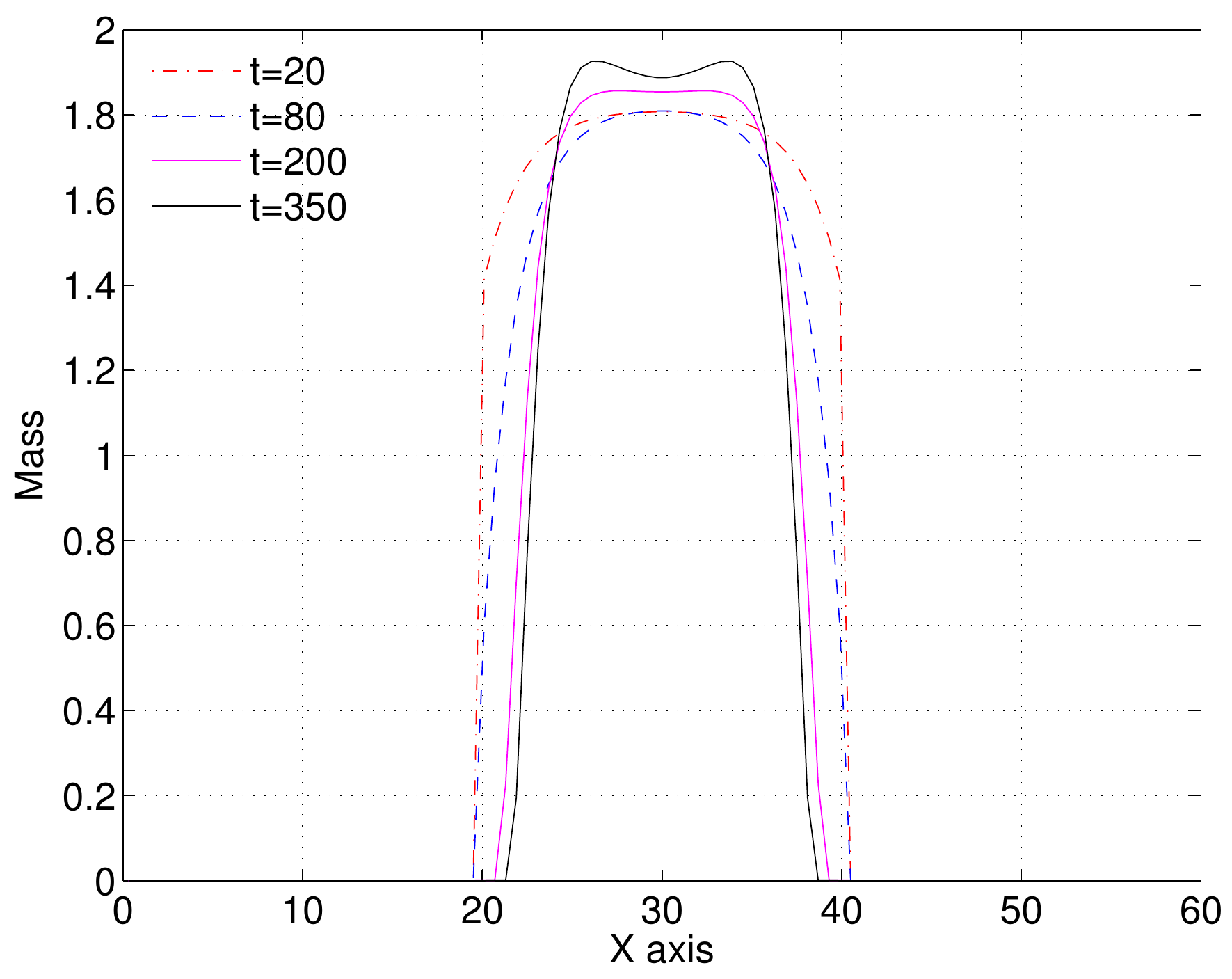}  &   \includegraphics[width=7cm]{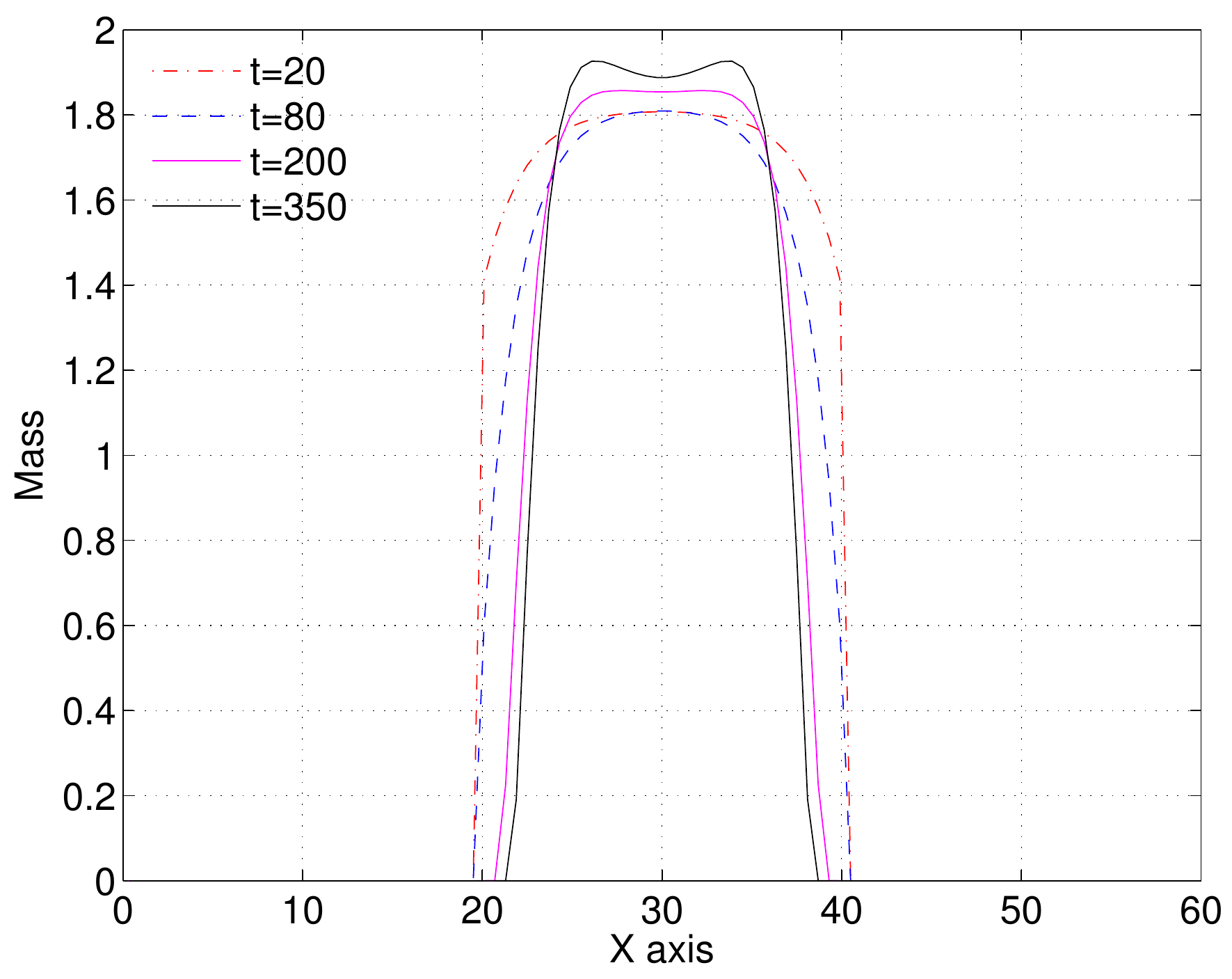} \\
  (c)  WENO scheme                &    (d) Hybrid Method
 \end{tabular}
\caption{\label{fig:LS2D_concentration_mass}Polymerization/depolymerization test in non-homogeneous space: {\it Evolution of the monomers concentration $c(t,x)$ (first row) and the mass of the marco-particles $\int_0^{\infty}\xi f(t,x,\xi)d\xi$ (second row) corresponding to the  equation~\eqref{eq:lif-sly1d1d} with the irregular initial data~\eqref{eq:lif-sly1d1d_init_concentration},~\eqref{eq:lif-sly1d1d_init_irreg}. Mesh size is $n_x\times n_{\xi}=100\times800$, $\Delta t=0.1$, $CFL\approx0.13$.}}
 \end{center}
\end{figure}

\begin{figure}
\begin{center}
 \begin{tabular}{cc}
  \includegraphics[width=6cm]{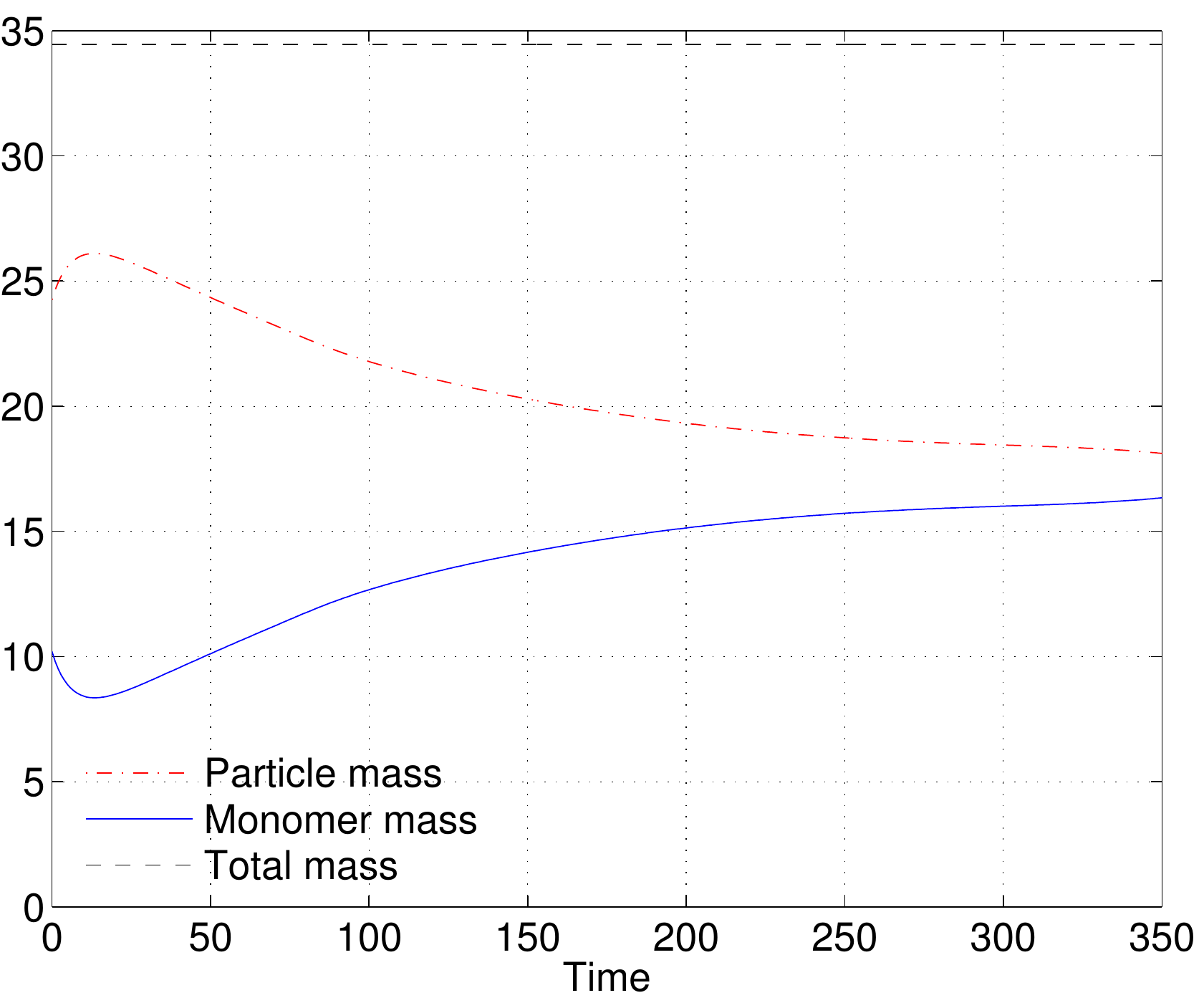}  &   \includegraphics[width=6cm]{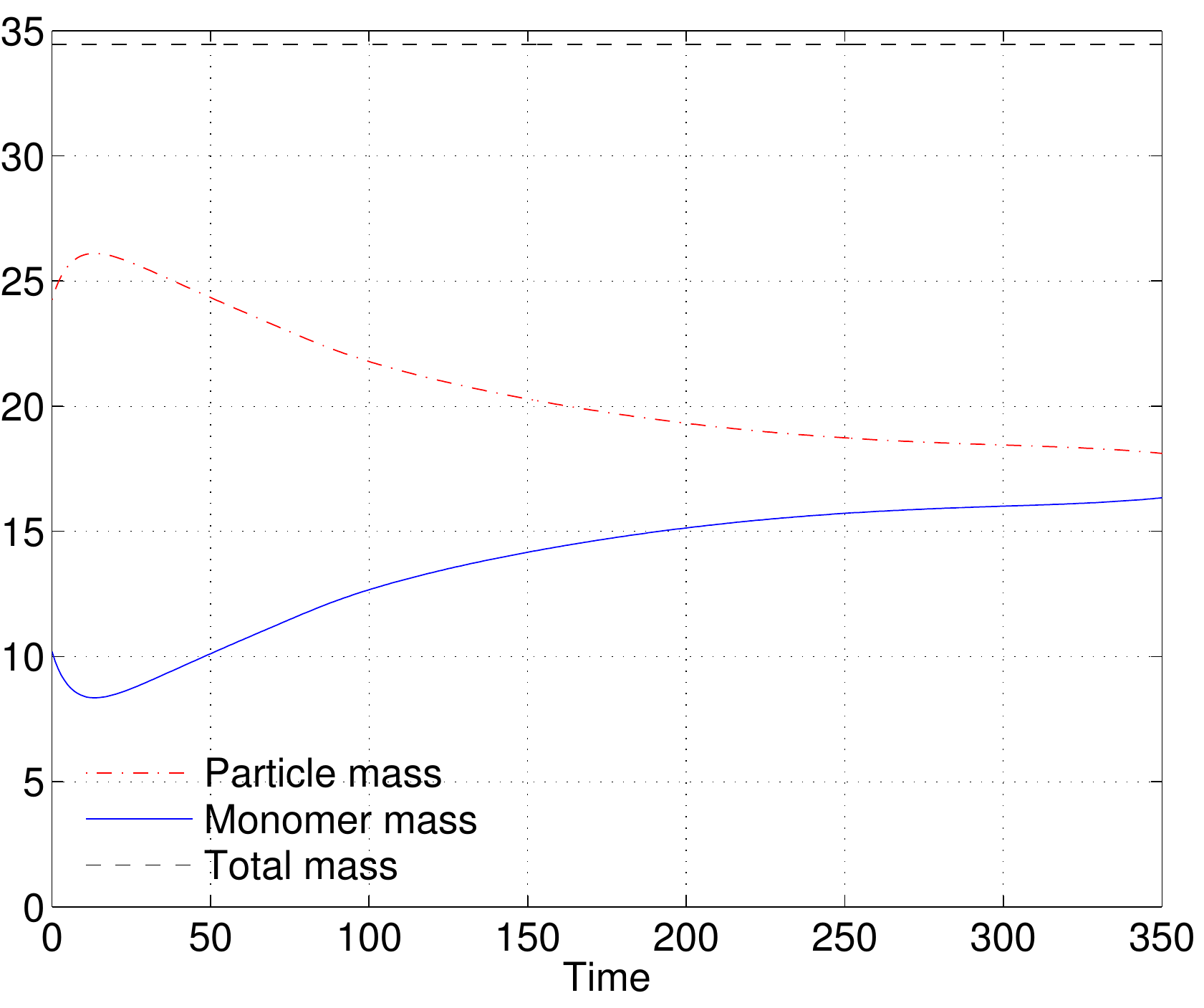}\\
  (a)  WENO method                &    (b) Hybrid method
 \end{tabular}
\caption{\label{fig:growth_nonhomo_conserv}Polymerization/depolymerization test in homogeneous space: {\it Mass conservation property corresponding to~\eqref{eq:mass_conservation}. Mesh size is $n_x\times n_{\xi}=100\times800$, $\Delta t=0.1$, $CFL\approx0.13$.}}
 \end{center}
\end{figure}

\begin{figure}
\begin{center}
 \begin{tabular}{ccc}
  \includegraphics[width=4.5cm]{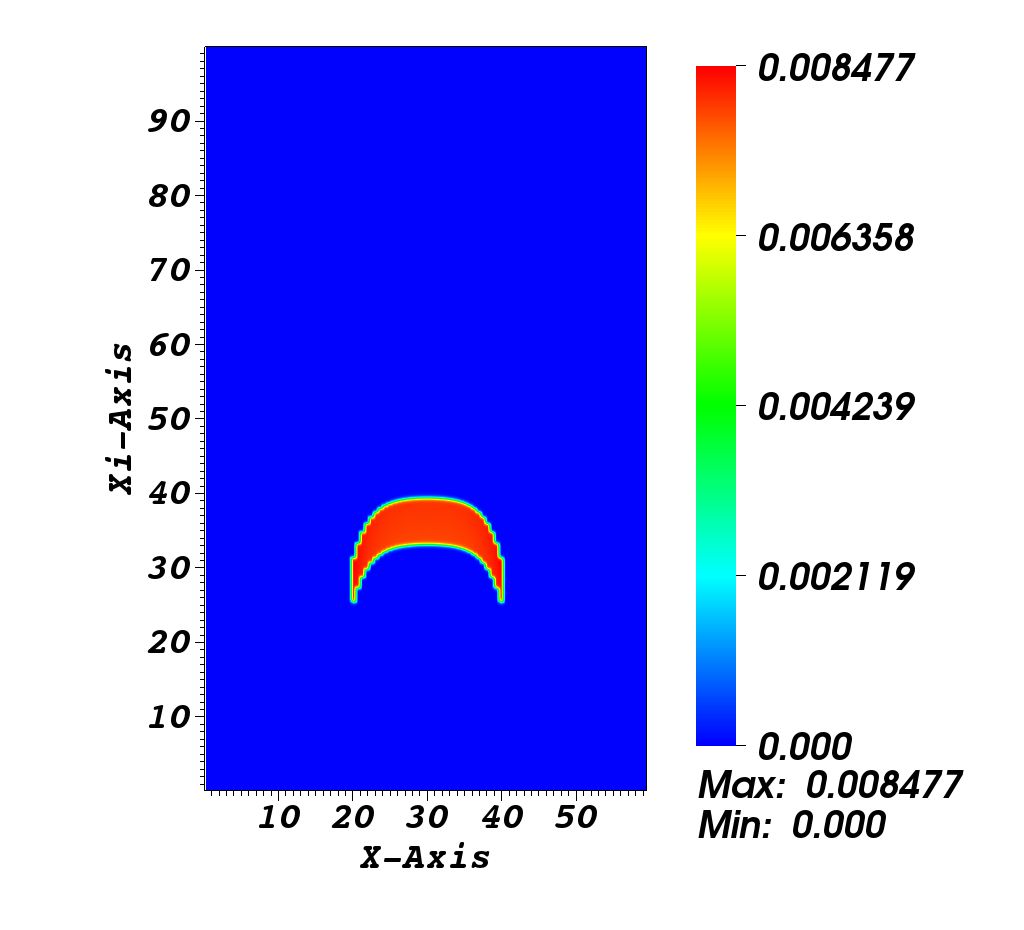}  &   \includegraphics[width=4.5cm]{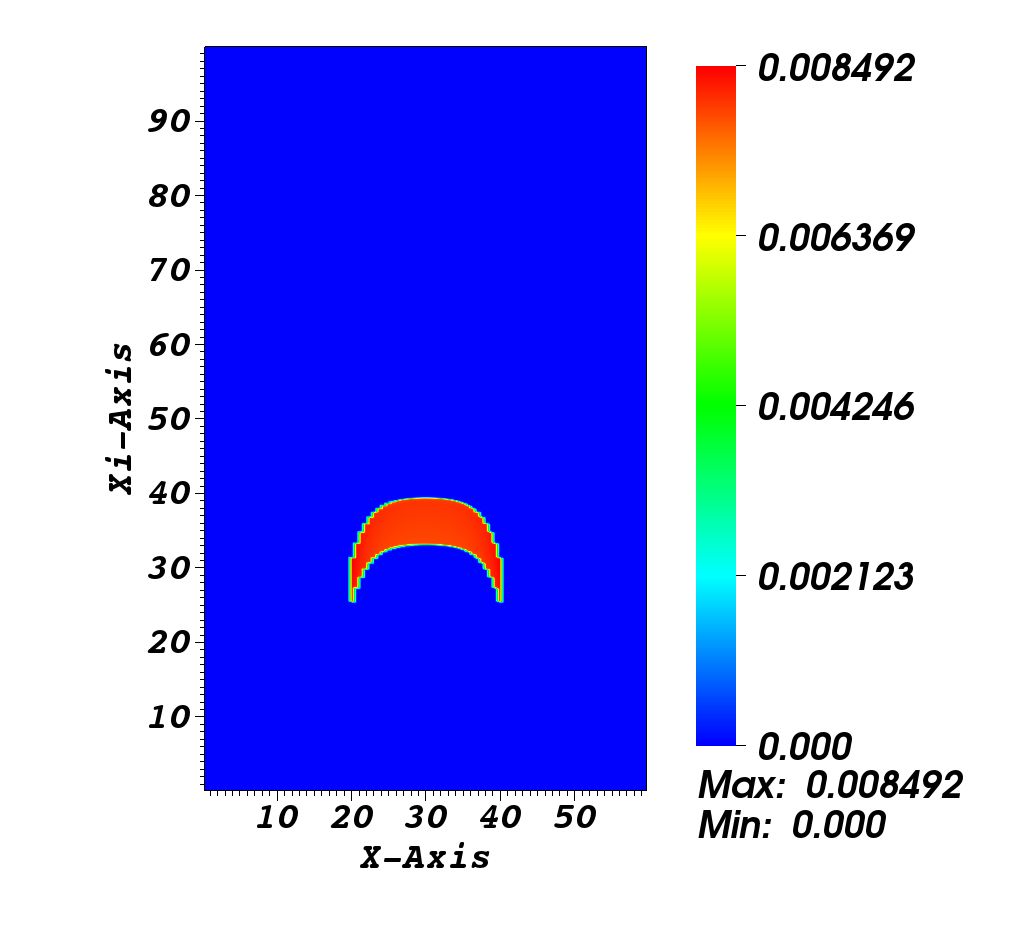}  &  \includegraphics[width=5.3cm]{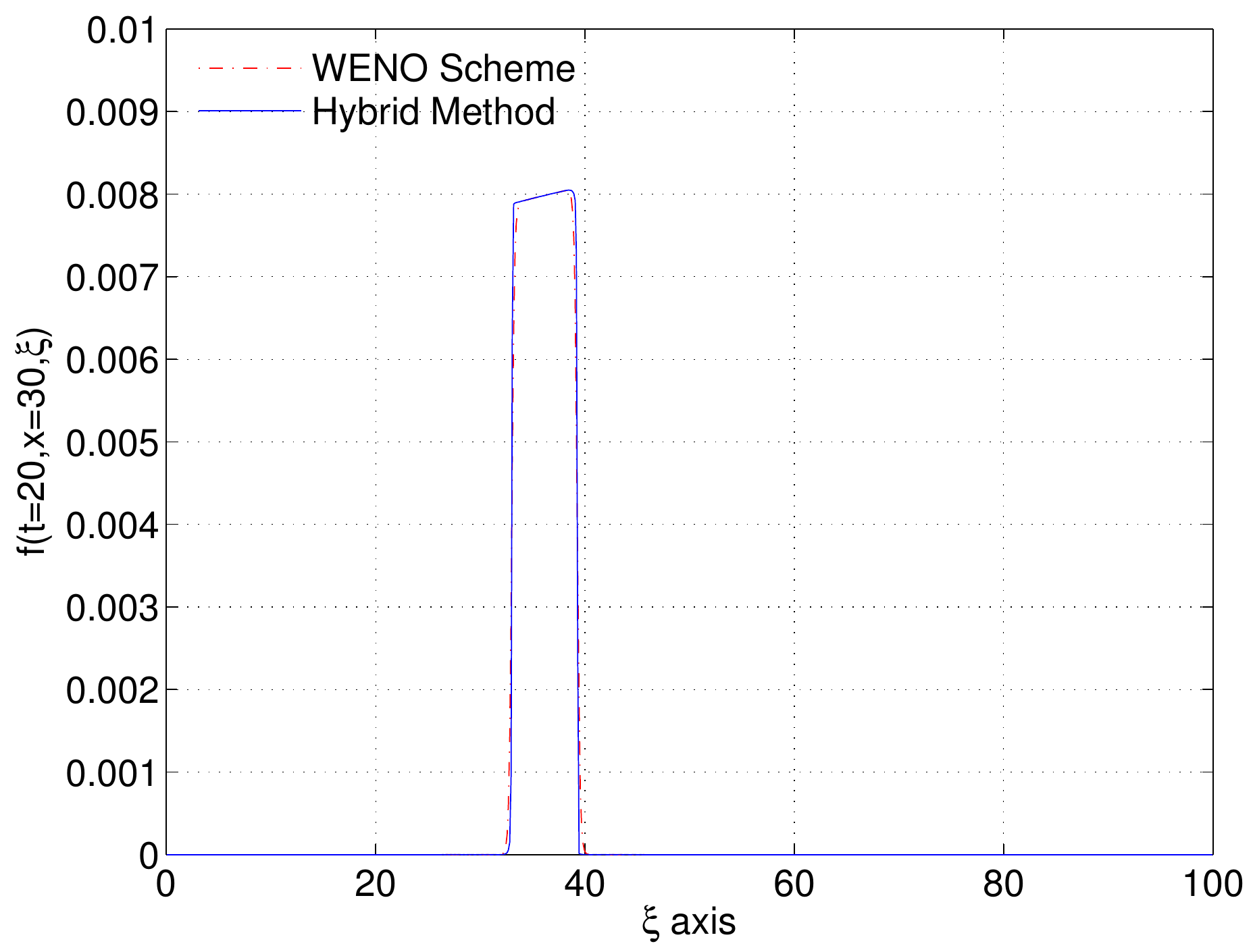}\\
  (a)  WENO Scheme                &    (b) Hybrid Method &    (c) $x=30,t=20$ \\
   \includegraphics[width=4.5cm]{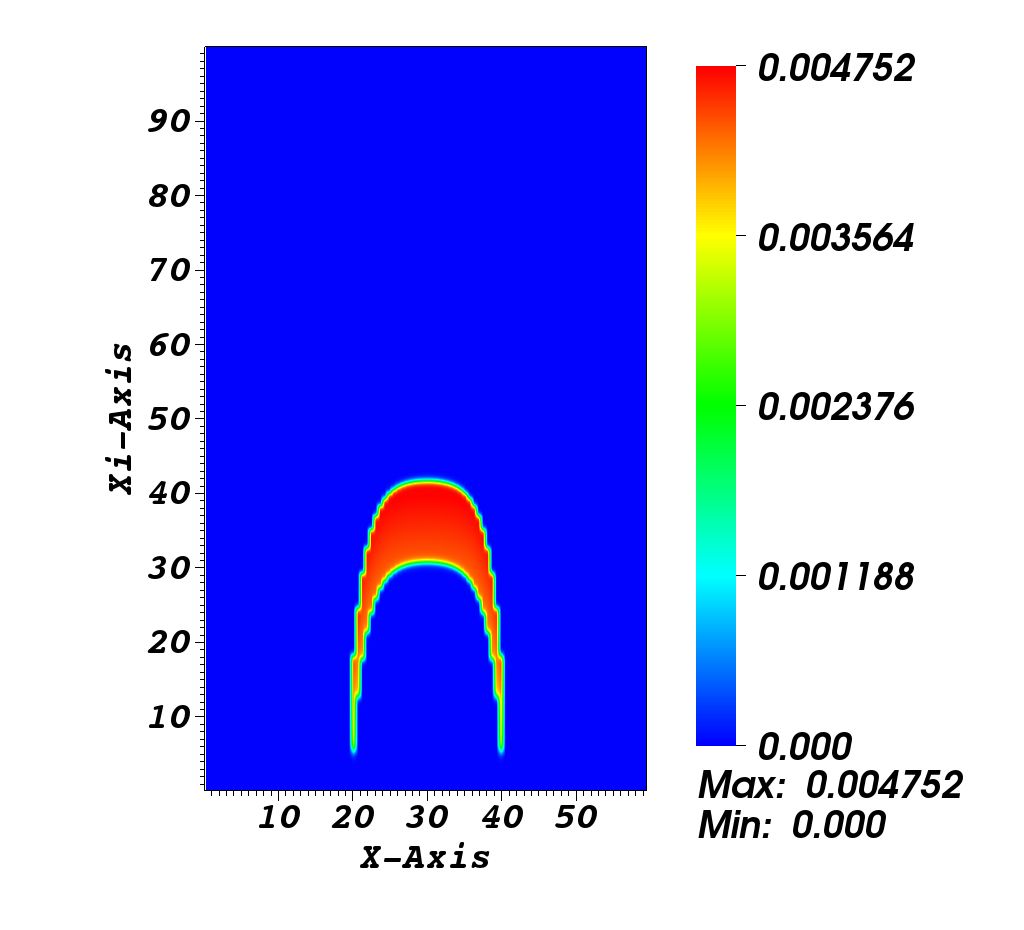}  &   \includegraphics[width=4.5cm]{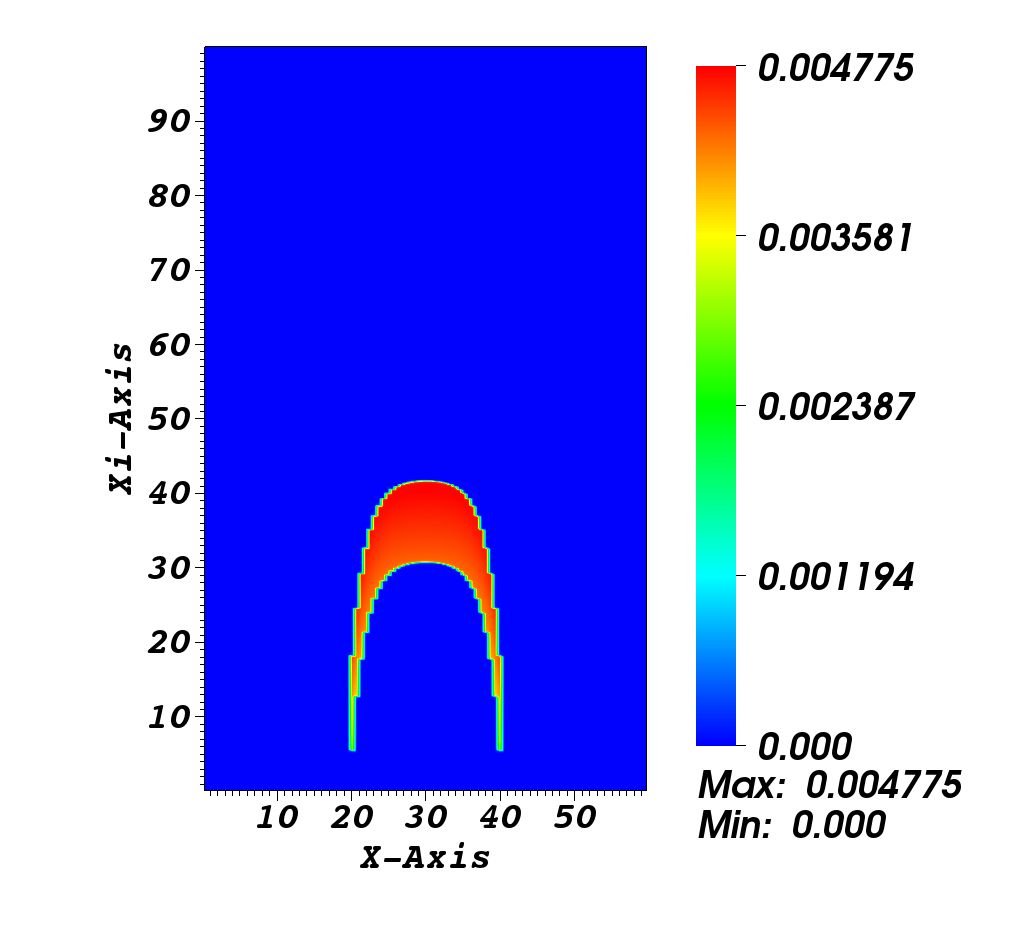}  &  \includegraphics[width=5.3cm]{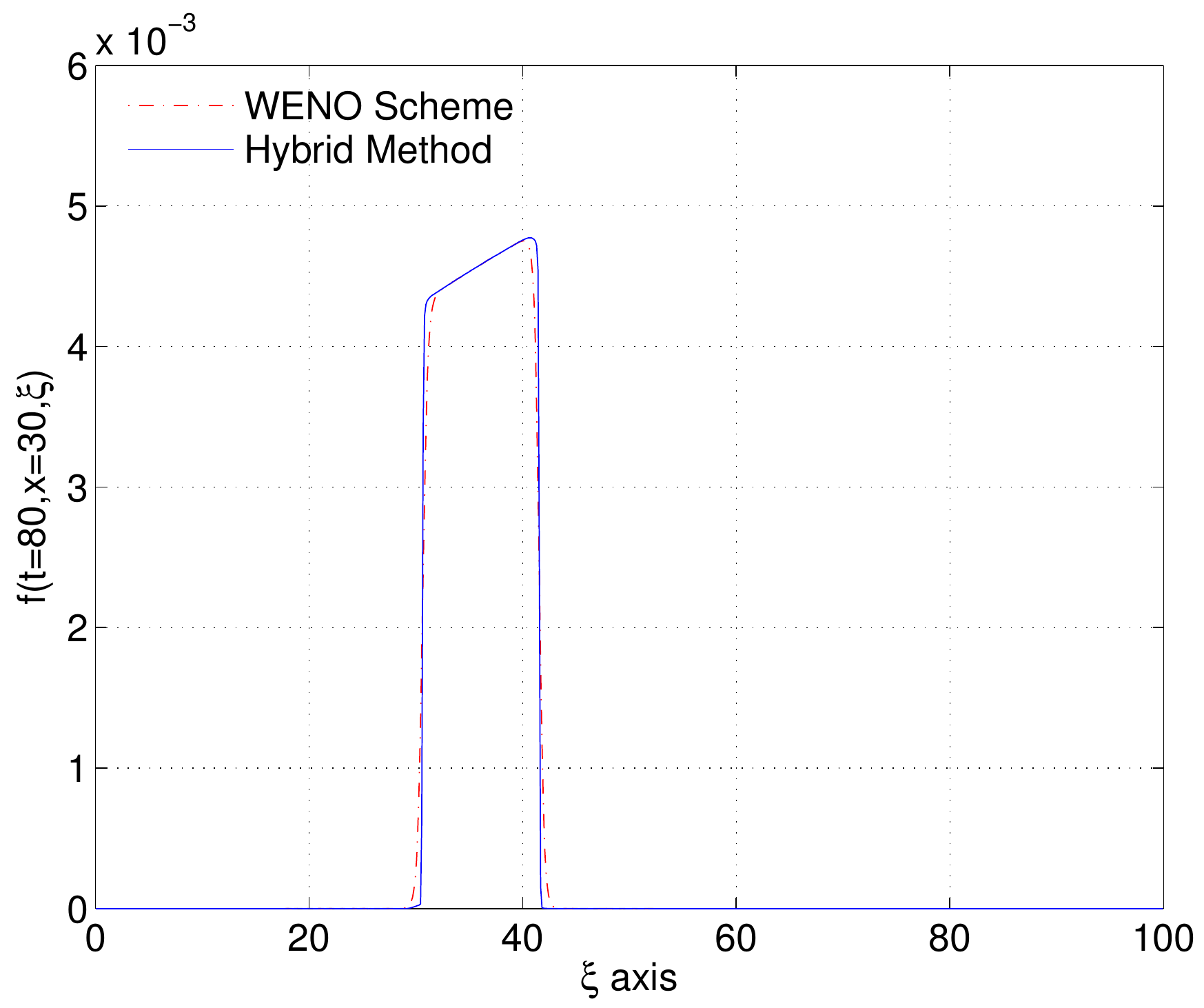}\\
  (d)  WENO Scheme                &    (e) Hybrid Method &    (f) $x=30,t=80$ \\
    \includegraphics[width=4.5cm]{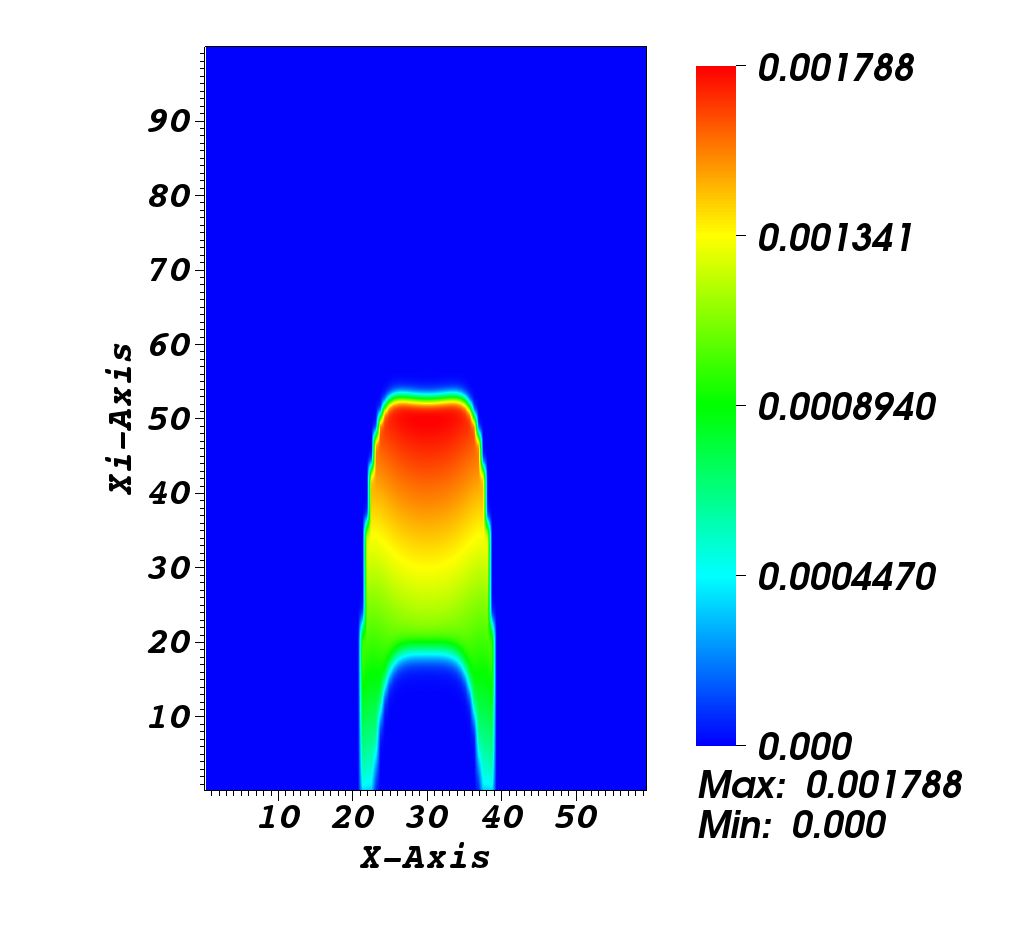}  &   \includegraphics[width=4.45cm]{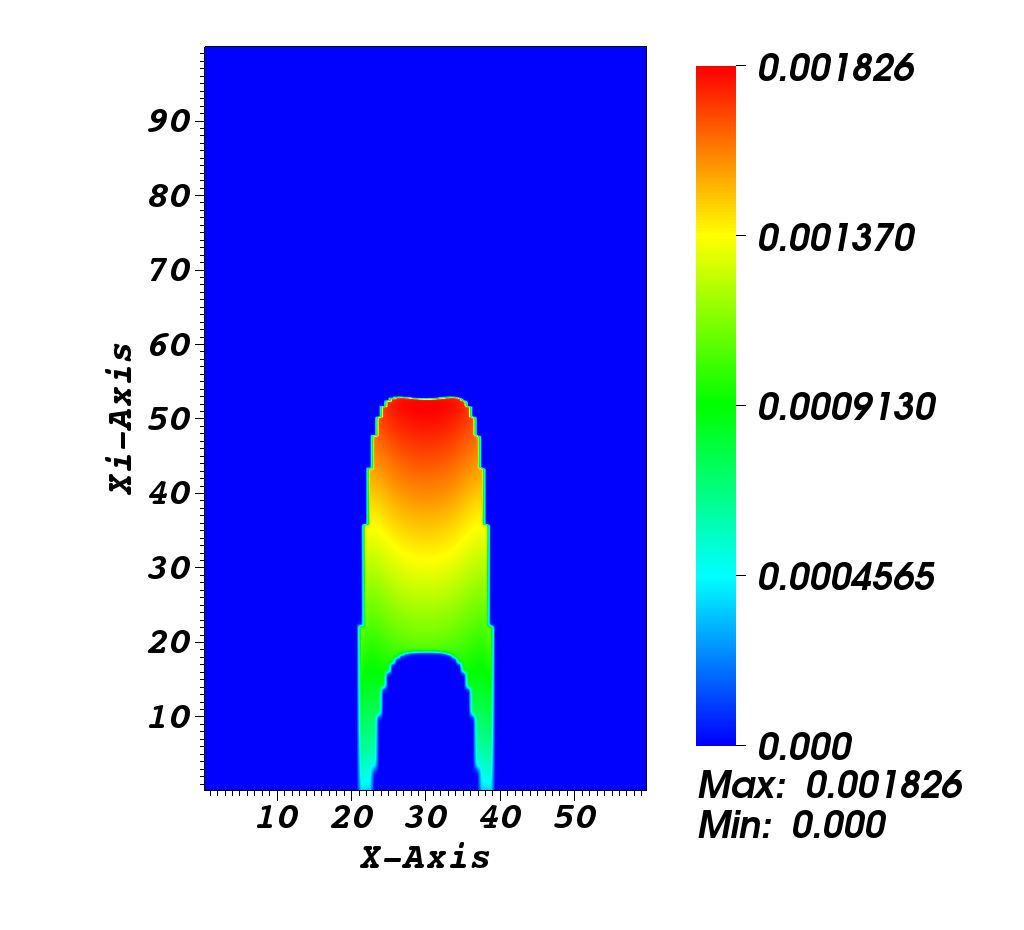}  &  \includegraphics[width=5.3cm]{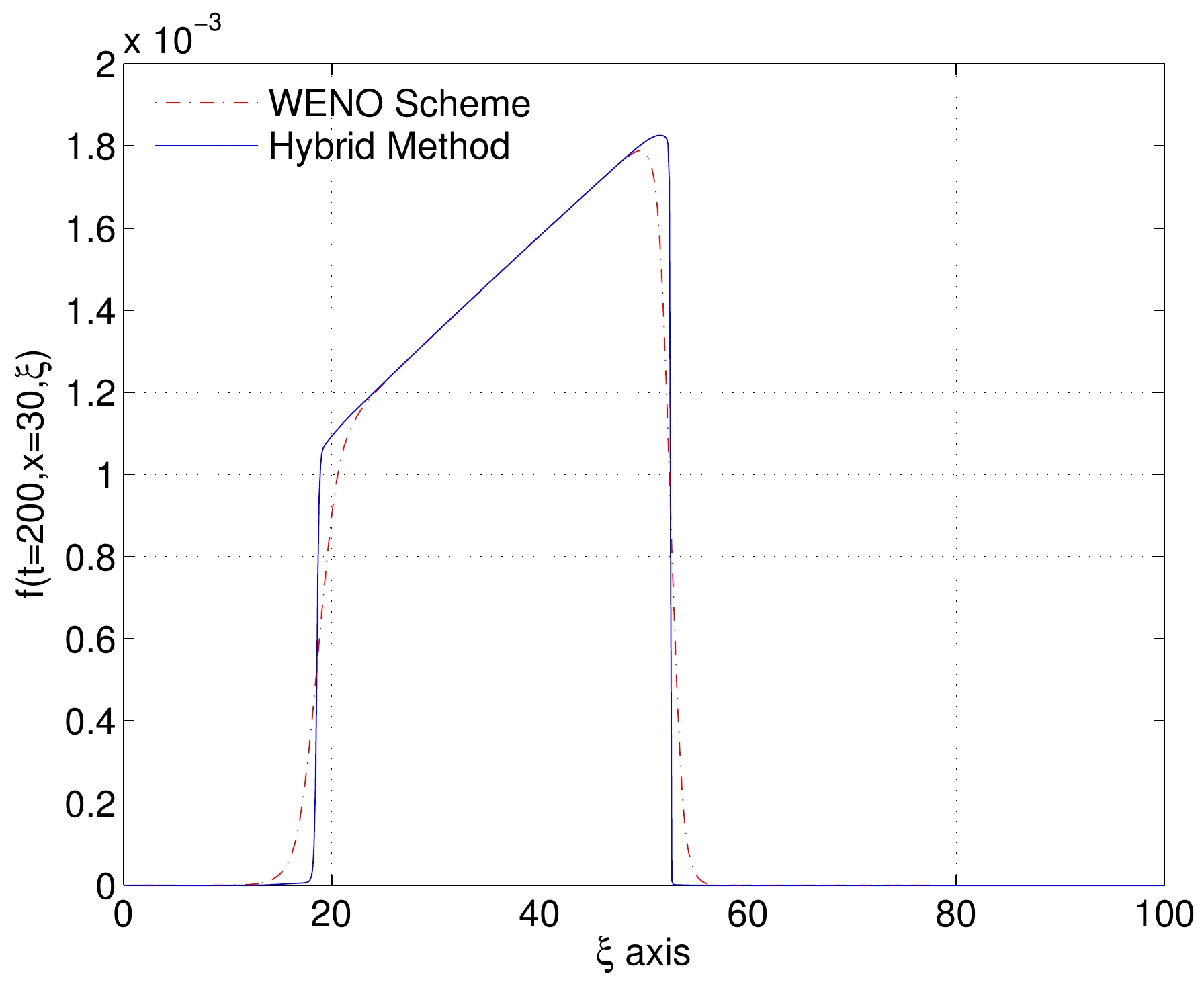}\\
  (g)  WENO Scheme                &    (h) Hybrid Method &    (i) $x=30,t=200$ \\
   \includegraphics[width=4.5cm]{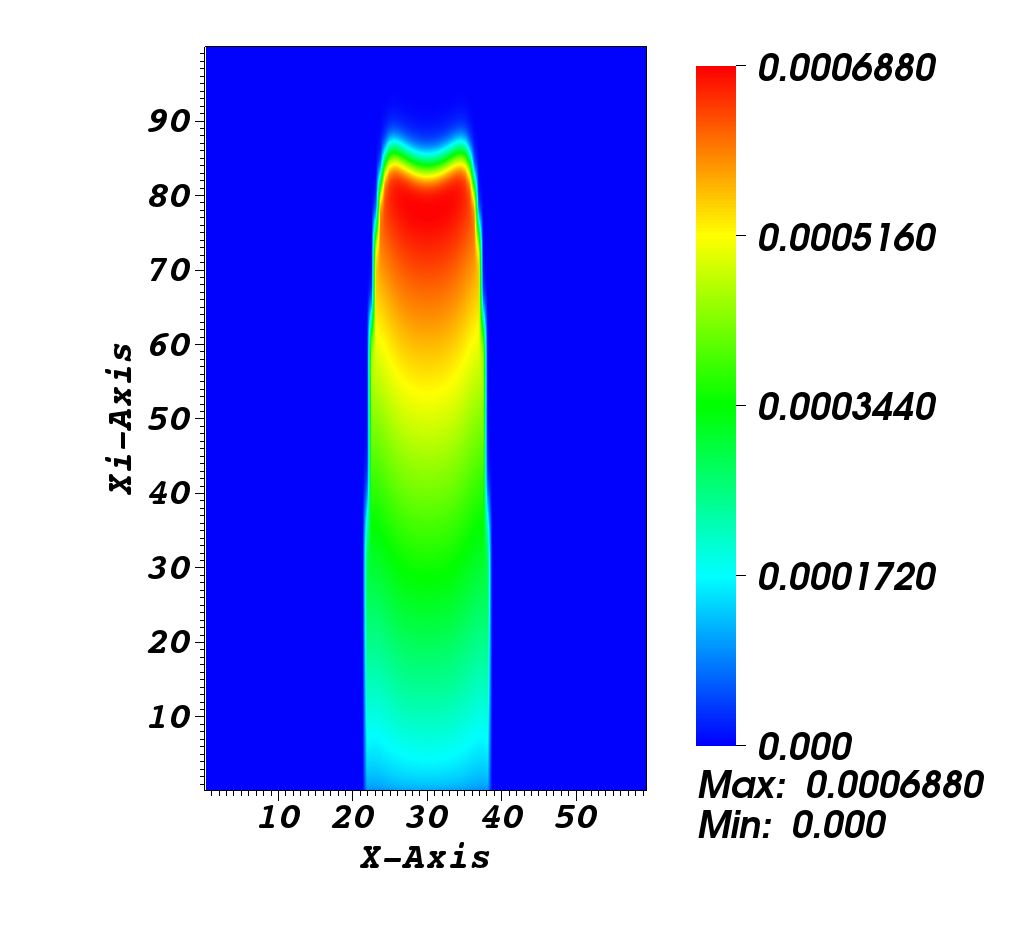}  &   \includegraphics[width=4.5cm]{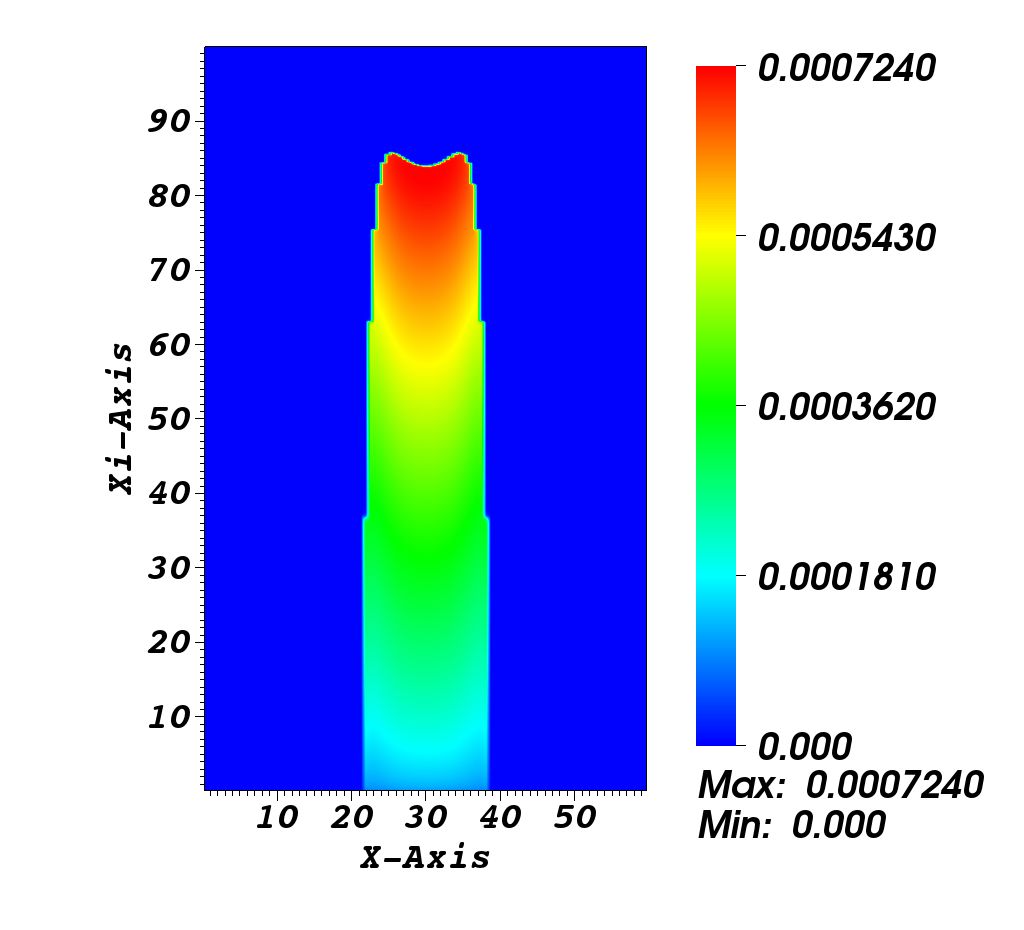}  &  \includegraphics[width=5.3cm]{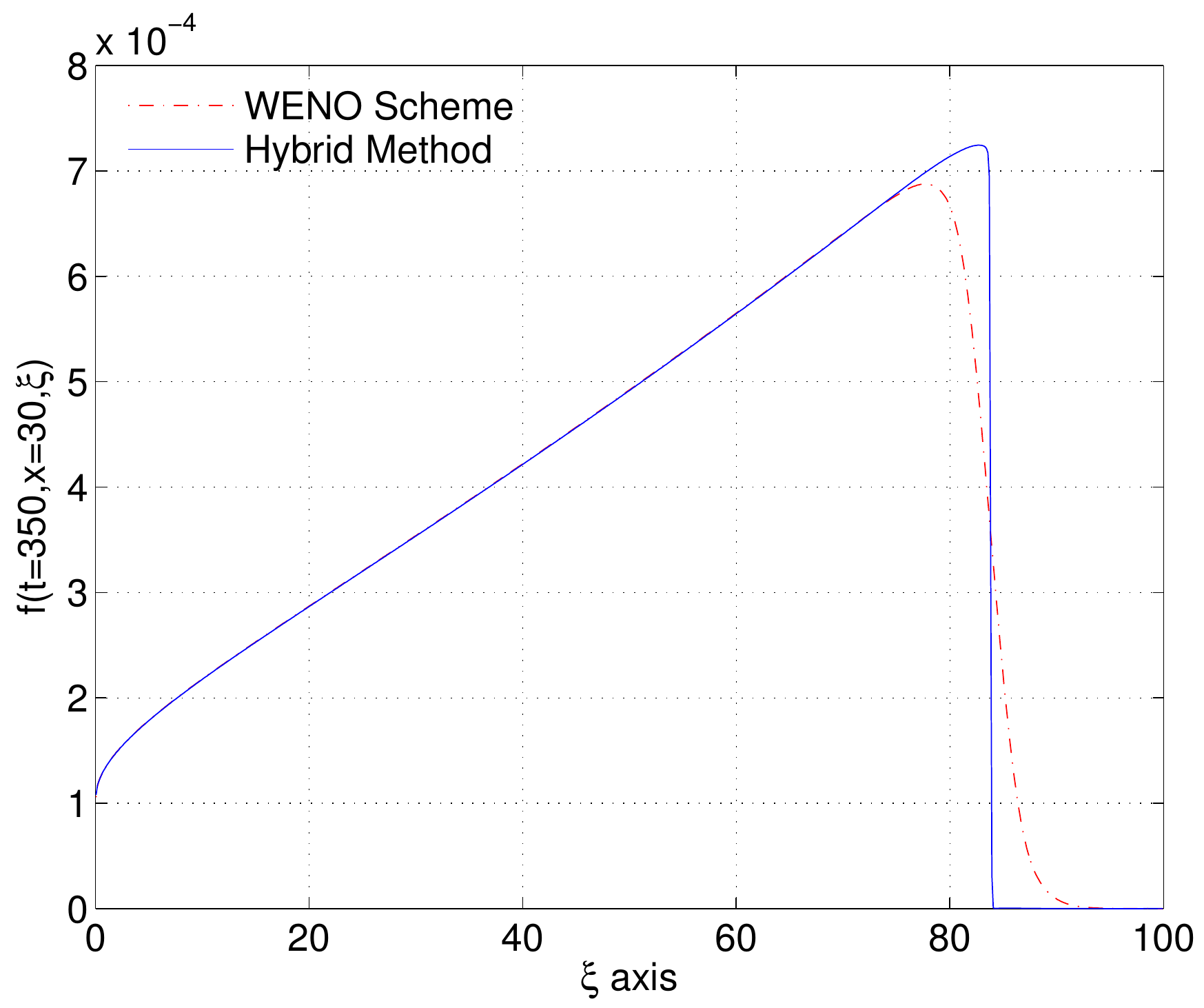}\\
  (j)  WENO Scheme                &    (k) Hybrid Method &    (l) $x=30,t=350$ \\
 \end{tabular}
\caption{\label{fig:LS2D}Polymerization/depolymerization test in non-homogeneous space: {\it Plot the size distribution function of the  equation~\eqref{eq:lif-sly1d1d} with the irregular initial data~\eqref{eq:lif-sly1d1d_init_concentration},~\eqref{eq:lif-sly1d1d_init_irreg} at different times. Mesh size is $n_x\times n_{\xi}=100\times800$, $\Delta t=0.1$, $CFL\approx0.13$.}}
 \end{center}
\end{figure}

   Secondly, we consider the irregular data~\eqref{eq:lif-sly1d1d_init_irreg}. The micro-organisms (nutrients) concentration $c(t,x)$ and the mass of the cells $\int_0^{\infty}\xi f(t,x,\xi)d\xi$ are presented in Figure~\ref{fig:LS2D_concentration_mass}. From these two quantities, we see there are no difference at different times. Moreover, the total mass for both methods preserves well in whole time evolution (see Figure~\ref{fig:growth_nonhomo_conserv}). However, we observe clearly that numerical dissipation appears with the WENO scheme, while our hybrid method preserves a sharp profile (see Figure~\ref{fig:LS2D}).\\
    
In general, having a scheme able to reconstruct accurately the profile of the solution is of paramount importance in population dynamics; indeed this profile is often used for parameter estimation and comparison with experimental measurements for prediction purposes such as prediction of the evolution of bacterial populations followed by genetic algorithms or the assessment of the average time for the balance of a product (micro-organisms or monomers in our case) in the considered mixture. That's the case in this paper because following the initial distribution, the expected solution does not have the same profile and thus does not lead to same predictions.    

\section{Conclusion and perspective}
\label{sec:conc}
\setcounter{equation}{0}

In this paper we have proposed a hybrid finite volumes scheme based on the flux convex combination between the Anti-Dissipative Method
 (ADM) \cite{bibdespres,bibGLT} and the WENO order 5 method \cite{bibJS}. The obtained numerical results show a good accuracy in 
 reconstructing the solution of transport type equation even in case of discontinuous initial data. Indeed the simulations in figure \ref{fig:growth_homo_reg} 
 show results that are as good as the ones in WENO scheme and better than ADM scheme which fail for regular initial distribution. In reverse, for irregular 
 initial distribution, the hybrid scheme show better numerical results than de ADM scheme in the sense that it advects very well the solution without ``stairs'' 
 and it shows also better results than WENO scheme which develop numerical diffusion for irregular distribution as depicted in figure \ref{fig:growth_homo}. So, when the WENO order 5 scheme fails because of numerical diffusion artifact, the hybrid method remains anti-dissipative and when ADM scheme 
 develops ``stairs" like oscillations, the hybrid scheme corrects them. This property is very suitable for long term asymptotic behavior of the solution of population dynamics, as presented in the numerical simulations for the polymerization/depolymerization type models.

 Finally, lets point out that the flux construction of our hybrid scheme depends on an empirical parameter $\alpha=\frac{3}{4}$ which is related to 
the weights in the convex combination. Nevertheless we can't explain yet why exactly $\alpha=\frac{3}{4}$ but we are working in a preprint paper 
which is focused on the theoretical and numerical explanation of this parameter $\alpha$.

\section*{Acknowledgment} 
Chang Yang is supported by National Natural Science Foundation of China (Grant No. 11401138) and Heilongjiang Postdoctoral Scientific Research Development Fund (No. LBH-Q14080).

L\'eon M. Tine is supported by the ``BQR" funding set up by Lyon 1 University for the new associate professors (DR\textbackslash
CG\textbackslash
CD\textbackslash
GG/2014-096).


\bibliographystyle{plain}

\end{document}